\documentclass[11pt,leqno]{article}

\usepackage{amsmath,amsthm}
\usepackage {latexsym}
\usepackage{amssymb}

\newcommand{\esssup}{\mbox{\rm ess sup}}

\newcommand{\bse}{\begin{equation}}

\newcommand{\be}{\begin{eqnarray}}
\newcommand{\ee}{\end{eqnarray}}

\newcommand{\supp}{\mbox{\rm supp}}
\newcommand{\spec}{{\rm sp}}

\newcommand{\half}{\frac{1}{2}}

\newcommand{\eps}{{\varepsilon}}
\newcommand{\R}{{\mathbb R}}

\newcommand{\Z}{{\mathbb Z}}

\newcommand{\Compl}{{\mathbb C}}

\newcommand{\calS}{{\mathcal S}}
\newcommand{\meas}{{\mathcal M}}

\newcommand{\sign}{\mbox{sign}}

\newcommand{\Laplace}{\triangle}

\newcommand{\Kop}{{\mathcal K}_s}
\newcommand{\kato}{{\mathcal K}}
\newcommand{\drei}{|\!|\!|}
\newcommand{\la}{\lambda}
\newcommand{\si}{\sigma}
\renewcommand{\b}{\beta}
\def\einbet{\hookrightarrow}
\def\pr{\partial}
\def\Lap{\Delta}
\def\nn{\nonumber}

\newtheorem{theorem}{Theorem}[section]
\newtheorem{lemma}[theorem]{Lemma}
\newtheorem{defi}[theorem]{Definition}
\newtheorem{cor}[theorem]{Corollary}
\newtheorem{prop}[theorem]{Proposition}
\newtheorem{proposition}[theorem]{Proposition}

\newtheorem{example}[theorem]{Example}

\theoremstyle{remark}
\newtheorem{remark}[theorem]{Remark}

\def\a{\alpha}
\def\ga{\gamma}
\def\de{\delta}
\def\il{\int\limits}

\numberwithin{equation}{section}

\begin{document}

\title{Time decay for solutions of Schr\"odinger equations with rough and time-dependent potentials.}
\author{Igor Rodnianski and Wilhelm Schlag}
\maketitle

\begin{abstract} 
\noindent In this paper we establish dispersive  estimates for solutions to
the linear Schr\"odinger equation in three dimension
\begin{equation}
\label{eq:schrabstract} 
\frac{1}{i}\partial_t \psi - \Laplace \psi + V\psi = 0,\qquad \psi(s)=f
\end{equation}
where $V(t,x)$ is a time-dependent potential that satisfies the conditions
\begin{equation} 
\nn
\sup_{t}\|V(t,\cdot)\|_{L^{\frac32}(\R^3)} + \sup_{x\in\R^3} \int_{\R^3} \int_{-\infty}^\infty \frac{|V(\hat{\tau},x)|}{|x-y|}\,d\tau\,dy < c_0.
\end{equation}
Here $c_0$ is some small constant and $V(\hat{\tau},x)$ denotes the Fourier transform
with respect to the first variable in the distributional sense. Examples of such potentials  
are of the form $V(t,x) =T(t) V_0(x)$,
where $T$ is quasiperiodic in time and $V_0\in L^{\frac 32+}(\R^3)\cap L^{\frac 32-}(\R^3)$ 
is small. We show that under these conditions~\eqref{eq:schrabstract} admits solutions 
$\psi(\cdot)\in L^{\infty}_t(L^2_x(\R^3))\cap L^2_t(L^6_x(\R^3))$ for any $f\in L^2(\R^3)$ 
satisfying the dispersive inequality
\begin{equation}
\label{eq:disabstract} 
\|\psi(t)\|_{\infty} \le C|t-s|^{-\frac32}\,\|f\|_1 \text{\ \ for all times $t,s$.} 
\end{equation}
For the case of time independent potentials $V(x)$, \eqref{eq:disabstract} remains true if 
\begin{equation}
\nn
\int_{\R^6} \frac{|V(x)|\;|V(y)|}{|x-y|^2} \, dxdy < (4\pi)^2\text{\ \ \ and\ \ \ }\|V\|_{\kato}:=\sup_{x\in\R^3}\int_{\R^3} \frac{|V(y)|}{|x-y|}\,dy<4\pi.  
\end{equation}
We also establish the dispersive estimate  with an $\eps$-loss for large energies provided $\|V\|_{\kato}+\|V\|_2<\infty$.
Finally, we prove 
Strichartz estimates for the Schr\"odinger equations with potentials that decay like~$|x|^{-2-\eps}$  in dimensions~$n\ge3$, thus solving an
open problem posed by Journ\'e, Soffer, and~Sogge.
\end{abstract}

\section{Introduction}\label{sec:intro}

\noindent It follows from the explicit expression for the kernel of $e^{-it\Laplace}$ that 
the free Schr\"odinger evolution in $\R^n$, $n\ge1$,  satisfies the dispersive inequality
\begin{equation}
\label{eq:freedisperse} \|e^{-it\Laplace} f\|_{L^\infty_x} \le C t^{-\frac{n}{2}}\|f\|_{L^1_x}.
\end{equation}
Closely related are the classical Strichartz estimate~\cite{strich}
\[ \|e^{-it\Laplace} f\|_{L^{2+\frac{4}{n}}(\R^{n+1})} \le C\|f\|_{L^2(\R^n)}\]
or more generally
\begin{equation}
\label{eq:genstr}
\|e^{-it\Laplace} f\|_{L_t^p L_x^q(\R^n)} \le C\|f\|_{L^2(\R^n)}
\end{equation}
for any $\frac{n}{q}+\frac{2}{p}=\frac{n}{2}$, $2\le p\le\infty$. The case $p=\infty$, $q=2$ 
is the energy estimate (in fact $\|e^{-itH}f\|_2=\|f\|_2$), whereas the range $2<p<\infty$ can be obtained
from the case $p=2$ and~\eqref{eq:freedisperse} by means of a well-known argument (see for example~\cite{KT}). 
The endpoint $p=2, q=\frac{2n}{n-2}$ result, which in fact fails in dimension $n=2$, is more difficult and was recently settled for $n\ge 3$ 
by Keel and Tao~\cite{KT}. 

  The question whether these bounds also hold for more general Schr\"odinger equations has been considered by various authors. From a physical perspective it is of course natural to consider the case of $e^{itH}$ with $H=-\Laplace+V$. For the purposes of the present discussion we assume that the potential~$V$ is real and has enough regularity to ensure that~$H$ is a self-adjoint operator on~$L^2(\R^n)$, see Simon's review~\cite{barry} for explicit conditions on~$V$. 
One obstacle to having decay in time for $e^{itH}$ are eigenvalues of the operator $H=-\Laplace+V$ and a result as in~\eqref{eq:freedisperse} and~\eqref{eq:genstr} therefore requires that~$f$ be orthogonal to any eigenfunction of~$H$. In fact, Journ\'e, Soffer, and~Sogge~\cite{JSS} have shown that, with $P_{c}$ being the projection onto the continuous subspace of~$L^2(\R^n)$ with respect to~$H$,
\begin{equation}
\label{eq:JSS} 
\|e^{it(-\Laplace +V)} P_{c} f\|_{\infty} \le C\,t^{-\frac{n}{2}}\|f\|_{L^1(\R^n)}
\end{equation}
for all dimensions $n\ge3$ provided that zero is neither an eigenvalue nor a resonance of~$H$. In addition, they need to assume that,  roughly speaking, $|V(x)|\lesssim (1+|x|)^{-n-4}$ and~$\hat{V}\in L^1(\R^n)$.   
Recall that a resonance is a distributional solution of $H\psi=0$ so that $\psi\not\in L^2$ but~$(1+|x|^2)^{-\frac{\sigma}{2}}\psi(x)\in L^2$ for any~$\sigma>\half$, see~\cite{JenKat}. 
It is well-known that under the assumptions on~$V$ used in~\cite{JSS} the spectrum $\sigma(H)$ satisfies
\[ \sigma(H)=[0,\infty)\cup \{\lambda_j\:|\: j=1,\ldots,N\} \]
where $[0,\infty)=\sigma_{a.c.}(H)$ and $\lambda_N< \lambda_{N-1}< \ldots <\lambda_1\le 0$ are a discrete and finite set of eigenvalues of finite multiplicity. Indeed, since~$V$ is bounded and decays at infinity Weyl's criterion (Theorem~XIII.14 in~\cite{RS4}) implies that $\sigma_{ess}(H)=\sigma_{ess}(-\Laplace)=[0,\infty)$, whereas the finiteness follows from the Cwikel-Lieb-Rosebljum bound. Furthermore, since $V$ is bounded and decays faster than $|x|^{-1}$ at infinity it follows from Kato's theorem (Theorem~XIII.58 in~\cite{RS4}) that there are no positive eigenvalues of~$H$. Finally, since any $V$ as in~\cite{JSS} is an Agmon potential, $\sigma_{sing}(H)=\emptyset$ by the Agmon-Kato-Kuroda theorem (Theorem~XIII.33 in~\cite{RS4}). 

The work by Journ\'e, Soffer, and~Sogge was preceded by related results of Rauch~\cite{rauch}, Jensen, Kato~\cite{JenKat}, and Jensen~\cite{Jensen1},\cite{Jensen2}. The fact that one cannot have~$t^{-\frac32}$ decay in the presence of a  resonance at zero energy was observed by these authors. Moreover, the small energy asymptotic expansions of the resolvent developed in~\cite{JenKat}, \cite{Jensen1}, \cite{Jensen2} are used in~\cite{JSS}. However, the actual time decay estimates obtained by Rauch, Jensen, and Kato are formulated in terms of weighted $L^2$-spaces rather than in the much stronger $L^1\to L^\infty$ sense of Journ\'e, Soffer, and~Sogge. 
The appearance of weighted $L^2$ spaces is natural in view of the so called limiting absorption principle. 
This refers to boundedness of the resolvents $(-\Laplace - \lambda \pm i0)^{-1}$ for $\lambda>0$ on certain weighted $L^2$ spaces as proved by Agmon~\cite{agmon} and Kuroda~\cite{kur1}, \cite{kur2}. 
It is also with respect to these weighted norms that the asymptotic expansions of the resolvents $(H-z)^{-1}$ as $z\to0$ with 
$\Im(z)\ge0$, $\Re(z)>0$ in~\cite{JenKat},~\cite{Jensen1},\cite{Jensen2} hold. Jensen and Kato need to assume that~$|V(x)|\lesssim (1+|x|)^{-\beta}$ 
for certain $\beta>1$ (most of their results require $\beta>3$). For a more detailed discussion of the limiting absorption principle see
our Strichartz estimates in Section~\ref{sec:strichartz}.

Another approach to decay estimates for $e^{it(-\Laplace+V)}$ 
was taken by Yajima~\cite{Y1}, \cite{Y2}, and Artbazar and Yajima~\cite{Y4}, who relied on scattering theory. 
Recall that if the so called wave-operator
\[ W = s-\lim_{t\to\infty} e^{-it(-\Laplace + V)} e^{-it\Laplace} \]
exists, where the limit is understood in the strong $L^2$ sense, then it is an isometry that intertwines the evolutions, i.e.,
\[ W e^{-it\Laplace} = e^{it(-\Laplace + V)} W \text{\ \ for all times $t$}.\]
In \cite{Y1} Yajima proved that the wave operators $W$ are bounded from $L^p(\R^n)\to L^p(\R^n)$ with $n\ge3$  
for $1\le p\le\infty$ provided $V$ has a certain explicit amount of decay, and provided zero is neither an eigenvalue nor a resonance. Since $WW^*=P_{a.c.}$, he concludes from the free dispersive estimate~\eqref{eq:freedisperse} that 
\[ 
\|e^{it(-\Laplace + V)}P_{a.c.}\|_{L^1\to L^\infty} = \|W e^{-it\Laplace} W^*\|_{L^1\to L^\infty} \le C\,t^{-\frac{n}{2}}
\]
under the usual assumption on the zero energy but imposing weaker conditions on $V$ than~\cite{JSS}. 
Moreover, \cite{Y2} contains the first results on dispersive and Strichartz estimates for $e^{it(-\Laplace+V)}$ in two dimensions. 

The relatively strong decay and regularity assumptions that appear in all aforementioned works are by far sufficient to ensure scattering, i.e, the existence of wave operators on $L^2$, even though Yajima was the first to exploit this link explicitly in the context of dispersive estimates. The connection with scattering is of course natural, as the decay of~$V$ (and possibly that of derivatives of~$V$) at infinity allows one to reduce matters to the free equation by methods that are to a large extent perturbative. 

On the other hand, the existence of scattering (in the traditional $L^2$ sense) is known for potentials that are small in some global sense, 
but without any explicit rate of decay. Indeed, it is a classical result of Kato~\cite{kato} that under the sole assumption that the real potential $V$ satisfies
\begin{equation}
\label{eq:kleinrol} 
\int_{\R^6} \frac{|V(x)|\;|V(y)|}{|x-y|^2} \, dxdy < (4\pi)^2 
\end{equation}
the operator $H=-\Laplace+V$ on $\R^3$ is self-adjoint and unitarily equivalent to~$-\Laplace$ via the wave operators. The left--hand side of~\eqref{eq:kleinrol} is usually referred to as the Rollnik norm, see~\cite{simon1}.  Observe that \eqref{eq:kleinrol} roughly corresponds to the potential decaying 
at infinity as $|x|^{-2-\eps}$.

The appearance of the Rollnik norm in the context of small potentials is natural from several perspectives, one of which is scaling. 
The Rollnik norm is invariant under the scaling  $R^2 V(Rx)$ forced by 
the Schr\"odinger operator~$H$  onto the potential~$V$. 
It is well--known that the Rollnik norm defines a class of potentials that is slightly wider than $L^{\frac32}(\R^3)$, which is also scaling invariant. Another natural occurrence of a scaling invariant condition arises in connection with bounds on the number of negative eigenstates. Indeed, in dimension $n$
it is precisely the scaling invariant $L^{\frac n2}$ norm of the negative part of the potential
that governs the number of negative eigenvalues of $-\Delta +V$ via the Cwikel-Lieb-Rosebljum bound. 

We show in this work that dispersive estimates lead naturally to what we call the ``global Kato norm'' of the potential. Recall that the Kato norm of~$V$ is defined to be 
$$
\sup_{x\in\R^3} \int_{|x-y|\le 1} \frac{|V(y)|}{|x-y|}\, dy,
$$
whereas the scaling invariant analogue is given by~\eqref{eq:kleinkato} below. 
The Kato norm, or more precisely the closely related Kato class, arise in the study of self-adjoint extensions
of~$H$, as well as in the study of the properties of the heat semigroup $e^{-tH}$, see~\cite{AS}, 
\cite{barry}, and Section~\ref{sec:epsloss} below.

One of the goals of our paper is to bridge the gap between the ``classical'' perturbation results of spectral theory that involve Rollnik and Kato classes
of potentials (or other scaling invariant classes) and the results concernning the dispersive properties of the time-dependent Schr\"odinger
equation.

In our first result, see Theorem~\ref{thm:high} below, we show that the dispersive estimates are stable
under perturbations by small potentials that belong to the intersection of the Rollnik and the global Kato classes. 

\begin{theorem} 
\label{thm:introthm1}
Suppose $V$ is real and satisfies~\eqref{eq:kleinrol}. Suppose in addition that
\begin{equation}
\label{eq:kleinkato} 
\sup_{x\in\R^3} \int_{\R^3} \frac{|V(y)|}{|x-y|}\, dy < 4\pi.
\end{equation}
Then one has the estimate
\[ \| e^{it(-\Laplace + V)}\|_{L^1\to L^\infty} \lesssim t^{-\frac32} \]
for all $t>0$. 
\end{theorem}

The proof relies  on a Born series expansion for the resolvent with a subsequent estimate of an arising oscillatory integral. 
The convergence of the resulting geometric series is guaranteed by \eqref{eq:kleinkato}. See Section~\ref{sec:stat} for details.

The main focus of this paper is on the dispersive properties of solutions of the Schr\"odinger equation \eqref{eq:schrabstract} with time dependent potentials, 
see Sections~\ref{sec:time1}--\ref{sec:time4}. It appears that not much is known on the long time behavior of
solutions to Schr\"odinger equations with time dependent potentials. See, however, Bourgain~\cite{master2}, \cite{master2} on the issue of slow growth of higher Sobolev norms in the space-periodic setting. 
In this paper we establish  dispersive and Strichartz estimates for a class of scaling invariant small potentials on~$\R^3$. 

\begin{theorem}
Let $V(t,x)$ be a real-valued measurable function on~$\R^4$ such that 
\begin{equation}
\sup_{t}\|V(t,\cdot)\|_{L^{\frac32}(\R^3)} + \sup_{y\in\R^3} \int_{\R^3}\int \frac{|V(\hat{\tau},x)|}{|x-y|}\,d\tau\,dx < c_0
\end{equation}
for some small constant $c_0>0$. Here $V(\hat{\tau},x)$ denotes the Fourier transform in the first variable, and if $V(\hat{\tau},x)$ happens to be a measure then the $L^1$--norm in~$\tau$ gets replaced with the norm in the sense of measures. Then for every initial time~$s$ and every $\psi_s\in L^2(\R^3)$ the equation
\be
&&\frac{1}{i}\pr_t \psi - \Lap \psi + V(t,x) \psi=0,\label{eq:schr_t}\\
&&\psi|_{t=s} (x) =\psi_s(x)\nonumber
\ee
admits a (weak) solution $\psi(t,\cdot)= U(t,s)\psi_s$ (via the Duhamel formula). The propagator $U(\cdot,s)$ satisfies  $U(\cdot,s):L^2(\R^3)\to L^{\infty}_t(L^2_x(\R^3))\cap L^2_t(L^6_x(\R^3))$,  $t\mapsto \psi(t,\cdot)$ is weakly continuous as a map into~$L^2(\R^3)$, and $\|U(t,s)\psi_s\|_2\le \|\psi_s\|_2$. Finally, $U(t,s)$ satisfies the dispersive inequality
\begin{equation}
\label{eq:disintro} 
\|U(t,s)\psi_s\|_{L^\infty} \le C|t-s|^{-\frac32}\,\|\psi_s\|_{L^1} \text{\ \ for all times $t,s$ and any $\psi_s\in L^1$.} 
\end{equation}
\label{thm:introthm2}
\end{theorem}

\noindent Examples of potentials to which the theorem applies are $V(t,x)=\cos(t) V_0(x)$ where $\|V_0\|_{L^{\frac32}(\R^3)}<c_0$, and for which~\eqref{eq:kleinkato} holds. More generally, one can take potentials that are quasi-periodic in time, such as $V(t,x)= \phi(t) V_0(x)$ with
\[ \phi(t) = \sum_{\nu \in Z^d} c_\nu e^{2\pi i t \omega\cdot \nu} \]
and $\sum_{\nu\in \Z^d} |c_\nu| < \infty$, $\omega\in [0,1)$ arbitrary. 

\noindent Note that Theorem~\ref{thm:introthm2} also applies to time independent potentials $V_0(x)$  via $V(t,x):=V_0(x)$. 
Clearly, in that case the conditions become
\[ \|V_0\|_{L^{\frac32}(\R^3)} + \sup_{x\in\R^3} \int_{\R^3} \frac{|V(y)|}{|x-y|}\, dy < c_0.\]
Since by fractional integration
\[ \int_{\R^6} \frac{|V(x)|\;|V(y)|}{|x-y|^2} \, dxdy \le C\|V\|^2_{L^{\frac32}(\R^3)},\]
it follows that these conditions are strictly stronger than those in Theorem~\ref{thm:introthm1}. 

\noindent Whereas our main emphasis is of course on the decay estimate~\eqref{eq:disintro}, it appears that
even the easier question of solvability of equation~\eqref{eq:schr_t} for rough potentials that do not decay in time had not been addressed before, at least under the conditions of Theorem~\ref{thm:introthm2}. 
Yajima~\cite{Y3} considered the problem of existence of solutions to the Schr\"odinger equation with time-dependent potentials. In his paper he proves the existence of the strongly continuous semigroup $U(t,s)$ 
on $L^2(\R^n)$ provided that the potential satisfies $V\in L^q_t L^p_x$ for 
$0\le \frac 1q<1-\frac n{2p}$. Notice that in our case $q=\infty$, $p=\frac n2$, which corresponds
to the endpoint of this condition not covered in~\cite{Y3}.  
We use the endpoint Strichartz estimate~\cite{KT} for the {\em free} problem for that purpose, which automatically yields the endpoint Strichartz estimate in the context of Theorem~\ref{thm:introthm2}. 

\noindent For time-dependent potentials the analogue of  Kato's scattering result~\cite{kato} was proved by Howland~\cite{Howl1}. More precisely,  under the condition that for a sufficiently large time $t_0>0$,
$V(t,x)\le V_0(x)$ for some {\em time independent} potential $V_0(x)$ obeying the small Rollnik condition \eqref{eq:kleinrol},
there exist a unitary wave operator $W$ intertwining $U(t,s)$ and $e^{it(-\Delta)}$. 
In case $V(t,x)$ does decay in time (in the sense of a small amount of integrability), wave operators were constructed by Howland~\cite{Howl2}  and Davies~\cite{Da}. In contrast to Theorem~\ref{thm:introthm2} they do not require smallness (the latter being replaced by time decay of the potential) and they also obtain strong continuity of the evolution.

\noindent 
One of the difficulties in this case is the absence of the connection between
the semigroup generated by the Schr\"odinger equation and the spectral properties of the operator $-\Delta +V$. Recall that for time independent potentials~$V$,
$$
e^{itH} f = \int e^{it\lambda} dE(\lambda)f
$$
where $dE(\lambda)$ is the spectral measure of the operator $-\Delta + V$. This is no longer available for time-dependent potentials. 

\noindent The proof of \eqref{eq:disintro} is similar to that of Theorem~\ref{thm:introthm1} but much more involved. Since we can no longer rely on the spectral theorem, resolvents, and Born series to construct the evolution of~\eqref{eq:schr_t}, we use the Duhamel formula instead (we note in passing that the Fourier transform in the spectral parameter establishes an equivalence between the representation of the evolution in terms of a Born series and an infinite expansion of the solution by means of Duhamel's principle). One of the novelties in our paper is the formula representing the time evolution of the Schr\"odinger
equation with a time-dependent potential as an infinite series of oscillatory integrals involving the resolvents of the {\em free} problem.
Most of the work in the proof of Theorem~\ref{thm:introthm2} is devoted to estimating these oscillatory integrals, whose phases typically have a critical point with degeneracies of the third order. See Sections~\ref{sec:time1}--\ref{sec:time4} for details.

\noindent Two sections of this paper are devoted to  time independent potentials without any restrictions on their sizes. In Section~\ref{sec:epsloss} we prove the following result. As before, $H=-\Laplace+V$ and~$P_{a.c.}$ refers to the projection onto the absolutely continuous subspace of~$L^2$ relative to~$H$. 

\begin{prop} 
\label{prop:introprop}
Let  
\[ \drei V \drei :=\|V\|_2 + \sup_{x\in\R^3} \int_{\R^3} \frac{|V(y)|}{|x-y|}\,dy < \infty.\]
Then for every $\eps>0$ there exists some positive $\lambda_0=\lambda_0(\drei V\drei,\eps)$ so that
\begin{equation}
\label{eq:introhigh}
\| e^{itH} \chi(H/\lambda_0) P_{a.c.} \|_{L^1\to L^\infty} \le Ct^{-\frac32+\eps}
\end{equation}
for all $t>0$.
\end{prop}

\noindent The proof is again perturbative. For the case of large energies, {\em and for those only}, 
the required smallness is provided by the following estimate on the resolvents, which can be viewed as some instance of the limited absorption principle:
\begin{equation}
\label{eq:introST}
\|(-\Laplace-\lambda+i0)^{-1}f\|_{L^4(\R^3)} \le C\lambda^{-\frac14}\|f\|_{L^{\frac43}(\R^3)}.
\end{equation}
The proof of \eqref{eq:introST} is an immediate consequence of the Stein-Tomas theorem~\cite{stein}.
The appearance of the Stein-Tomas theorem in this context is most natural, as the resolvent 
$(-\Laplace-\lambda+i0)^{-1}$ of the free problem is closely related to the restriction of the 
Fourier transform to the sphere $|x|=\sqrt{\lambda}$ for $\lambda>0$. In contrast to~\eqref{eq:introST},
which heavily relies on the nonvanishing Gaussian curvature of the sphere,  the
classical limiting absorption principle of Agmon and Kuroda~\cite{agmon}, \cite{kur2}, and \cite{kur1}
only uses the most elementary restriction property of the Fourier transform to arbitrary surfaces which leads
to a loss of $\half+\eps$ derivatives in~$L^2$ (on the physical side this translates into the weights $|x|^{\half+\eps}$ in~$L^2$ that appear in~\cite{agmon}, \cite{JenKat} etc.). For further details of the proof of Proposition~\ref{prop:introprop} we refer the reader to Section~\ref{sec:epsloss}. 

\noindent It is common knowledge that the case of large energies should be the most
accessible one. From the perspective of scattering the intuition is that particles with high energies
will escape the scatterer and thus lead to extended states (absolutely continuous spectrum) whereas
particles with smaller energies can be trapped and create bound states (pure point spectrum). 
It is of course a most interesting problem to extend Proposition~\ref{prop:introprop} to small energies under
similar conditions. Recall that~\cite{JSS} and particularly~\cite{Y1} have accomplished exactly that, but 
under conditions on~$V$ that are by far stronger than those in Proposition~\ref{prop:introprop}.

\noindent We also address the question of Strichartz estimates for $e^{it(-\Laplace + V)}$ in dimensions greater or equal than three. Traditionally the mixed norm Strichartz estimates~\eqref{eq:genstr} are shown to be a consequence of the dispersive estimates. In fact,  in~\cite{JSS}, Journ\'e,  Soffer, and Sogge 
establish the $L^1\to L^\infty$ dispersive bound and therefore also Strichartz estimates under strong decay and regularity assumptions on~$V$, see~\eqref{eq:JSS}. However, they conjecture that 
 Strichartz estimates hold for potentials that decay only faster than~$(1+|x|)^{-2}$. 
In this paper we prove this conjecture assuming only this rate of decay. In particular, we do not require
any regularity. More precisely, the following theorem holds.

\begin{theorem}
\label{thm:introthm3}
Suppose that for some $\eps>0$ one has $|V(x)|\lesssim (1+|x|)^{-2-\eps}$ for all $x\in\R^n$ with $n\ge3$. Then
\[ \|e^{itH} P_{c}f\|_{L^q_t L^r_x(\R^n)} \lesssim \|f\|_{L^2_x(\R^n)}  
\qquad \forall (q,r,n),\quad
\frac 2q=n(\frac 12 -\frac 1r)
\]
provided the zero energy is neither an eigenvalue nor a resonance of the operator $H=-\Lap +V$.
Here $P_c$ denotes the spectral projection onto the continuous states.
\end{theorem}

\noindent The decay condition $|V(x)|\lesssim (1+|x|)^{-2-\eps}$ is very natural from the perspective
of Kato's smoothing theory~\cite{kato}. In contrast to~\cite{JSS} we prove the Strichartz estimates directly, i.e., without relying on dispersive estimates. In fact, we do not know if the $L^1\to L^\infty$ estimates hold
under the conditions of Theorem~\ref{thm:introthm3}. 
It is known that (local in time) Strichartz estimates can hold even if the $L^1\to L^\infty$ dispersive property fails, see Bourgain~\cite{master3} for the case of the torus, Staffilani, Tataru~\cite{ST} for variable coefficients, and Burq, Gerard, Tzvetkov~\cite{BGT} for the case of equations on Riemannian manifolds. 

This paper is organized as follows: Section~\ref{sec:stat} to~\ref{sec:strichartz} deal with time independent potentials. Section~\ref{sec:stat} establishes dispersive estimates for small Rollnik potentials in~$\R^3$. Section~\ref{sec:epsloss} considers the high energy case for low regularity potentials, and in 
Section~\ref{sec:strichartz} we establish mixed norm Strichartz estimates for potentials that decay like~$(1+|x|)^{-2-\eps}$. 
The remaining sections~\ref{sec:time1}-\ref{sec:time4} are devoted to small time-dependent potentials. In Section~\ref{sec:time1} we show that solutions exist for potentials that do not necessarily decay in time by means of the Keel-Tao~\cite{KT} endpoint. We then proceed to represent the solution by means of an infinite Duhamel expansion and we derive a formula for each term in the Duhamel series. The most technical part are Sections~\ref{sec:time2} and~\ref{sec:time3} that provide the necessary bounds on the oscillatory integrals that arise in this context. We combine all the pieces in the final Section~\ref{sec:time4}.

\underline{Acknowledgments:} The authors thank Alexander Pushnitski for valuable discussions on the Agmon-Kato-Kuroda theory, 
Thomas Spencer for his interest in the problem
of time dependent potentials, as well as Elias Stein for a discussion on Section~\ref{sec:time2}.
The first author was supported by an NSF grant DMS-0107791.
The second author was supported by an NSF grant DMS-0070538 and a Sloan fellowship. 
Part of this work was done while he was a member at the Institute for Advanced Study, Princeton and
partially supported by an NSF grant DMS-9729992.

\section{Small time independent potentials in $\R^3$}\label{sec:stat}

\noindent The purpose of this section is to prove the $L^1(\R^3)\to L^\infty(\R^3)$  dispersive inequality for $e^{itH}$ 
where $H=-\Laplace+V$ in $\R^3$. The following definition states the properties of the real potential~$V$ 
that we will need.

\begin{defi}
\label{def:V}
We require that both 
\be
\label{eq:rollnik} \|V\|_{R}^2:=\int_{\R^3\times \R^3} \frac{|V(x)|\,|V(y)|}{|x-y|^2}\,dx\,dy &<& (4\pi)^2  
\text{\ \ \ and}\\
\label{eq:kato} \|V\|_{\kato}:=\sup_{x\in\R^3} \int_{\R^3} \frac{|V(y)|}{|x-y|}\;dy &<& 4\pi.
\ee
\end{defi}

\noindent The norm $\|\cdot\|_R$ on the left-hand side of~\eqref{eq:rollnik} is usually referred to as the {\em Rollnik norm}. Kato~\cite{kato} showed that under the condition~\eqref{eq:rollnik} the operator~$H$ admits a self-adjoint extension which is unitarily equivalent to $H_0=-\Laplace$. In particular, the spectrum of~$H$ is purely absolutely continuous. 
Many properties of the Rollnik norm, which can be seen to be majorized by the norm of~$L^{\frac32}(\R^3)$ via fractional integration,  
can be found in Simon's monograph~\cite{simon1}.  
The norm $\|\cdot\|_{\kato}$ in~\eqref{eq:kato} is closely related to the well-known {\em Kato norm}, 
see Aizenman and Simon~\cite{AS}, \cite{barry} and we refer to it as the {\em global Kato norm}. 

\noindent The main result in this section is Theorem~\ref{thm:high}. 
The proof splits into several lemmas, the first of which presents some well-known properties of the resolvents
$R_V(z)=(-\Laplace+V-z)^{-1}$ under the condition~\eqref{eq:rollnik}. We begin by recalling 
that a potential with finite (but not necessarily small) Rollnik norm is {\em Kato smoothing}, i.e.,
\begin{equation}
\label{eq:katosm} 
\sup_{\eps>0}
\|\,|V|^{\frac 12} R_0(\la\pm i\eps) f\|_{L^2_\la L^2_x}\le C \|f\|_{L^2},\qquad
\sup_{\eps>0}\|\,\,R_0(\la\pm i\eps) \,|V|^{\frac 12} f\|_{L^2_\la L^2_x}\le C \|f\|_{L^2}
\end{equation}
for any $f\in L^2(\R^3)$ and with $R_0(z)=(-\Laplace-z)^{-1}$. 
This implies, in particular, that ${\mathcal D}(|V|^{\half})\supset H^2$.
The Rollnik norm arises in this context as a majorant for the Hilbert-Schmidt norm $\|\cdot\|_{HS}$ of the operators
\begin{equation}
\label{eq:fundop} 
K(\lambda\pm i\eps):= |V|^{\half}R_0(\lambda\pm i\eps)|V|^{\half}.
\end{equation}
Indeed, it is well-known that the resolvent $R_0(z)$ for $\Im z\ge0$ has the kernel
\begin{equation}
\label{eq:R0kernel} 
R_0(z)(x,y) = \frac{\exp(i\sqrt{z}|x-y|)}{4\pi|x-y|}
\end{equation}
where $\Im(\sqrt{z})\ge0$. Thus
\begin{equation}
\label{eq:Klam} 
\|K(z)\|_{L^2\to L^2}\le \|K(z)\|_{HS}\le (4\pi)^{-1}\|V\|_R,
\end{equation}
for every $z\in\Compl$ with $\Im z\ge0$.
This allows one to check immediately that $S_z:=|V|^{\half}R_0(z):L^2\to L^2$ for every $z\in{\mathbb C}\setminus\R$. Indeed, by the resolvent identity,
\[ S_zS_z^* = \frac{1}{-2i\Im z}\bigl[ |V|^{\half}R_0(z)|V|^{\half}-|V|^{\half}R_0(\bar{z})|V|^{\half}\bigr].
\] 
In view of \eqref{eq:Klam} therefore
\begin{equation}
\label{eq:Sz} 
\|S_z\|^2=\| S_zS_z^* \| \lesssim \frac{1}{|\Im z|} \|K(z)\| \lesssim \frac{1}{|\Im z|} \|V\|_R,
\end{equation}
as desired. One of the main observations of Kato~\cite{kato} was the relation between 
this pointwise condition in~$z=\lambda\pm i\eps$ and the $L^2_\lambda$ boundedness that appears in~\eqref{eq:katosm}.  We present a short proof of this fact for the sake of completeness. Although it is standard, 
 the following argument is somewhat different from the usual one which can be found in basic references like Kato~\cite{kato} and Reed, Simon~\cite{RS4}. Denote $T_\eps:= |V|^{\frac 12} R_0(\la+ i\eps)$ for $\eps>0$. Truncating the large
values of $V$ and then passing to the limit we may assume that $V$ is bounded. Then $T_\eps:L^2\to L^2_\lambda L^2_x$ for every~$\eps>0$ and one checks that
\[ T^*F = \int R_0(\la-i\eps)|V|^{\frac 12} F(\la)\,d\la\]
for every $F\in L^2_\lambda(L^2_x)$. Thus
\be 
T_\eps T_\eps^*F &=& \int |V|^{\frac 12} R_0(\la+ i\eps) R_0(\mu - i\eps) |V|^{\frac 12}F(\mu)\,d\mu \nn \\
 &=& -\int |V|^{\frac 12} \frac{R_0(\la+ i\eps)- R_0(\mu - i\eps)}{\lambda-\mu+2i\eps} |V|^{\frac 12} F(\mu)\,d\mu \label{eq:secresol} \\
 &=& -\int \frac{K(\la+ i\eps)F(\mu)}{\lambda-\mu+2i\eps} \,d\mu + \int \frac{K(\mu- i\eps)F(\mu)}{\lambda-\mu+2i\eps} \,d\mu, \label{eq:hilbklam}
\ee
where we used the resolvent identity to pass to \eqref{eq:secresol}. 
By the $L^2$ boundedness of the (vector valued) Hilbert transform,
\[
\sup_{\eps>0} \Bigl\| \int \frac{F(\mu)}{\lambda+i\eps-\mu}\, d\mu \Bigr\|_{L^2_\lambda L_x^2} \lesssim 
 \|F\|_{L^2_\lambda(L_x^2)}.
\]
Using this bound and \eqref{eq:Klam} in \eqref{eq:hilbklam} yields
\[ \sup_{\eps>0}\|T_\eps T_\eps^*F\|_{L^2_x} \lesssim \|F\|_{L^2_\lambda L_x^2}\|V\|_R\]
which implies \eqref{eq:katosm} with a constant of the form $C\|V\|_R^{\half}$.
 
\begin{lemma}
\label{lem:resolvent}
Let $\|V\|_R<4\pi$ as in Definition~\ref{def:V}. Then for all $f, g\in L^2(\R^3)$ 
\begin{equation}
\label{eq:Resseries}
\langle R_V(\la\pm i\eps) f,\,g\rangle - \langle R_0(\la\pm i\eps) f,\,g\rangle = \sum_{\ell=1}^\infty (-1)^\ell \langle R_0(\la\pm i\eps)(VR_0(\la\pm i\eps))^\ell f,g\rangle
\end{equation}
where the right-hand side of \eqref{eq:Resseries} is an absolutely convergent series in 
the norm of~$L^1(d\la)$ uniformly in $\eps>0$. Furthermore, if $\|V-V_m\|_R\to 0$ as $m\to\infty$, then
\begin{equation}
\label{eq:contR} 
\sup_{\eps>0}\int\Bigl|\langle R_{V_m}(\la\pm i\eps) f,\,g\rangle - \langle R_V(\la\pm i\eps) f,\,g\rangle\Bigr| \,d\lambda \to 0
\end{equation} 
as $m\to\infty$.
\end{lemma}
\begin{proof}
We start from the resolvent identity
\begin{equation}
\label{eq:secres} R_V(z)-R_0(z)= -R_0(z) V R_V(z) = - R_V(z)V R_0(z)
\end{equation}
which holds in the sense of bounded operators on $L^2$ for any $\Im z\not=0$, see~\eqref{eq:Sz}.  
It is a standard fact, see~\cite{kato},  that the Kato smoothing property~\eqref{eq:katosm} remains 
valid with $R_V$ instead of $R_0$ provided that $\|V\|_R<4\pi$. Indeed, multiplying~\eqref{eq:secres}
by~$|V|^{\half}$ leads to
\begin{equation}
\label{eq:Q} 
(1+Q(z))A R_V(z) = A R_0(z)
\end{equation} 
where $Q(z):=A R_0(z) B$, $A=|V|^{\frac 12}$, and $B=|V|^{\frac 12} \sign V$. 
In view of~\eqref{eq:Klam} one has 
\begin{equation}
\label{eq:1plusQ}
\sup_{\Im z\not=0}\|Q(z)\|_{L^2\to L^2}=:\rho<1\text{\ \  so that\ \ } \sup_{\Im z\not=0}\|(1+Q(z))^{-1}\|_{L^2\to L^2}\le (1-\rho)^{-1}.
\end{equation}
In conjunction with~\eqref{eq:Q} and~\eqref{eq:katosm} this implies that
\begin{equation}
\label{eq:Vkatosm} 
\sup_{\eps>0}
\|\,A R_V(\la\pm i\eps) f\|_{L^2_\la L^2_x}\le C \|f\|_{L^2},\qquad
\sup_{\eps>0}\|\,\,R_V(\la\pm i\eps) \,B f\|_{L^2_\la L^2_x}\le C \|f\|_{L^2}
\end{equation}
for any $f\in L^2$. Fix $f,g \in L^2$. Iterating \eqref{eq:secres} leads to
\begin{equation}
\label{eq:neuman} 
\langle R_V(\la\pm i\eps) f,\,g\rangle = \sum_{\ell=0}^N (-1)^\ell \langle R_0(\la\pm i\eps)(VR_0(\la\pm i\eps))^\ell f,g\rangle + (-1)^{N+1} \langle R_V(\la\pm i\eps)(VR_0(\la\pm i\eps))^{N+1} f,g\rangle 
\end{equation}
for any positive integer $N$. By \eqref{eq:Vkatosm} the  error term is 
\[ 
\langle R_V(\la\pm i\eps) B (AR_0(\la\pm i\eps)B)^N AR_0(\la\pm i\eps)f,g\rangle  = \langle  R_V(\la\pm i\eps) B Q(\la\pm i\eps)^N AR_0(\la\pm i\eps)f,g\rangle
\] 
and thus has $L^1(d\lambda)$ norm bounded by $C\,\rho^N$, see~\eqref{eq:Vkatosm} and~\eqref{eq:1plusQ}.
Similarly, each of the terms in the sum for $1\le \ell\le N$ has $L^1(d\lambda)$ norm at 
most $C\,\rho^{\ell-1}$. Thus~\eqref{eq:Resseries} holds for any~$V$ which satisfies~\eqref{eq:rollnik}. 
If $m$ is sufficiently large, then the series expansion~\eqref{eq:Resseries} holds for both $V$ and~$V_m$. 
Subtracting these series termwise and invoking the previous bounds yields that the left-hand side of~\eqref{eq:contR} is bounded by
\[ \sum_{\ell=1}^\infty C\ell \rho^{\ell-1}\,\|V-V_m\|_R \le C(1-\rho)^{-2}\,\|V-V_m\|_R,\]
and the lemma follows.
\end{proof}

\noindent The following technical corollary deals with the case $\eps=0$ in Lemma~\ref{lem:resolvent}. 
We state it in the form in which it is used later on. In particular, we did not strive for the 
greatest generality. Below $C_b^0(\R)$ refers to the bounded continuous functions on~$\R$ 
with the supremum norm.

\begin{cor}
\label{cor:epszero}
Let $V\in C_0^\infty(\R^3)$ satisfy $\|V\|_R<4\pi$. Then for all $f,g\in C_0^\infty(\R^3)$ the limit 
\[ \langle R_V(\lambda+i0)f,g \rangle = \lim_{\eps\to0} \langle R_V(\lambda+i\eps)f,g \rangle\]
exists for every~$\lambda\in\R$ and is a continuous function in~$\lambda$. 
Moreover, for each $\lambda$ one can pass to the limit $\eps\to0$ in all other 
terms in~\eqref{eq:Resseries} and 
\begin{equation}
\label{eq:eps0series}
\langle R_V(\la+ i0) f,\,g\rangle - \langle R_0(\la+ i0) f,\,g\rangle = \sum_{\ell=1}^\infty (-1)^\ell \langle R_0(\la+ i0)(VR_0(\la+ i0))^\ell f,g\rangle
\end{equation}
holds for every~$\lambda$ and the series converges absolutely in the norm of $C^0_b(\R)\cap L^1(d\lambda)$. 
\end{cor}
\begin{proof}
Fix $f,g\in C_0^\infty(\R^3)$. By our assumptions on~$V$ and the explicit representation~\eqref{eq:R0kernel},
$VR_0(z)f\in C_0^\infty$, and thus also $R_0(z)(VR_0(z))^\ell f$ for every $z\in \Compl$ with~$\Im z\ge 0$.  
Moreover, $z\mapsto \langle R_0(z)(VR_0(z))^\ell f,g \rangle$ is a continuous function in~$\Im z\ge 0$ for every $\ell\ge0$.  As in the previous proof one obtains the Kato smoothing bound
\begin{equation}
\label{eq:katoagain} 
\sup_{\eps\ge0}\int \Bigl| \langle R_0(\la+i\eps)(VR_0(\la+i\eps))^\ell f,g \rangle \Bigr| \,d\lambda \le C(\|V\|_R/4\pi)^{\ell-1}
\|f\|_{L^2_x}\|g\|_{L^2_x}
\end{equation}
for each $\ell\ge1$ (note that the case $\eps=0$ is included here). Moreover, see~\eqref{eq:fundop} and~\eqref{eq:Klam}, 
\be
 \sup_{\Im z\ge0}\Bigl| \langle R_0(z)(VR_0(z))^\ell f,g \rangle \Bigr| &\le& 
 \sup_{\Im z\ge0} \|\,|V|^{\half} R_0(z)g\|_2\,\|K(z)\|^{\ell-1}\, \|\,|V|^{\half} R_0(z)f\|_2  \nn \\
&\le& C(f,g,V)\, (\|V\|_R/4\pi)^{\ell-1}. \nn 
\ee
This implies that 
\begin{equation}
\label{eq:Sfgdef} 
S_{f,g}(\la) := \sum_{\ell=0}^\infty  \langle R_0(\la+i0)(VR_0(\la+i0))^\ell f,g \rangle,
\end{equation}
converges uniformly and thus defines a continuous function. Furthermore, one concludes that the series in~\eqref{eq:Resseries} converges uniformly in the closed upper half-plane (i.e., for all $\lambda\in\R$ and $\eps\ge0$) and therefore defines the limit~$\langle R_V(\lambda+i0) f,g \rangle$ pointwise in~$\lambda\in\R$. 
Also note that, by~\eqref{eq:katoagain}, the series for $S_{f,g}(\la)-\langle R_0(\la+i0)f,g \rangle$ converges absolutely in~$L^1(d\lambda)$, and similarly for every $z\in \Compl$ with~$\Im z\ge 0$. 
In view of \eqref{eq:Resseries}, \eqref{eq:katoagain}, and with an arbitrary~$N\ge1$, 
\be
&& \int \liminf_{\eps\to0} |\langle R_V(\la+ i\eps) f,\,g\rangle - S_{f,g}(\la)| \, d\la \nn \\
&\le&  \int \liminf_{\eps\to0} |\langle S_{f,g}(\la+ i\eps) f,\,g\rangle - S_{f,g}(\la)| \, d\la \nn \\
&\le& \int  \sum_{\ell=1}^N \limsup_{\eps\to0} \Bigl| \langle R_0(\la + i\eps)(VR_0(\la+i\eps))^\ell f,g \rangle
- \langle R_0(\la + i0)(VR_0(\la+i0))^\ell f,g \rangle \Bigr| \, d\la \nn \\
&& \qquad\qquad\qquad + C \sum_{\ell=N+1}^\infty (\|V\|_R/4\pi)^{\ell-1} \|f\|_2\|g\|_2 \nn \\
&\le& C  (1-\|V\|_R/4\pi)^{-1}\,(\|V\|_R/4\pi)^{N} \|f\|_2\|g\|_2, \nn
\ee
and we are done.
\end{proof}

\noindent Next we turn to a simple lemma that is basically an instance of stationary phase.

\begin{lemma}
\label{lem:statphas}
Let $\psi$ be a smooth, even  bump function
with $\psi(\lambda)=1$ for $-1\le\lambda\le 1$ and $\supp(\psi)\subset[-2,2]$. 
Then for all $t\ge1$ and any real~$a$,
\begin{equation}
\label{eq:decay}
\sup_{L\ge 1}\Bigl| \int_0^\infty e^{it\lambda} \sin(a\sqrt{\lambda})\,
\psi(\frac{\sqrt\lambda}{L})\,d\lambda\Bigr| \le
C \,t^{-\frac32}\,|a|
\end{equation}
where $C$ only depends on ~$\psi$.
\end{lemma}
\begin{proof}
Denote the integral in \eqref{eq:decay} by $I_L(a,t)$. 
Clearly, $I_L(a,t)$ is a smooth function of $a,t$ for any $L>0$
and~$I_L(0,t)=0$. The change of variables $\la\to \la^2$ leads to the expression
$$
I_L(a,t) = 2 \int_0^\infty \la \,e^{it\lambda^2} \sin (a\lambda)\,\psi(\lambda/L)\,d\lambda 
$$
Integrating by parts we obtain
$$
I_L(a,t) = \frac i{t}\int_0^\infty e^{it\lambda^2} 
\bigg(a\,\cos (a\lambda)\,\psi(\lambda/L) +\frac {1}{L} \sin (a\lambda)\,\psi'(\lambda/L)\bigg) 
\,d\lambda. 
$$
Since $\psi$ is assumed to be even, $\psi'$ is odd. Hence, 
$$
\aligned
I_L(a,t)& = \frac i{2t}\int_{-\infty}^\infty e^{it\lambda^2} 
\bigg(a\,\cos (a\lambda)\,\psi(\lambda/L) +\frac {1}{L} \sin (a\lambda)\,\psi'(\lambda/L)\bigg) 
\,d\lambda \\ &=  \frac {a}{4t}\,i\int_{-\infty}^\infty e^{it\lambda^2} 
\big(\,e^{i a\lambda} + e^{-i a\la}\big)\,\psi(\lambda/L)\,d\la + 
\int_0^a  \frac {i}{4t}\,\int_{-\infty}^\infty e^{it\lambda^2} 
\big(\,e^{i b\lambda} + e^{-i b\la}\big)\,\frac {\la}L\,\psi'(\lambda/L)\,d\la \,db.
\endaligned
$$ 
Thus it suffices to show that 
$$
J_L(a,t)= \int_{-\infty}^\infty e^{i(t\lambda^2 + a\lambda)}\,\phi(\lambda/L)\,d\la
$$
obeys the estimate $|J_L(a\,t)|\le C t^{-\frac 12}$ for any smooth bump function $\phi$ 
satisfying the same properties as $\psi$.
The change of variables $\la\to \la/L$ further reduces the problem to the 
estimate $|J(a,t)|\le C t^{-\frac 12}$ with 
$$
J(a,t)= \int_{-\infty}^\infty e^{i(t\lambda^2 + a\lambda)}\,\phi(\lambda)\,d\la
$$
for all $t\ne 0$ and all real $a$.
Observe that $J(a,t)$ is a smooth solution of the 1-dimensional Schr\"odinger
equation 
$$
\aligned
\ &i\, \frac{\partial}{\partial t} J(a,t) - \frac{\partial^2}{\partial a^2} J(a,t)=0,\\
\ &J(a,0)= \int_{-\infty}^\infty e^{-ia\la} \phi(\la)\,d\la.
\endaligned
$$
By the explicit representation of the kernel of the fundamental solution
$$
J(a,t)=  (-4\pi it)^{-\frac 12}\int_{-\infty}^\infty 
e^{-i\frac {|a-b|^2}{4t}} J(b,0)\,db
$$
which implies that $J(a,t)$ obeys the standard one-dimensional decay estimate 
$$
|J(a,t)|\le C t^{-\frac 12} \|J(\cdot,0)\|_{L^1}.
$$
Since the function $J(a,0)$ is the Fourier transform of the smooth bump function~$\phi$,
the desired estimate on $J(a,t)$ follows.
\end{proof}

\noindent The following lemma explains to some extent why condition~\eqref{eq:kato} is needed. 
Iterated integrals as in~\eqref{eq:multint} will appear in a series expansion of the spectral resolution of $H=-\Laplace+V$.

\begin{lemma}
\label{lem:iter} 
For any positive integer $k$ and $V$ as in Definition~\ref{def:V}
\begin{equation}
\label{eq:multint}
\sup_{x_0,x_{k+1}\in\R^3}\int_{\R^{3k}} \frac{\prod_{j=1}^k |V(x_j)|}{\prod_{j=0}^k|x_j-x_{j+1}|}\sum_{\ell=0}^k |x_\ell-x_{\ell+1}|
\; dx_1\ldots\,dx_k \le (k+1) \|V\|_{\kato}^k.
\end{equation}
\end{lemma}
\begin{proof} 
Define the operator ${\mathcal A}$ by the formula
$$
{\mathcal A} f (x) = \int_{\R^3} \frac{|V(y)|}{|x-y|} \,f(y)\,dy.
$$
Observe that the assumption \eqref{eq:kato} on the potential $V$ implies
that ${\mathcal A}:\,L^\infty\to L^\infty$ and $\|{\mathcal A}\|_{L^\infty\to L^\infty}\le c_0$ where we have set $c_0:=\|V\|_{\kato}$ for convenience. 
Denote by $<,>$ the standard $L^2$ pairing.
In this notation the estimate \eqref{eq:multint} is equivalent to proving 
that the operators ${\mathcal B}_k$ defined as  
$$
{\mathcal B}_k f = \sum_{m=0}^k <f,{\mathcal A}^{k-m} 1>{\mathcal A}^m 1 
$$
are bounded as operators from $L^1\to L^\infty$ with the bound
$$
\|{\mathcal B}_k\|_{L^1\to L^\infty}\le (k+1) c_0^k.
$$
For arbitrary $f\in L^1$ one has
$$
\aligned
\|{\mathcal B}_k f\|_{L^\infty}\le  \sum_{m=0}^k |<f,{\mathcal A}^{k-m} 1>|\,\,
\|{\mathcal A}^m 1\|_{L^\infty}\ & \le \sum_{m=0}^k \|{\mathcal A}^{k-m}\|_{L^\infty\to L^\infty}
 \|{\mathcal A}^{m}\|_{L^\infty\to L^\infty} \|f\|_{L^1} \\ &\le 
\sum_{m=0}^k c_0^k \|f\|_{L^1}\le (k+1) c_0^k \|f\|_{L^1},
\endaligned
$$ 
as claimed.
\end{proof}

\noindent We are now in a position to prove the main result of this section.

\begin{theorem}
\label{thm:high}
With $H=-\Laplace+V$ and~$V$ 
satisfying the conditions in Definition~\ref{def:V} one has the bound
\[ \Bigl\|e^{itH}\Bigr\|_{L^1\to L^\infty} \le C\, t^{-\frac32}\]
in three dimensions. 
\end{theorem}
\begin{proof}
Let $\psi$ be a smooth cut-off function as in Lemma~\ref{lem:statphas}. We will show that there is an absolute constant~$C$
such that 
\begin{equation}
\label{eq:reduc}
\sup_{L\ge1} \Bigl|\bigl \langle e^{itH}\psi(\sqrt{H}/L) f,g \bigr\rangle \Bigr| \le C t^{-\frac32} \|f\|_1\|g\|_1
\end{equation}
for any $f,g \in C^\infty_0(\R^3)$, which proves the theorem.
It will be convenient to assume that the potential~$V$ belongs to~$C_0^\infty(\R^3)$, 
in addition to satisfying~\eqref{eq:rollnik} and~\eqref{eq:kato}.
In case of a general potential~$V$ as in Definition~\ref{def:V}, one approximates $V$ by $V_j\in C_0^\infty$
via the usual cut-off and mollifying process. Clearly,  $\|V-V_j\|_R\to 0$ as $j\to\infty$ and  $\|V_j\|_{\kato}\le \|V\|_{\kato}<4\pi$.  
Since the spectral resolution $E_V$ of $H$ satisfies (recall that the spectrum of $H$ is purely absolutely
continuous)
\begin{equation}
\label{eq:specres}
E_V'(\lambda):=\frac{d}{d\lambda} E_V(\lambda) = \Im R_V(\lambda+i0),
\end{equation}
one concludes from Lemma~\ref{lem:resolvent} that
\[ \int \Bigl| \langle E'_V(\lambda) f,g \rangle -  \langle E'_{V_j}(\lambda)f,g \rangle\Bigr|\,d\lambda \to 0\]
as $j\to\infty$. In particular, with $H_j := -\Laplace+V_j$, 
\[ \bigl \langle e^{itH_j} \psi(H_j/L) f,g \bigr\rangle \to \bigl \langle e^{itH} \psi(H/L) f,g \bigr\rangle\]
as $j\to\infty$ for any $f,g \in C^\infty_0(\R^3)$. 
It therefore suffices to prove \eqref{eq:reduc} under the additional assumption 
that~$V\in C_0^\infty(\R^3)$.  
Fix such a potential~$V$, as well as any $L\ge1$, and real $f,g\in C_0^\infty(\R^3)$.
Then applying \eqref{eq:specres}, Corollary~\ref{cor:epszero}, \eqref{eq:R0kernel}, Lemma~\ref{lem:statphas}, and Lemma~\ref{lem:iter} in this order, 
\be
&& \sup_{L\ge1} \Bigl|\bigl \langle e^{itH}\psi(\sqrt{H}/L) f,g \bigr\rangle \Bigr| \nonumber\\
&\le& 
\sup_{L\ge1} \Bigl| \int_0^\infty e^{it\lambda}\;\psi(\sqrt{\lambda}/L) \langle E'(\lambda)f,g 
\rangle \,d\lambda \Bigr| \nonumber \\
&=& \sup_{L\ge1} \Bigl|\int_0^\infty e^{it\la}\, \psi(\sqrt{\lambda}/L) \Im 
\langle R_V(\lambda+i0)f,g \rangle \,d\lambda \Bigr| \nonumber\\
&=& \sup_{L\ge1} \Bigl| \int_0^\infty e^{it\lambda}\;\psi(\sqrt{\lambda}/L) 
\sum_{k=0}^\infty \Im \langle R_0(\lambda+i0)(VR_0(\lambda+i0))^k\,f,g \rangle \,
d\lambda \Bigr| \nonumber\\
&\le& \sum_{k=0}^\infty \int_{\R^6} |f(x_0)||g(x_{k+1})|
\int_{\R^{3k}} \frac{\prod_{j=1}^k |V(x_j)|}{\prod_{j=0}^k 4\pi |x_j-x_{j+1}|}\cdot\nonumber\\
&& \qquad\qquad\qquad \cdot \sup_{L\ge1} \Bigl| \int_0^\infty e^{it\lambda}\; 
\psi(\sqrt{\lambda}/L) \sin\Bigl(\sqrt{\lambda}\sum_{\ell=0}^k |x_\ell-x_{\ell+1}|\Bigr)\,
d\lambda \Bigr|
\; d(x_1,\ldots,x_k)\,dx_0\,dx_{k+1} \label{eq:gross}\\
&\le& Ct^{-\frac32} \sum_{k=0}^\infty \int_{\R^6} |f(x_0)||g(x_{k+1})| 
\int_{\R^{3k}} \frac{\prod_{j=1}^k |V(x_j)|}{(4\pi)^{k+1}\prod_{j=0}^k|x_j-x_{j+1}|}\sum_{\ell=0}^k |x_\ell-x_{\ell+1}|
\; d(x_1,\ldots,x_k)\;dx_0\,dx_{k+1} \nonumber \\
&\le& Ct^{-\frac32} \sum_{k=0}^\infty \int_{\R^6} |f(x_0)||g(x_{k+1})|\; 
(k+1) (\|V\|_{\kato}/4\pi)^k \;dx_0\,dx_{k+1}\nonumber \\
&\le& Ct^{-\frac32} \|f\|_1\|g\|_1, \nonumber
\ee
since $\|V\|_{\kato}<4\pi$. In order to pass to \eqref{eq:gross} one uses the explicit representation 
of the kernel of $R_0(\lambda+i0)$, see~\eqref{eq:R0kernel}, which leads to a $k$-fold integral. Next, one interchanges the
order of integration in this iterated integral. This is legitimate, since the corresponding $L^1$-integral (i.e., with absolute values on everything)
is finite ($V$ is bounded and compactly supported). 
The theorem follows.
\end{proof}

\section{The high energy case in $\R^3$ with an $\epsilon$ loss}\label{sec:epsloss}

\noindent The purpose of this section is to prove a dispersive inequality for $e^{itH}\chi(H)P_{a.c}$ where $\chi$ is a cut-off to large energies and~$P_{a.c}$ is the projection onto the absolutely continuous part of~$L^2(\R^3)$ with respect to $H=-\Laplace+V$. We will assume that $V$ satisfies the following properties:
\begin{equation}
\label{eq:highEcond}
\drei V \drei := \|V\|_2 + \sup_{x\in\R^3} \int_{\R^3} \frac{|V(y)|}{|x-y|}\,dy < \infty.
\end{equation}

\noindent Under these conditions we will prove the following result. As usual, we let
$\chi\in C^\infty$ with $\chi(\lambda)=0$ if $\lambda\le1$ and $\chi(\lambda)=1$ for $\lambda\ge2$.

\begin{prop} 
\label{prop:higheps}
Let $\drei V\drei <\infty$ as in \eqref{eq:highEcond}. Then for every $\eps>0$ there exists some positive $\lambda_0=\lambda_0(\drei V\drei,\eps)$ so that
\begin{equation}
\label{eq:highdisp}
\| e^{itH} \chi(H/\lambda_0) P_{a.c.} \|_{L^1_x\to L^\infty_x} \le Ct^{-\frac32+\eps}
\end{equation}
for all $t>0$.
\end{prop}

\noindent Previously, convergence of the Born series was guaranteed by a smallness assumption on the potential~$V$. The following lemma will allow us to sum the Born series for large potentials in~$L^2(\R^3)$, but only for large energies. This lemma is an immediate consequence of the Stein-Tomas theorem in the formulation due to Stein~\cite{stein}. 

\begin{lemma}
\label{lem:stein} 
Let $R_0(z)= (-\Laplace -z)^{-1}$ for $\Im(z)>0$ be the resolvent of the free Laplacean.
Then there is an absolute constant $C$ so that for any $\lambda>0$ 
\begin{equation}
\label{eq:ST} 
\|R_0(\lambda+i0)f\|_{L^4(\R^3)} \le C\lambda^{-\frac14}\|f\|_{L^{\frac43}(\R^3)}
\end{equation}
for all $f\in\calS$ .
\end{lemma}
\begin{proof}
It is well-known that the resolvent $R_0(z)= (-\Laplace -z)^{-1}$ for $\Im(z)>0$ has the kernel
\begin{equation}
\label{eq:kernel} 
K_0(z)(x,y) = \frac{\exp(i\sqrt{z}|x-y|)}{4\pi|x-y|}
\end{equation}
where $\Im(\sqrt{z})>0$. By the Stein-Tomas theorem in Stein's version~\cite{stein} one has
\begin{equation}
\label{eq:eli} 
\Big\| \int_{\R^3} \frac{\exp(i|x-y|)}{4\pi|x-y|}\,f(y)\,dy \Big\|_{L^4(\R^3)} \le C \|f\|_{L^{\frac43}(\R^3)}.
\end{equation}
Passing to \eqref{eq:ST} only requires changing variables $x\mapsto \sqrt{\lambda} x$ and $y\mapsto \sqrt{\lambda}y$, which we skip.
\end{proof}

\noindent It is well-known, see Simon~\cite{barry} Theorem~A.2.9, that for any $V\in L^2_{loc}(\R^3)$ that 
satisfies the so called {\em Kato condition}
\begin{equation}
\label{eq:3dkato} 
\lim_{r\to 0}\sup_{x\in\R^3} \int_{|x-y|< r} \frac{|V(y)|}{|x-y|}\, dy = 0,
\end{equation}
the operator $-\Laplace+V$ with domain $C_0^\infty(\R^3)$ is essentially self-adjoint with $\spec(H)\subset[-M,\infty)$ for some $0< M<\infty$, and that~$H$ is the generator of a semi group $e^{-tH}$ that is bounded from~$L^p$ to~$L^q$ for any choice of~$1\le p\le q\le\infty$, see Theorem~B.1.1 in~\cite{barry}. Moreover, explicit bounds for these norms are of the form
\[ \|e^{-tH}\|_{L^p_x\to L^q_x} \le C\,t^{-\gamma}\,e^{At}\]
with $\gamma=\frac32(p^{-1}-q^{-1})$ and any $A>M$ with~$M$ as before, see~(B11) in~\cite{barry}. These bounds imply the Sobolev inequalities
\begin{equation}
\label{eq:sobH} 
\|(H+2M)^{-\beta} \|_{L^p_x\to L^q_x} <\infty \text{\ \ for any $1\le p,q\le\infty$ and with $\beta > \frac32(p^{-1}-q^{-1})$,}
\end{equation}
as can be seen from writing $(H+2M)^{-\beta}$ as the Laplace transform of the heat-semigroup, see Theorem~B.2.1.~in~\cite{barry}. 
Since we are assuming that $V\in L^2$, Cauchy-Schwarz implies that~\eqref{eq:3dkato} holds, and thus so do all aforementioned properties. In addition, we will use the following result of Jensen and Nakamura, see~\cite{JN} Theorem~2.1: Suppose that $V\in L^2_{loc}$ satisfies~\eqref{eq:3dkato}. Let $g\in C_0^\infty(\R)$ and $1\le p\le\infty$. Then there exists a constant~$C$ such that
\begin{equation}
\label{eq:unifH} 
\|g(\theta H)\|_{L^p_x\to L^p_x} \le C \text{\ \ uniformly in $0<\theta\le 1$.} 
\end{equation}
Moreover, the constant~$C$ is uniform for $g$ ranging over bounded sets of~$C_0^\infty(\R)$. 
As an immediate corollary of \eqref{eq:sobH} and \eqref{eq:unifH} one obtains that for any 
$g\in C_0^\infty(\R)$, any $1\le p\le q\le\infty$, $\beta > \frac32(p^{-1}-q^{-1})$, there is a constant~$C$ depending on~$g,V$, and~$\beta$, such that
\begin{equation}
\label{eq:LitPalsob} 
\| g(H/\lambda_0) \|_{L^p_x\to L^q_x} \le C \lambda_0^{\beta} \text{\ \ uniformly in $\lambda_0\ge1$.}
\end{equation}

\noindent This bound is needed for the following lemma. Recall that $R_0(z)$ denotes the resolvent of the free Laplacean. 

\begin{lemma}
\label{lem:igor} 
Let $\eta\in C_0^\infty(\R)$ be fixed. 
Then for any $\lambda, \lambda_0\ge 1$ and any nonnegative integer~$k$ one has the estimate
\[ 
\|\eta(H/\lambda_0) R_0(\lambda+i0)(VR_0(\lambda+i0))^k \eta(H/\lambda_0)\|_{L^1_x\to L^\infty_x} \le C\lambda_0^{\frac34+}\,
\lambda^{-\frac14}(\|V\|_2 \lambda^{-\frac14})^k 
\]
where the constant $C$ only depends on $g$ and $V$. 
\end{lemma}
\begin{proof} 
By Lemma~\ref{lem:stein} and H\"older's inequality,
\begin{equation}
\label{eq:smallVR} 
\|VR_0(\lambda+i0)f\|_{L^{\frac43}_x} \le \|R_0(\lambda+i0)f\|_{L^4_x}\|V\|_{L^2_x} \le C\,\|V\|_{L^2_x}\,\lambda^{-\frac14}\,\|f\|_{L^{\frac43}_x}
\end{equation}
for any $f\in\calS$. Hence
\be 
&& \| \eta(H/\lambda_0) R_0(\lambda+i0)(V R_0(\lambda+i0))^k\, \eta(H/\lambda_0)\|_{L^1_x\to L^\infty_x} \nonumber \\
&& \le \|\eta(H/\lambda_0) R_0(\lambda+i0)\|_{L^{\frac43}_x\to L^\infty_x}\, \|(VR_0(\lambda+i0))^k\|_{L^{\frac43}_x\to L^{\frac43}_x}\,
\|\eta(H/\lambda_0)\|_{L^1_x\to L^{\frac43}_x} \nonumber \\
&& \le \|\eta(H/\lambda_0)\|_{L^4_x\to L^\infty_x} \|R_0(\lambda+i0)\|_{L^{\frac43}_x\to L^4_x}\,\|VR_0(\lambda+i0)
\|^k_{L^{\frac43}_x\to L^{\frac43}_x}\, \|\eta(H/\lambda_0)\|_{L^1_x\to L^{\frac43}_x} \nonumber \\
&& \le C\lambda_0^{\frac38+} \lambda^{-\frac14} (\|V\|_2\lambda^{-\frac14})^k \,\lambda_0^{\frac38+}, \nonumber 
\ee
as claimed.
\end{proof}

\begin{proof}[Proof of Proposition \ref{prop:higheps}] 
We start with a justification of the Born series expansion for high energies. 
Let $\lambda_0>0$ be chosen so that $\|V\|_{L^2_x}\,\lambda_0^{-\frac14}<1$. 
By~\eqref{eq:smallVR}, the operator $1+VR_0(\lambda+i0)$ is invertible in $L^{\frac43}(\R^3)$ provided
$\lambda>\lambda_0$  and the Neumann series
\begin{equation}
\label{eq:neumann}
(1+VR_0(\lambda+i0))^{-1} = \sum_{k=0}^\infty (-1)^k (VR_0(\lambda+i0))^k 
\end{equation}
converges in $L^{\frac43}(\R^3)$. Therefore,  the resolvent $R_V(z):=(-\Laplace + V -z)^{-1}$ satisfies 
\[ R_V(\lambda+i0) = R_0(\lambda+i0) (1+VR_0(\lambda+i0))^{-1}= \sum_{k=0}^\infty (-1)^k R_0(\lambda+i0)(VR_0(\lambda+i0))^k \]
for all $\lambda>\lambda_0$ and is thus a bounded operator from $L^{\frac43}(\R^3)\to L^4(\R^3)$. 
Furthermore, since the spectral resolution $E(\cdot)$ of $H=-\Laplace + V$ satisfies $P_{a.c.}\,E(d\lambda) = \Im R_V(\lambda+i0) \,d\lambda$, one has
\[ 
\langle P_{a.c.} E(d\lambda) f, g \rangle = \sum_{k=0}^\infty (-1)^k \langle \Im [R_0(\lambda+i0) (VR_0(\lambda+i0))^k] f,g \rangle\, d\lambda
\]
for any $f,g\in L^{\frac43}(\R^3)$. Now define $\eta(\lambda):= \chi(\lambda)-\chi(\lambda/2)$. 
Clearly, $\eta\in C_0^\infty(\R)$, and also 
\[ \sum_{j=0}^\infty \eta(\lambda 2^{-j}) = \chi(\lambda) \text{\ \ for all $\lambda$}.\] 
Observe that at most three terms in this sum can be nonzero for any given~$\lambda$.  
Now let $\tilde{\eta}\in C_0^\infty(0,\infty)$ have the property that $\eta\tilde{\eta}=1$. 
Then for any $f,g\in \calS$, one has the expansion
\be
&& |\langle e^{itH} \chi(H/\lambda_0) f , \chi(H/\lambda_0) g \rangle|  \nn \\
&=& \left| \int_0^\infty e^{it\lambda} \langle E(d\lambda) \eta(H/(2^j\lambda_0)) f, \eta(H/(2^\ell\lambda_0)) g \rangle \, d\lambda \right| \nn \\
&\le& \biggl| \sum_{\substack{j,\ell=0 \\ |j-\ell|\le 1}}^\infty  \int_0^\infty e^{it\lambda}  \langle E(d\lambda) \eta(H/(2^j\lambda_0)) f, \eta(H/(2^\ell\lambda_0)) g \rangle \, \tilde{\eta}(\lambda/(2^j\lambda_0))\, d\lambda \biggr| \nn \\
&\le& \sum_{k=0}^\infty \biggl| \sum_{\substack{j,\ell=0 \\ |j-\ell|\le 1}}^\infty 
 \int_0^\infty e^{it\lambda} \langle R_0(\lambda+i0)\,(VR_0(\lambda+i0))^k \eta(H/(2^j\lambda_0)) f,\eta(H/(2^\ell\lambda_0)) g \rangle  \nn \\
&& \qquad\qquad\qquad\qquad\qquad\qquad\qquad\qquad\qquad\qquad\qquad\qquad\qquad  \tilde{\eta}(\lambda/(2^j\lambda_0))\, d\lambda 
\biggr| \label{eq:langform} \\
&=& \sum_{k=0}^\infty \left| 
 \int_0^\infty e^{it\lambda} \langle R_0(\lambda+i0)\,(VR_0(\lambda+i0))^k \chi(H/\lambda_0) f, 
\chi(H/\lambda_0) g \rangle \, d\lambda \right|. \label{eq:langform2}
\ee
From the previous section one has the dispersive bounds
\be
\label{eq:3/2}
&& \quad \left| \int_0^\infty e^{it\lambda} \langle R_0(\lambda+i0)\,(VR_0(\lambda+i0))^k \chi(H/\lambda_0) f, 
\chi(H/\lambda_0) g \rangle \, d\lambda \right| \le C \, t^{-\frac32} \drei V\drei ^k\,\|f\|_{L^1_x}\,\|g\|_{L^1_x}, \\ 
&& \quad \left| \int_0^\infty e^{it\lambda} \langle R_0(\lambda+i0)\,(VR_0(\lambda+i0))^k \eta(H/(2^j\lambda_0)) f, \eta(H/(2^\ell\lambda_0)) g \rangle \,\tilde{\eta}(\lambda/(2^j\lambda_0))\, d\lambda \right| \nn \\
&& \qquad\le C \, t^{-\frac32} \drei V\drei ^k\, \|f\|_{L^1_x}\,\|g\|_{L^1_x}, 
\label{eq:3/2'}
\ee
where we have also used \eqref{eq:LitPalsob} to remove the $\chi$ and $\eta$ cutoffs. On the other hand, Lemma~\ref{lem:igor} shows that
\be
&&\left| \int_0^\infty e^{it\lambda} \langle R_0(\lambda+i0)\,(VR_0(\lambda+i0))^k \eta(H/(2^j\lambda_0)) f, \eta(H/(2^\ell\lambda_0)) g \rangle \, \tilde{\eta}(\lambda/(2^j\lambda_0))\, d\lambda \right| \nn \\ 
&\le& C \, (2^j\lambda_0)^{\frac32+} (\|V\|_{L^2_x} (2^j\lambda_0)^{-\frac14})^k \|f\|_{L^1_x}\,\|g\|_{L^1_x}. \label{eq:noosc}
\ee
Combining \eqref{eq:3/2'} and \eqref{eq:noosc} yields that for any $0<\theta<1$ 
\be
&& \sum_{k=7}^\infty \sum_{\substack{j,\ell=0 \\ |j-\ell|\le 1}}^\infty 
 \left| \int_0^\infty e^{it\lambda} \langle R_0(\lambda+i0)\,(VR_0(\lambda+i0))^k \eta(H/(2^j\lambda_0)) f, \eta(H/(2^\ell\lambda_0)) g \rangle \, \tilde{\eta}(\lambda/(2^j\lambda_0))\, d\lambda \right| \nn \\ 
&&\qquad\le C \, \sum_{j=0}^\infty \sum_{k= 7}^\infty t^{-\frac32(1-\theta)}\,\drei V\drei^{k(1-\theta)}\, (2^j\lambda_0)^{\theta \frac32+} 
(\|V\|_{L^2_x} (2^j\lambda_0)^{-\frac14})^{\theta k} \|f\|_1\,\|g\|_1 \nn \\
&&\qquad\le C \, \sum_{k=7}^\infty t^{-\frac32(1-\theta)} \,\drei V\drei^{k(1-\theta)}\, \lambda_0^{\theta \frac32+} (\|V\|_2 \lambda_0^{-\frac14})^{\theta k} \|f\|_{L^1_x}\,\|g\|_{L^1_x} \nn \\
&&\qquad\le C\, t^{-\frac32(1-\theta)} \lambda_0^{\theta \frac32+} \sum_{k=0}^\infty \drei V\drei^k \lambda_0^{-\frac{k}{4}\theta} \|f\|_1\,\|g\|_1 \le C\, t^{-\frac32(1-\theta)}\,\lambda_0^{\theta \frac32+} \|f\|_{L^1_x}\|g\|_{L^1_x}
\label{eq:uffeps}
\ee
provided $\drei V\drei \lambda_0^{-\frac{\theta}{4}}<1$. The choice of $k\ge7$ was made to ensure summability over~$j$.  The bound~\eqref{eq:uffeps} yields the desired bounds for the terms with $k\ge 7 $ in~\eqref{eq:langform}. For the remaining cases of~$k$, one simply invokes the estimate~\eqref{eq:3/2}, and the proposition follows.
\end{proof}

\begin{remark}
It seems clear that the condition $\|V\|_{L^2_x}<\infty$ can be weakened to a condition closer to~$L^{\frac32}(\R^3)$.
The reason for this is the ``slack'' in the Stein-Tomas bound that yields~$\lambda^{-\frac14}$, whereas
the high energies argument only requires~$\lambda^{-\gamma}$ for some~$\gamma>0$. It appears that
a complex interpolation argument allows one to exploit this slack, but we do not pursue this here.
\end{remark}

\section{Strichartz estimates for $(1+|x|^2)^{-1-\eps}$ potentials}\label{sec:strichartz} 

In this section we settle a problem posed by Journ\'e, Soffer, Sogge \cite{JSS} concerning Strichartz estimates for the solutions
of the Schr\"odinger equation with potentials decaying at the rate of $|x|^{-2-\eps}$ at infinity.
To obtain the result we prove a more general statement relating an $L^q_t L^p_x$ estimate for the semigroup $e^{it H_0}$ to the
corresponding estimate for $e^{it H}$ with $H= H_0 +V$. The conditions of the result involve the notion of Kato's smoothing for 
the multiplication operator $|V|^{\frac 12}$ relative to $H_0$ and $H$. 
Applying the abstract result to $H_0=-\Delta$, $H=-\Delta + V$ with $V$ obeying the estimate 
$|V(x)|\le C (1+|x|^2)^{-1-\eps}$ requires appealing to the Agmon-Kato-Kuroda theory on the absence of positive singular continuous spectrum for $H$ and a separate argument that deals with the point~$0$ in the spectrum of $H$.

We start with the preliminaries.
Consider a self-adjoint operator $H_0$  on $L^2(\R^n)$ with domain 
${\mathcal D}(H_0)$.  Let $e^{it H_0}$ be the
associated unitary semigroup, which is a solution operator for the Schr\"odinger equation 
$$
\frac 1i \partial_t \psi - H_0 \psi =0,\qquad \psi|_{t=0} = \psi_0.
$$
 We denote by $R_0(z)$ the resolvent of $H_0$.
For complex $z$ with $\Im z >0$ we have that
\begin{equation}
\label{eq:respro}
R_0(z)= \il_0^\infty e^{izt} e^{itH_0}\,dt
\end{equation}
as well as the inverse: for any $\b >0$ and $t\ge 0$,
$$
e^{-\b t}  e^{itH_0} = \il_{-\infty}^{\infty} e^{-it\la} R_0(\la+i\b)\,d\la
$$
Let $A$ and $B$ be a pair of bounded operators\footnote{The assumption of boundedness
is a convenience that is sufficient for our main application below.}  on $L^2(\R^n)$
and consider a self-adjoint operator $H=H_0 + B^*A$ 
with domain  ${\mathcal D}(H_0)$,
corresponding semigroup $e^{itH}$, and the resolvent $R(z)$.
The resolvent $R(z)$ and $R_0(z)$ for $\Im z\ne 0$
are connected via the second resolvent identity
\begin{equation}
\label{eq:Resid}
R(z) = R_0(z) - R_0(z) B^* AR(z)
\end{equation}
On the other hand,
the semigroups $e^{itH}$ and $e^{it H_0}$ are related via the Duhamel
formula
\begin{equation}
\label{duh}
e^{it H} \psi_0 = e^{it H_0} \psi_0 - 
i \il_0^t e^{i(t-s) H_0} B^*A e^{is H}\psi_0\,ds.
\end{equation}
which holds for any $\psi_0\in L^2_x$.
We recall that for a self-adjoint operator $\bar H$, 
an operator $\Gamma$ is called $\bar H$-smooth in Kato's sense if 
for any $f\in  {\mathcal D}(H_0)$
\begin{equation}
\label{eq:Csmooth}
\|\Gamma e^{it \bar H} f\|_{L^2_t L^2_x}\le C_{\Gamma}(\bar H) \|f\|_{L^2_x}
\end{equation}
or equivalently, for any $f\in L^2_x$ 
\begin{equation}
\label{eq:CsmoothR}
\sup_{\b >0} \|\Gamma R_{\bar H}(\la\pm i\b) f\|_{L^2_\la L^2_x}
\le C_{\Gamma}(\bar H)\|f\|_{L^2_x}.
\end{equation}
We shall call $C_{\Gamma}(\bar H)$ the smoothing bound of $\Gamma$ relative to $\bar H$.
Let $\Omega\subset \R$ and let $P_\Omega$ be a spectral projection of $\bar H$ 
associated with a set $\Omega$. We say that  $\Gamma$ is $\bar H$-smooth on 
$\Omega$ if $\Gamma P_{\Omega}$ is $\bar H$-smooth. We denote the corresponding
smoothing bound by $C_\Gamma(\bar H,\Omega)$. It is not difficult to
show (see e.g. \cite{RS4}) that, equivalently,  $\Gamma$ is $\bar H$-smooth on 
$\Omega$ if
\begin{equation}
\label{eq:smon}
\sup_{\b >0, \la\in \Omega} \|\Gamma R_{\bar H}(\la\pm i\b ) f\|_{L^2_\la L^2_x}
\le C_{\Gamma}(\bar H, \Omega)\|f\|_{L^2_x}.
\end{equation}
We now are ready to state the main result of this section.

\begin{theorem}
Let $H_0$ and $H= H_0 + B^* A$ be as above.  We assume that 
that $B$ is $H_0$ smooth with a smoothing bound $C_B(H_0)$ and 
that for some $\Omega\subset \R$ the operator $A$ is 
$H$-smooth on $\Omega$ with the smoothing bound
$C_A(H,\Omega)$. 
Assume also that the unitary semigroup
$e^{it H_0}$ satisfies the estimate
\begin{equation}
\label{eq:strH0}
\|e^{it H_0} \psi_0\|_{L^q_t L^r_x} \le C_{H_0} \|\psi_0\|_{L^2_x}
\end{equation}
for some $q\in (2,\infty]$ and $r\in [1,\infty]$. 
Then the semigroup $e^{itH}$ associated with  
$H = H_0 + B^*A$,
restricted to the spectral set $\Omega$,
also verifies the estimate \eqref{eq:strH0}, i.e.,
\begin{equation}
\label{eq:strHk}
\|e^{it H}P_\Omega \psi_0\|_{L^q_t L^r_x} \le C_{H_0}C_B(H_0) C_A(H,\Omega) 
\|\psi_0\|_{L^2_x}
\end{equation}
\label{th:Katoth}
\end{theorem}
\begin{proof}
We start with the Duhamel formula \eqref{duh}
$$
e^{it H} \psi_0 = e^{it H_0} \psi_0 - 
i \il_0^t e^{i(t-s) H_0} B^*A e^{is H}\psi_0\,ds.
$$
We have the following estimate with the exponents $q,r$ described in \eqref{eq:strH0}:
\begin{align}
\|e^{it H} P_\Omega \psi_0\|_{L^q_t L^r_x}&\le \|e^{itH_0} P_\Omega\psi_0\|_{L^q_t L^r_x} +
\|\il_0^t e^{i(t-s) H_0} B^*A e^{is H}P_\Omega \psi_0\,ds\|_{L^q_t L^r_x}\nn\\ &
\le  C_{H_0} \|\psi_0\|_{L^2_x} +
\|\il_0^t e^{i(t-s) H_0} B^*A e^{is H} P_\Omega \psi_0\,ds\|_{L^q_t L^r_x}\label{eq:dbound}
\end{align}
To handle the Duhamel term we recall the Christ-Kiselev lemma.
The following version is from Sogge, Smith~\cite{sogge}

\begin{lemma}[CK]
Let $X, Y$ be Banach spaces and let $K(t,s)$ be the kernel of the 
operator $K: L^p([0,T]; X)\to L^q([0,T]; Y)$. Denote by $\|K\|$ the 
operator norm of $K$. Define the lower diagonal operator 
$\tilde K:\,L^p([0,T]; X)\to L^q([0,T]; Y)$
$$
\tilde Kf(t) = \il_0^t K(t,s) f(s)\,ds
$$
Then the operator $\tilde K$ is bounded from $L^p([0,T]; X)$ to $L^q([0,T]; Y)$
and its norm $\|\tilde K\|\le c \|K\|$, provided that $p<q$.
\end{lemma}
We shall apply this lemma to the operator with kernel
$K(t,s) = e^{i(t-s)H_0} B^*$ acting between the spaces
$L^2([0,\infty);  L^2_x)$ and  $L^q([0,\infty);  L^r_x)$.
Observe that by the assumptions of Theorem \ref{th:Katoth},  
$q>2$ and thus the condition $q>p$ in Lemma [CK] is verified.

\noindent We can rewrite the Duhamel
term 
$$ 
D= \il_0^t e^{i(t-s) H_0} B^*A e^{is H}P_\Omega \psi_0\,ds
$$ in the form
$D = \tilde K \bigg (Ae^{i\cdot H} P_\Omega \psi_0\bigg )$.
Therefore,
\begin{equation}
\label{D}
\|D\|_{L^q_t L^r_x} \lesssim \|K\|_{L^2([0,\infty);  L^2_x)\to L^q([0,\infty);  L^r_x)}\,\, 
\|A e^{is H}\psi_0\|_{L^2_t L^2_x}
\end{equation}
We now need to estimate the norm of the operator $K$.
$$
\aligned
\|KF\|_{L^q_t L^r_x}& =  
\|\il_0^\infty e^{i(t-s)H_0} B^* F(s)\,ds\|_{L^q_t L^r_x} = 
\|e^{it H_0} \il_0^\infty e^{-is H_0} B^* F(s)\,ds\|_{L^q_t L^r_x}  \nn \\
\\ &\le C_{H_0} \|\il_0^\infty e^{-is H_0} B^* F(s)\,ds\|_{L^2_x}.
\endaligned
$$
The last inequality is the estimate \eqref{eq:strH0} for $e^{it H_0}$.
By duality
$$
\aligned
 \|\il_0^\infty e^{-is H_0} B^* F(s)\,ds\|_{L^2_x} &= 
\sup_{\|\phi\|_{L^2_x}=1} <\il_0^\infty e^{-is H_0} B^* F(s)\,ds,\phi>\\ &=
\sup_{\|\phi\|_{L^2_x}=1}\il_0^\infty ds <F(s), B e^{is H_0} \phi>\\ &\le
\|F\|_{L^2_t L^2_x}  \sup_{\|\phi\|_{L^2_x}=1} 
\| B e^{is H_0} \phi\|_{L^{2}_t L^2_x}\le C_B(H_0) \|F\|_{L^2_t L^2_x} \| \phi\|_{L^2_x},
\endaligned
$$
where the last inequality follows from $H_0$-smoothness of the operator $B$.
Thus 
the operator $K(t,s) = e^{i(t-s)H_0} A$ is bounded 
from $L^2([0,\infty); L^2_x)$ to $L^q([0,\infty); L^r_x)$. 
Therefore, back to \eqref{D}
\begin{equation}
\label{D1}
\|D\|_{L^q_t L^r_x} \le C_{H_0} C_B(H_0) \|\,A e^{is H} P_\Omega \psi_0\|_{L^2_t L^2_x}. 
\end{equation}
It remains to observe that since the operator $A$ is $H$-smooth on $\Omega$,
we have 
\begin{equation}
\label{eq:Aboun}
\|\,A e^{is H} P_\Omega \psi_0\|_{L^2_t L^2_x}\le C_A (H,\Omega) \|\psi_0\|_{L^2_x}.
\end{equation}
Thus, combining \eqref{eq:dbound}, \eqref{D1}, and  \eqref{eq:Aboun} 
we finally obtain
$$ 
\|e^{it H} \psi_0\|_{L^q_t L^r_x}\le C_{H_0}
C_B(H_0) C_A(H,\Omega) \|\psi_0\|_{L^2_x},
$$
as claimed. 
\end{proof}

\noindent We apply Theorem \ref{th:Katoth} in the situation where 
$H_0=-\Lap$ and $H=H_0 + V(x)$. We have the following family of Strichartz 
estimates for the semigroup $e^{-it\Lap}$ associated with $H_0=-\Lap$: 
\begin{equation}
\label{eq:striH0}
\|e^{-it \Lap} \psi_0\|_{L^q_t L^r_x} \le C \|\psi_0\|_{L^2_x},
\qquad \forall (q,r,n)\ne (2,\frac {2n}{n-2}, n), \quad 
\frac 2q=n(\frac 12 -\frac 1r),
\end{equation}
which hold for any $\psi_0\in L^2(\R^n)$. We introduce the factorization 
$$V= B^* A,\qquad B=|V|^{\frac 12},\qquad A=|V|^{\frac 12} \mbox{sgn} V, 
$$
and restrict our attention to the class of potentials satisfying the 
assumption that for all $x\in \R^n$
\begin{equation}
\label{eq:Vass}
|V(x)|\le C_V (1+|x|^2)^{-1 -\eps} 
\end{equation}
with some constants $C_V,\eps>0$. This assumption, in particular, places 
us in the framework of the Agmon-Kato-Kuroda and the Agmon-Kato-Simon 
theorems guaranteeing the absence of the positive singular continuous 
spectrum and positive eigenvalues. In fact, one only needs the $|x|^{-1-\eps}$
decay for their results to apply. We should note that for potentials satisfying
\eqref{eq:Vass} the absence of the singular continuous spectrum was established
by Ikebe \cite{Ik}.

\noindent In addition, the Weyl criterion 
implies that the essential spectrum of $H$ is the half-axis $[0,\infty)$.
However, without an appropriate smallness or sign assumption on $V$, the operator 
$H=-\Lap + V$ can have negative eigenvalues, thus destroying any hope to have
Strichartz estimates for $e^{it H}\psi_0$ for all initial data $\psi_0\in L^2$.
Therefore, we shall assume that the initial data are orthogonal to the eigenfunctions
corresponding to the possible negative eigenvalues. We achieve this in the following
simple manner. Let $P$ be a spectral projection of $H$ corresponding to 
the interval $\Omega =[0,\infty)$. Our goal is to prove Strichartz inequalities
for $e^{itH}$ restricted to the absolutely continuous spectrum of $H$.
We now state the result.
\begin{theorem}
Let $V$ be a potential verifying  \eqref{eq:Vass}.
In addition, we impose the condition that the point $\la=0$ 
in the spectrum of the operator $H=-\Lap +V$ is neither
an eigenvalue nor a resonance (see the discussion below, in particular Definition~\ref{de:eire}). 
Then if $P$ is the spectral projection of $H$ corresponding to 
the interval $[0,\infty)$ (on which $H$ is purely absolutely continuous), 
\begin{equation}
\label{eq:striH}
\|e^{it H} P \psi_0\|_{L^q_t L^r_x} \le C \|\psi_0\|_{L^2_x},
\qquad \forall (q,r,n),\,\,n\ge 3,\quad 
\frac 2q=n(\frac 12 -\frac 1r).
\end{equation}
\label{th:MaStr}
\end{theorem}

To apply Theorem \eqref{th:Katoth} we need to verify that $B$ is an $H_0$-smooth
operator and that $A$ is an $H$-smooth operator on $[0,\infty)$.
The first condition is easy to verify since by a result of Kato~\cite{kato} any function
$f\in  L^{p_1}\cap L^{p_2}$ with $1\le p_1 < n  <p_2\le \infty$ and
$n\ge 3$ is a $-\Lap$-smooth multiplication operator.  Since $B=|V|^{\frac 12}$
is an $L^\infty$ function decaying at infinity as $|x|^{-1-\eps}$, it falls precisely
under these conditions.

\noindent The condition that $A$ is an $H$-smooth operator on $[0,\infty)$ is much more subtle.
First, we can show that $A$ is $H$-smooth on the interval $[\de,\infty)$ for 
any $\de>0$. This is a consequence of the results of Agmon-Kato-Kuroda on the 
absence of the positive singular continuous spectrum, (see \cite{agmon}, also 
Theorem XIII.33 and Lemma 2 XIII.8 in \cite{RS4}). 
In fact, even half of the assumed decay would be
sufficient to prove this.  To deal with the remaining spectral interval 
$[0,\de)$, according to~\eqref{eq:smon}, one needs to understand the behavior
of the resolvent $R(\la\pm i\b)$ of the operator $H$ near the point $\la=0, \,\b=0$. 
We introduce the following 

\begin{defi} 
We say that $0$ is a regular point of the spectrum of $H$ if it is neither an
eigenvalue nor a resonance of $H$, i.e., the equation $-\Lap u + V(x) u=0$ has 
no solutions $u\in \cap_{\alpha >\frac 12} L^{2,-\alpha}$. 

\noindent Here, $L^{2,\alpha}$ is the weighted $L^2$ space of functions $f$ such that
$(1+|x|^2)^{\frac{\alpha}2} f\in L^2$. The 0 eigenvalue, of course, would correspond to 
an $L^2$ solution $u$.
\label{de:eire}
\end{defi}

\noindent The presence of a $0$ eigenvalue and most likely that of a resonance would violate the 
validity of the Strichartz estimates \eqref{eq:striH} for $e^{itH}$. Their appearance
cannot be ruled out by merely strengthening the regularity and decay assumptions
on the potential $V$. We therefore impose an additional condition that $0$ is 
a regular point. There are several situations where this condition, or at least
part of it, is automatically satisfied. In particular, for any non-negative potential
$0$ is a regular point. In addition, it is well-know (see e.g. \cite{JenKat}) that 
$0$ is not a resonance in dimensions $n\ge 5$.
The behavior of the resolvent near $0$ in the spectrum and even its
asymptotic expansions was extensively studied in \cite{JenKat}, \cite{Jensen1}, \cite{Jensen2},
but their assumptions are too strong for our purposes. 

\begin{proposition}
Suppose $V$ satisfies the assumption \eqref{eq:Vass} and assume, in addition, that~$0$ is 
a regular point of the spectrum of $H=-\Lap +V$. Then the operator 
$A=|V|^{\frac 12}\mbox{sgn}\,( V)$ is $H$-smooth on $[0,\de)$ for some 
sufficiently small $\de>0$.
\label{pr:LH}
\end{proposition}    

The first observation, which follows from \eqref{eq:smon},
is that since the potential $V$ decays at the rate
$|x|^{-2-2\eps}$ it suffices to prove the following property of the 
resolvent $R(\la\pm i\b)$ formulated in the language of the weighted 
spaces $L^{2,\alpha}$:
\begin{equation}
\label{eq:ressm}
\sup_{0<\b<\de, \la \in [0,\de)} \|R(\la\pm i\b) f\|_{L^{2,-1-\ga}}\le 
C \|f\|_{L^{2,1+\ga}}
\end{equation}
for any $f\in L^{2,1+\ga}$ and some sufficiently small $\de$ and $\ga$ such 
that $\ga <\eps$. The restriction of the range of $\b$ to the interval 
$(0,\de)$ is justified since for $\b\ge \de$ the resolvent $R(\la\pm i\b)$
in fact maps $L^2$ into the Sobolev space $W^2_2$ with a constant dependent only on $\de$.

\noindent One should compare \eqref{eq:ressm} with the standard limiting absorption 
principle which states that on the interval $[\de,\infty)$ the resolvent
$R(\la\pm i\b)$ is a bounded map between $L^{2,\frac 12 +\ga}$ and 
$L^{2,-\frac 12 -\ga}$ for any $\ga>0$. As in that case we reduce the proof to the 
same estimates for the free resolvent $R_0(\la\pm i\b)$.
The connection is established via the resolvent identity
$$
R(\la\pm i\b) = R_0(\la\pm i\b) -  R_0(\la\pm i\b) V R(\la\pm i\b).
$$
Thus formally we can solve for $R(\la\pm i\b)$,
\begin{equation}
\label{eq:RV}
R(\la\pm i\b) = (I +  R_0(\la\pm i\b) V)^{-1} R_0(\la\pm i\b).
\end{equation}
We identify the boundary value of the free resolvent at $\la=0,\, \b=0$ as 
the operator with the kernel given by the Green's function (up to constants),
$$
G (x,y):=R_0(0)(x,y) = \frac 1{|x-y|^{n-2}},\qquad n\ge 3
$$
and break the proof into a series of lemmas.
\begin{lemma}
$G$ is a bounded map 
from the weighted space $L^{2,1+\ga}$ into $L^{2,-1-\ga}$
for any $\ga >0$. Moreover, for any positive $\sigma <\frac 12$ we have
$G :\, L^{2,1+\ga +\sigma }\to L^{2,-1-\ga +\sigma }$.
\label{le:LAP}
\end{lemma} 
\begin{lemma}
The resolvent $R_0(\la\pm i\b)$ is continuous at $\la=0, \b=0$ in 
the topology of the bounded operators between  $L^{2,1+\ga}$ and $L^{2,-1-\ga}$
for any $\ga >0$.
\label{le:Cont} 
\end{lemma}
\begin{lemma}
Under the assumptions of Proposition \ref{pr:LH} the operator 
$G V$ is compact as an operator on $L^{2,-1-\ga}$ for any $\ga >0$, 
and  $(I+G V)$ is invertible on   $L^{2,-1-\ga}$.
\label{le:Vinv}
\end{lemma}
\begin{proof}[Proof of Proposition \ref{pr:LH}]
According to Lemma \ref{le:Vinv} the operator  $(I+G V)$ is invertible on   
$L^{2,-1-\ga}$ for any $\ga >0$. Therefore, by continuity of  
$R_0(\la\pm i\b)$ at $\la=0, \b=0$ asserted in Lemma \ref{le:Cont} and the 
fact that $V$ maps $L^{2,-1-\ga}$ to $L^{2, 1+\ga}$ provided that $\ga <\b$, there
exists a small neighborhood $\de$ of $0$ such that $(I+R_0(\la\pm i\b) V)$ is
uniformly invertible for all $|\la|, \b<\de$ on the space $L^{2,-1-\ga}$. 
It follows that for such $\la, \b$ the resolvent $R(\la\pm i\b)$ is well-defined
via the identity \eqref{eq:RV} and acts between the spaces 
$L^{2,1+\ga}$ and $L^{2,-1-\ga}$ as desired.
\end{proof}

\begin{proof}[Proof of Lemma \ref{le:LAP}]
The resolvent $G=R_0(0)$ is a multiplier with symbol $|\xi |^{-2}$. 
Therefore, after passing to the Fourier variables, 
$G:L^{2,1+\ga +\sigma }\to L^{2,-1-\ga +\sigma }$  is equivalent
to showing that multiplication by $|\xi|^{-2}$ acts between the Sobolev
spaces $W^{1+\ga+\sigma}_2$ and $W^{-1-\ga +\si}_2$.  
We consider the end-points of the desired values of $\ga$ and $\si$ 
corresponding to $\ga=0$ and $\si=0, \frac 12$ and prove that 
$$
|\xi|^{-2}\,:\,\,W^{1+}_2\to W^{-1-}_2, \qquad
|\xi|^{-2}\,:\,\,W^{\frac 32 +}_2\to W^{-\frac 12-}_2,
$$
where $\pm$ represent the fact that we do not prove the end-point
results themselves. Since  $|\xi|^{-2}$ is smooth away from $\xi=0$
we can consider instead the operator of multiplication by $\chi(\xi)|\xi|^{-2}$
where $\chi$ is a smooth cut-off function with support in a unit ball $B$. 
We have the standard Sobolev embeddings 
\begin{align}
&W^{1+}_2\einbet L^{\frac {2n}{n-2}+}, \label{eq:w12}\\
&W^{\frac 32 +}_2\einbet L^{\frac {2n}{n-3}} (L^\infty, \,n=3),\label{eq:w32}
\end{align}
the dual version of \eqref{eq:w12}, 
$L^{\frac {2n}{n+2}-} \einbet W^{-1-}_2$, 
and 
$L^{\frac{2n}{n+1}-}\einbet W^{-\frac 12-}_2$.
Therefore, we shall, in fact, prove a stronger result that 
$$
\chi(\xi)\,|\xi|^{-2}\,:\,\,L^{\frac {2n}{n-2}+}\to L^{\frac {2n}{n+2}-}, \qquad
\chi(\xi)\,|\xi|^{-2}\,:\,\,L^{\frac {2n}{n-3}+}\to L^{\frac{2n}{n+1}-},
$$  
Since $\frac {n-2}{2n} + \frac 2n=\frac {n+2}{2n}$ and 
$\frac {n-3}{2n} + \frac 2n= \frac {n+1}{2n}$, we have
$$
\|\,\chi(\xi) |\xi|^{-2} f\|_{ L^{\frac {2n}{n+2}-}}\le \|\,\,|\xi|^{-2}\|_{L^{\frac n2-}(B)}
\|f\|_{\frac {2n}{n-2}+}\le C \|f\|_{\frac {2n}{n-2}+}
$$
and
$$
\|\,\chi(\xi) |\xi|^{-2} f\|_{L^{\frac{2n}{n+1}-}}\le 
 \|\,\,|\xi|^{-2}\|_{L^{\frac n2-}(B)} \|f\|_{L^{\frac {2n}{n-3}}+}\le 
C \|f\|_{L^{\frac {2n}{n-3}}+}
$$
as desired.
\end{proof}

\begin{proof}[Proof of Lemma \ref{le:Cont}]
The result of Lemma \ref{le:Cont} is contained in Ginibre-Moulin 
\cite{GM} and can be traced to the earlier work of Kato \cite{kato}.
Here we essentially reproduce the proof in \cite{GM}.
 
\noindent Consider the resolvent $R(\la+i\b)$. We shall 
prove that it is continuous (in fact, H\"older  continuous)
in the upper half-plane $\overline{\mathbb C}_+$. The same statement
also holds for $R(\la-i\b)$.
We appeal to the representation \eqref{eq:respro},
$$
R_0(\la+i\b)= \il_0^\infty e^{i(\la+i\b)t} e^{itH_0}\,dt.
$$
Therefore, using the inequality 
$|e^{i(\la_2+i\b_2)t} -  e^{i(\la_1+i\b_1)t}| \le 
\min (2, (|\la_2-\la_1| + |\b_2-\b_1|)t)$ and the embedding $L^{\frac{2n}{n-2}+}(\R^n)\einbet L^{2,-1-\gamma}(\R^n)$,  
we obtain for arbitrary $\la_1, \la_2 \in \R$ and  $\b_1, \b_2 \in \R_+$,
\be
\|(R_0(\la_2+i\b_2) - R(\la_1 +i\b_1))f\|_{L^{2,-1-\ga}} 
&&\le \il_0^\infty  \| e^{itH_0} f\|_{L^{2,-1-\ga}} \,
\min (2, (|\la_2-\la_1| + |\b_2-\b_1|)t)\,dt \nn \\
&&\le \il_0^\infty  \| e^{itH_0} f\|_{L^{\frac {2n}{n-2}+}} \,
\min (2, (|\la_2-\la_1| + |\b_2-\b_1|)t)\,dt.  \label{eq:ginibre}
\ee 
We now recall that in addition to the Strichartz estimates \eqref{eq:striH0}
the semigroup $e^{it H_0}$ also verifies the dispersive estimates
$$
\|e^{it H_0} f\|_{L^p} \le \frac {C_p}{t^{n(\frac 12-\frac 1p)}}
\|f\|_{L^{p'}},\qquad \forall p\in [2,\infty].
$$
Inserting this bound into~\eqref{eq:ginibre} and invoking the embedding 
$L^{2,1+\gamma} \einbet L^{\frac {2n}{n+2}-}$, 
we infer that
$$
\aligned
\|(R_0(\la_2+i\b_2) - R(\la_1 +i\b_1))f\|_{L^{2,-1-\ga}} 
&\lesssim \il_0^M  \frac {(|\la_2-\la_1| + |\b_2-\b_1|)t)}{t^{1+}}\,
dt\, \|f\|_{L^{\frac {2n}{n+2}-}}  + \il_M^\infty  \frac 2{t^{1+}}\,\|f\|_{L^{\frac {2n}{n+2}-}}
\\ &\lesssim   \bigg((|\la_2-\la_1| + |\b_2-\b_1|) M^{1-} + 2 M^{0-}\bigg)
\|f\|_{L^{2,1+\ga}}
\endaligned
$$
for some constant $M>0$. Choosing $M= (|\la_2-\la_1| + |\b_2-\b_1|)^{-1}$
we finally conclude that
$$
\|(R_0(\la_2+i\b_2) - R(\la_1 +i\b_1))f\|_{L^{2,-1-\ga}} \lesssim
(|\la_2-\la_1| + |\b_2-\b_1|)^{0+}   \|f\|_{L^{2,1+\ga}},
$$
as claimed. 
\end{proof}

\begin{proof}[Proof of Lemma \ref{le:Vinv}]
Lemma \ref{le:LAP} implies the boundedness of $G: L^{2,1+\ga} \to L^{2,-1-\ga}$ for 
any positive $\ga$. Therefore, since  $|V(x)|\le C(1+|x|^2)^{-1-\eps}$, the potential~$V$
maps the space $L^{2,-1-\ga}$ into $L^{2, 1-\ga +2\eps}$. 
Thus, using the second conclusion of Lemma \ref{le:LAP}, 
we obtain that 
\begin{equation}
\label{eq:extra}
G V :  L^{2,-1-\ga} \to L^{2,-1-\ga + 2\eps}.
\end{equation} 
provided that $0<-\ga + 2\eps <\frac 12$.

Since $-\Lap (G V)= V$ and $V$, of course, maps  $L^{2,-1-\ga}\to L^{2,-1-\ga}$,
we conclude that $G V$ takes the space $L^{2,-1-\ga}$ to $H^2_{\text loc}$.
The compactness 
of $G V$ on $ L^{2,-1-\ga}$ then follows from the observation above and 
 the extra $2\eps$ decay at infinity
established in \eqref{eq:extra}.  We infer that $I + R_0 V$ is a Fredholm operator
on  $L^{2,-1-\ga}$ and it is therefore  invertible iff its null space is empty.

\noindent 
Let $\phi$ be an  $ L^{2,-1-\ga}$  solution of the equation 
\begin{equation}
\label{eq:nullspace}
\phi + G V \phi =0
\end{equation}
First we observe that by \eqref{eq:extra} function $\phi$, in fact, belongs 
to the space $L^{2,-1-\ga + 2\eps}$ for some $\ga:\,0<\ga<2\eps$. 
It then follows that $V\phi \in L^{2, 1-\ga + 4\eps}$.
Lemma \ref{le:LAP} then implies that, as long 
as $4\eps<\frac 12$, also $G V\phi \in  L^{2,-1-\ga + 4\eps}$ and, using~\eqref{eq:nullspace} again 
we have that $\phi \in  L^{2,-1-\ga + 4\eps}$.
We can continue this argument and obtain that $\phi\in  L^{2,-\frac 12-\a}$ for 
any $\a>0$. 

\noindent  Applying $-\Lap$ to both sides of the equation \eqref{eq:nullspace}
we conclude that $\phi$ is an $L^{2,-\frac 12-\a}$ solution of the equation
$$
-\Lap \phi + V\phi =0,
$$
and thus either an eigenfunction or a resonance corresponding to $\la=0$.
Since we assumed that $\la=0$ is a regular point, 
$\phi$ must be identically zero and the null space of $I+R_0V$ is empty.
\end{proof}

\section{Time dependent potentials: Reduction to oscillatory integrals}\label{sec:time1}

\begin{defi}
\label{def:Ydef}
Let $Y$ be the normed space of measurable functions $V(t,x)$ on~$\R^3$  that satisfy the following 
properties: $t\mapsto \|V(t,\cdot)\|_{L^{\frac32}(\R^3)} \in L^\infty(\R)$ and
for a.e.~$x\in\R^3$ the function
$t\mapsto V(t,x)$ is in $\mathcal S'(\R)$, the space of tempered distributions. Moreover, the Fourier
transform of this distribution, which we denote by~$V(\hat{\tau},x)$, is a (complex) measure whose 
norm satisfies
\begin{equation}
\label{eq:measnorm}
\sup_{y\in\R^3} \int_{\R^3} \frac{\|V(\hat{\tau},x)\|_{\mathcal M}}{|x-y|}\,dx < \infty.
\end{equation}
The norm in $Y$ is the sum of the expression on the left-hand side of~\eqref{eq:measnorm} and the norm
in~$L_t^\infty(L_x^{\frac32})$. 
\end{defi}

\noindent In what follows we study the Schr\"odinger equation
\be
&&i\pr_t \psi - \Lap \psi + V(t,x) \psi=0,\label{eq:schrt}\\
&&\psi|_{t=s} (x) =\psi_s(x)\nonumber
\ee
for potentials $V\in Y$ and with initial data $\psi_s\in L^2(\R^3)$. 
An interesting case is $V(t,x)=\cos(t)\,V(x)$ where~$V\in L^{\frac32}$ satisfies 
$\sup_y \int_{\R^3} \frac{|V(x)|}{|x-y|}\,dx<\infty$. 
Because of the limited regularity
of potentials in~$Y$, we define (weak) solutions $U(t,s)\psi_s$ of~\eqref{eq:schrt} via
Duhamel's formula:
\begin{equation}
\label{eq:duhamel}
U(t,s) \psi_s = 
e^{i(t-s)H_0}\psi_s + i \int_s^t e^{i(t-s_1)H_0} V(s_1,\cdot) U(s_1,s) \psi_s \,ds_1.
\end{equation}
In the following lemma we show by means of  Keel's and Tao's endpoint Strichartz estimate~\cite{KT} 
that such weak solutions exist and are unique provided the potential is small in an appropriate
sense. The proof is presented only in~$\R^3$, but it carries over to any dimension~$n\ge3$. 
We set 
\begin{equation}
\label{eq:Xdef}
X=L_t^\infty(L_x^{2}(\R^3))\cap L_t^2(L_x^6(\R^3))
\end{equation} 
and  define $H_0=-\Lap$ to be the unperturbed Schr\"odinger operator with evolution~$e^{itH_0}$.

\begin{lemma}
\label{lem:weaksol}
Assume that the potential $V(t,x)$ satisfies the smallness assumption
\begin{equation}
\label{eq:small}
\|V\|_{L^\infty_t L^{\frac 32}_x} = \sup_{t\in {\R}} 
\bigg(\int_{\R^3} |V(t,x)|^{\frac 32} \bigg)^{\frac 23} < c_0
\end{equation}
for some sufficiently small constant $c_0>0$. Then for any
$s\in {\mathbb R}$ and any $\psi_s\in L^2_x$ there exists a unique weak 
solution~$U(t,s)\psi_s$ of~\eqref{eq:duhamel} with the property that $U(\cdot,s)\psi_s \in X$
and so that $t\mapsto \langle U(t,s)\psi_s, g\rangle$ is continuous for any~$g\in L^2(\R^3)$.
Moreover, for any such~$g$ and any $t>s$, 
\begin{align}
\langle U(t,s)\psi_s,g \rangle  = \sum_{m=0}^\infty \quad i^m \idotsint\limits_{s\le s_m\le..\le s_1\le t}
&\langle e^{i(t-s_1)H_0} V(s_1,\cdot) e^{i(s_1-s_2) H_0} V(s_2,\cdot)\ldots\nn\\
&\qquad \qquad V(s_m,\cdot) e^{i(s_m-s)H_0}\, \psi_s, g \rangle\; ds_1\ldots ds_m \label{eq:series}
\end{align}
where the series converges absolutely. In the strong sense, i.e., without the pairing against~$g$, this
representation holds in the sense of norm convergence in the space~$X$ (and thus can only be assumed
for~a.e.~$t$). 
\end{lemma}
\begin{proof} For the purposes of this proof, we let $F=F(t,x)$ be a function of $(t,x)\in \R^{1+3}_{t,x}$.
For simplicity, we often write $F(t)$ for the function $x\mapsto F(t,x)$.
Recall the following end-point Strichartz estimates for the operator $H_0$ proved by Keel-Tao in any dimension $n\ge 3$, see~\cite{KT}: There exists some dimensional constant $C_1=C_1(n)$ so that for all 
$f\in L^2_x$ and $F\in L^2_t L_x^{\frac {2n}{n+2}}$ 
\begin{align}
&\|e^{it H_0} f \|_{L^2_t L^{\frac {2n}{n-2}}_x} \le C_1 \|f\|_{L^2_x},
\label{eq:homstr}\\
&\|\int_s^t e^{i(t-s_1)H_0} F(s_1)\,ds_1\|_{L^2_t L_x^{\frac {2n}{n-2}}} \le C_1
\|F\|_{L^2_t L_x^{\frac {2n}{n+2}}}.\label{eq:inhomstr}
\end{align}
Consider the operator ${\mathcal K}_s$ defined by
$$
({\mathcal K}_s F) (t,\cdot) = i \int_s^t e^{i(t-s_1)H_0} V(s_1,\cdot) F(s_1,\cdot)\,ds_1.
$$
Then definition~\eqref{eq:duhamel} takes the form
\begin{equation}
\label{eq:defagain} 
\bigl[(1-{\mathcal K}_s)U(\cdot,s)\psi_s\bigr] (t) = e^{i(t-s)H_0}\psi_s.
\end{equation}
Inequality~\eqref{eq:inhomstr} and the smallness assumption~\eqref{eq:small} imply that 
 the norm of the operator ${\mathcal K}_s: \,L^2_t L_x^{6} \to L^2_t L_x^{6}$ satisfies
\begin{equation}
\label{eq:Kmixed}
\|{\mathcal K}_s F \|_{L^2_t L_x^{6}}\le C_1\,
\|V\, F\|_{L^2_t L_x^{\frac {6}{5}}}\le C_1\, \|V\|_{L^\infty_t L_x^{\frac 32}}
\|F\|_{L^2_t L_x^{6}}\le C_1\,c_0  \|F\|_{L^2_t L_x^{6}}.  
\end{equation}
Moreover, for any $g\in L^2(\R^3)$, 
\be
|\langle ({\mathcal K}_s F )(t) , g \rangle | &=& \left | \int_s^t \langle V(s_1,\cdot) F(s_1,\cdot), 
e^{-i(t-s_1)H_0} g \rangle \,ds_1 \right| \nn \\
&\le& \int_s^t \| V(s_1,\cdot) F(s_1,\cdot)\|_{\frac65} \; \| e^{-i(t-s_1)H_0} g \|_{6} \,ds_1 \nn\\ 
&\le& \Bigl(\int_s^t \| V(s_1,\cdot) F(s_1,\cdot)\|^2_{\frac65}\,ds_1\Bigr)^{\half} \; \Bigl(\int_s^t \| e^{-i(t-s_1)H_0} g \|_{6}^2 \,ds_1 \Bigr)^{\half}  \nn \\
&\le& C_1\,\Bigl(\int_s^t \| V \|^2_{L_t^\infty(L_x^{\frac32})} \| F(s_1,\cdot)\|^2_{6}\,ds_1\Bigr)^{\half}\; \|g\|_2 \label{eq:ktagain} \\
&=& C_1\,\| V \|_{L_t^\infty(L_x^{\frac32})} \|F\|_{L_t^2(L_x^6)} \|g\|_2 \nn
\ee
where we used~\eqref{eq:homstr} to pass to~\eqref{eq:ktagain}. This shows that
\[ \esssup_{t}\| ({\mathcal K}_s F)(t) \|_2 \le C_1\, c_0 \|F\|_{L_t^2(L_x^6)}\] 
which in conjunction with~\eqref{eq:Kmixed} yields that
\begin{equation}
\label{eq:contract}
 \| {\mathcal K}_s \|_{X\to X} \le C_1\,c_0 < \frac12,
\end{equation}
provided $c_0$ is small (see~\eqref{eq:Xdef} for the definition of~$X$).
Therefore, the operator $I-{\mathcal K}_s$ is invertible on the space~$X$ and~$U(t,s)$ 
can be expressed via the Neumann series
$$
U(t,s) = \bigl[ (I-{\mathcal K}_s)^{-1}  e^{i(\cdot-s)H_0} \bigr] (t) =
\sum_{m=0}^\infty \bigl[ {\mathcal K}_s^m  e^{i(\cdot -s)H_0} \bigr] (t)
$$
which converges in the norm of $X$. Writing out $\langle U(t,s)\psi_s,g\rangle$ explicitly 
leads to~\eqref{eq:series}.
Next we check that for any $F\in L^2_t(L^6_x)$ the function $t\mapsto \langle \Kop F , g \rangle$ is
continuous for any choice of~$g\in L^2$. In fact, if $t_1<t_2$, then
\be
| \langle (\Kop F)(t_2), g \rangle  - \langle (\Kop F)(t_1), g \rangle  |
&\le& \int_s^{t_2} |\langle V(s_1)F(s_1), (e^{-i(t_1-s_1)H_0}- e^{-i(t_2-s_1)H_0}) g \rangle |\, ds_1 \nn \\
&& \qquad\qquad + \int_{t_1}^{t_2} |\langle V(s_1)F(s_1), e^{-i(t_1-s_1)H_0} g \rangle |\,ds_1 \nn \\
&\le& \|V\|_{L_t^\infty(L^{\frac32}_x)}\|F\|_{L_t^2(L_x^6)} \|g-e^{-i(t_2-t_1)H_0} g\|_2 \nn \\
&& \qquad\qquad + \|V\|_{L_t^\infty(L^{\frac32}_x)}\Bigl( \int_{t_1}^{t_2}\|F(s_1)\|_{L_x^6}^2\,ds_1\Bigr)^{\half} \;\|g\|_2.
\nn
\ee
Since the last expression tends to zero as $t_2\to t_1$, continuity follows. Hence $(\Kop^m F)(t)$ is also weakly continuous in~$t$, and thus $\Kop^m e^{i(\cdot-s)H_0} \psi_s$ is, too. Since $\langle U(t,s)\psi_s,g\rangle$ is a uniformly convergent series of these continuous functions, it follows that it is continuous.
\end{proof}
\begin{remark}
\label{re:weakc}
The proof of Lemma~\ref{lem:weaksol} shows that the operator 
${\mathcal K}_s: L^2_t L^6_x\to w-C^0_t(L^2_x)$ maps~$L^2_t L^6_x$ 
into the space of weakly continuous functions
with values in $L^2({\mathbb R}^3)$.
\end{remark}

\noindent For technical reasons connected with the functional calculus in the following section it will be
convenient to work with smooth potentials in~$Y$ rather than general ones. To approximate
a general potential~$V$ by means of smooth ones, choose nonnegative cut-off functions $\chi\in \calS(\R^3)$ and $\eta\in\calS(\R)$ so that $\chi$ and $\hat{\eta}$ have compact support and satisfy
$\int_{\R^3} \chi(x)\,dx=1$, $\int_{\R}\eta(t)\,dt=1$. In addition, let $\chi=1$ on a neighborhood of~$0$. 
For any $V\in Y$ and $R>1$ define
\[ V^{(1)}_R(t,\cdot) := V(t,\cdot)\chi(\frac{.}{R}) \ast R^3\chi(R\cdot)\]
where the convolution is in the $x$-variable only.
Note that $V^{(1)}_R$ is well-defined, smooth and compactly supported in~$x$,  and satisfies $\|V^{(1)}_R\|_{L^\infty_{t,x}}<\infty$ since~$\|V\|_{L^\infty_t(L^{\frac32}_x)}<\infty$. Moreover, it is standard to check that
\[ \sup_{R>0}\|V^{(1)}_R\|_Y \le \|\chi\|_\infty\, \|V\|_Y. \]
Indeed, 
\[ \|V^{(1)}_R(t,\cdot) \|_{L^\frac{3}{2}} \le \|\chi\|_\infty 
\|\, |V(t,\cdot)| \ast R^3\chi(R\cdot) \|_{L^{\frac{3}{2}}} 
\le \|\chi\|_\infty\, \|V(t,\cdot)\|_{L^\frac{3}{2}},\]
whereas with $\Gamma(x):= |x|^{-1}$ and $\meas$ denoting measures in the $\tau$-variable,
\be 
\Bigl(\|V^{(1)}_R(\hat{\tau},\cdot)\|_{\meas}\ast \Gamma\Bigr)(x)  &\le& 
   \sup_{y} \Bigl(\chi(\frac{.}{R})\,\|V(\hat{\tau},\cdot)\|_{\meas}\ast \Gamma\Bigr)(y) \nn \\
&\le& \|\chi\|_\infty\,\sup_{y} \Bigl(\|V(\hat{\tau},\cdot)\|_{\meas} \ast \Gamma\Bigr)(y), \nn
\ee
as claimed. To regularize in~$t$, define
\[ V_R(\cdot,x) := [V^{(1)}_R(\cdot,x) \ast R\eta(R\cdot)] \eta(\frac{.}{R})\]
where the convolution is in the $t$-variable only.
Again one checks that
\[ \|V_R\|_Y \le (\|\eta\|_\infty + \|\hat{\eta}\|_1)\, \|V^{(1)}_R\|_Y\le (\|\eta\|_\infty+\|\hat{\eta}\|_1)\,\,\|\chi\|_\infty\, \|V\|_Y \]
for any $R>0$. We will use that $V_R\to V$ as $R\to\infty$ in the following sense:
For a.e. $t$ one has
\begin{equation}
\label{eq:Vconv} 
\|V_R(t,\cdot)-V(t,\cdot)\|_{L^{\frac32}(\R^3)} \to 0 \text{\ \ as $R\to \infty$.}
\end{equation}
Firstly, it follows from standard measure theory that for a.e. $t$ 
\begin{equation}
\label{eq:V1conv}
\|V_R^{(1)}(t,\cdot)-V(t,\cdot)\|_{L^{\frac32}(\R^3)} \to 0 \text{\ \ as \ \ }R\to \infty.
\end{equation}
Secondly, with $\eta_R(t):=R\eta(Rt)$,
\be
&& \|V_R(t,\cdot) - V_R^{(1)}(t,\cdot) \|_{\frac32} \nn \\
 &\le& \| (V_R^{(1)} \ast \eta_R)(t,\cdot) - V_R^{(1)}(t,\cdot) \|_{\frac32} + \|V_R^{(1)}(t,\cdot)\|_{\frac32} |1-\eta(t/R)|    \nn\\
&\le& \int_{-\infty}^\infty \eta_R(s)\, \|V_R^{(1)}(t-s,\cdot) - V_R^{(1)}(t,\cdot) \|_{\frac32} \, ds
  + |1-\eta(t/R)|\, \|\eta\|_\infty \|V\|_{L_t^\infty(L_x^{\frac32})}   \nn \\
&\le& \|\chi\|_\infty\,\int_{-\infty}^\infty \eta_R(s)\, \|V(t-s,\cdot) - V(t,\cdot) \|_{\frac32}\, ds  + o(1) \label{eq:leb} \\
&\to& 0 \nn
\ee
for a.e. $t$ as $R\to \infty$. The conclusion~\eqref{eq:leb} follows from the vector-valued analogue
of the Lebesgue differentiation theorem (in this case ``vector-valued'' means with values in~$L^{\frac32}$).  
In combination with~\eqref{eq:V1conv} this yields~\eqref{eq:Vconv}.

\noindent We shall now prove the convergence of the approximate solutions $\psi_R(t,x)$ satisfying 
the equation 
\be
&&i\,\partial_t \psi_R - \Delta \psi_R + V_R(t,x)\psi_R = 0,\label{eq:psiVR}\\
&&\psi_R|_{t=s}=\psi_s\nn
\ee
to the solution $\psi(t,x)$ of the original problem corresponding to the 
potential $V(t,x)$. Note that due to the smoothness and boundedness of
the potentials $V_R$ the $L^\infty_t L^2_x$ function $\psi_R$ can be interpreted as a
distributional solution of equation~\eqref{eq:psiVR}. In fact, the left 
hand-side of~\eqref{eq:psiVR} belongs to the space~$L^\infty_t H^{-2}$. In addition,
$\psi_R$ is also a Duhamel solution as in~\eqref{eq:duhamel}.

\begin{lemma}
\label{Le:Conv}
Let $U_R(t,s)$ be the propagator~\eqref{eq:psiVR},
i.e.,~$U_R(t,s)\psi_s = \psi_R(t,s)$. Then for any $s,t\in \mathbb R$ such that
$t\ge s$, and arbitrary functions $\psi_s, g\in L^2({\mathbb R}^3)$,
$\|\psi_s\|_{L^2_x}=\|g\|_{L^2_x}=1$  we have
\begin{equation}
\label{eq:conv}
<U_R(t,s)\psi_s, g>\;\to\; <U(t,s)\psi_s, g>\qquad \mbox{as}\quad R\to \infty
\end{equation}
\end{lemma}
\begin{proof}
First observe that since the potential $V$ satisfies the smallness assumption 
\eqref{eq:small}, $V_R$ also obeys~\eqref{eq:small} for all $R>0$. 
According to Lemma~\ref{lem:weaksol}, 
\begin{align}
\langle U_R(t,s)\psi_s,g \rangle  = \sum_{m=0}^\infty \quad i^m 
\idotsint\limits_{s\le s_m\le..\le s_1\le t}
&\langle e^{i(t-s_1)H_0} V_R(s_1,\cdot) e^{i(s_1-s_2) H_0} V_R(s_2,\cdot)....\nn\\
&\qquad \qquad V_R(s_m,\cdot) e^{i(s_m-s)H_0}\, \psi_s, g \rangle\; ds_1...ds_m 
\label{eq:seriesR}
\end{align}
for any  $\psi_s, g\in L^2({\mathbb R}^3)$. Equivalently, 
$U_R(t,s)$ can be represented by the Neumann series
$$
U_R(t,s) = \bigl[ (I-{{\mathcal K}_R}_s)^{-1}  e^{i(\cdot-s)H_0} \bigr] (t) =
\sum_{m=0}^\infty \bigl[ {{\mathcal K}_R}_s^m  e^{i(\cdot -s)H_0} \bigr] (t)
$$
which converges in the norm of the space $X$ defined above.
The operators $ {{\mathcal K}_R}_s: X\to X$ are defined as
$$
({{\mathcal K}_R}_s F) (t,\cdot) = i \int_s^t e^{i(t-s_1)H_0} V_R(s_1,\cdot) F(s_1,\cdot)\,ds_1.
$$ 
To verify the conclusion of Lemma~\ref{Le:Conv} it suffices to show
that for an arbitrary positive~$\eps>0$, all positive integers $m\le m_0(\eps)$, 
and all sufficiently large $R=R_0(\eps,m_0)$ 
\begin{equation}
\label{eq:convK} 
|\langle \big({\mathcal K}_s^m - {{\mathcal K}_R}_s^m\big)e^{i(\cdot -s)H_0}\psi_s,g
\rangle (t)|
< \eps (C_1 c_0)^{m-1} m
\end{equation}
The positive integer $m_0(\eps)$ is chosen so that
$2(C_1 c_0)^m\le \eps$ which ensures the smallness of the ``tails'' 
of the series for $U(s,t)$ and $U_R(s,t)$. 

\noindent For the bounded operators ${\mathcal K}_s,  {{\mathcal K}_R}_s$ on the space
$X$ we have the following identity: 
\begin{equation}
\label{eq:Trotter}
{{\mathcal K}_R}^m_s - {{\mathcal K}}_s^m =\sum_{\ell=0}^{m-1} {{\mathcal K}_R}_s^{\ell} 
({\mathcal K}_s - {{\mathcal K}_R}_s) {{\mathcal K}}_s^{m-\ell-1} 
\end{equation}
We shall prove that for $\ell\in [0,m-1]$
\begin{equation}
\label{eq:convKprod}
|\langle  {{\mathcal K}_R}_s^{\ell} 
({\mathcal K}_s - {{\mathcal K}_R}_s) {{\mathcal K}}_s^{m-\ell-1}
 e^{i(\cdot -s)H_0}\psi_s,g\rangle (t)|
< \eps (C_1 c_0)^{m-1} 
\end{equation}
which immediately implies~\eqref{eq:convK}.

\noindent In view of Remark~\ref{re:weakc} the operator 
${\mathcal K}_s$, and thus also
${{\mathcal K}_R}_s$, maps $L^2_t L^6_x\to w-C^0_t (L^2_x)$.
Therefore, for an arbitrary fixed $t\ge s$ we can define the operator 
${{\mathcal K}_R}_{s,t}: L^2_t L^6_x \to L^2_x$ via the formula 
$$
{{\mathcal K}_R}_{s,t} F = ({\mathcal K}_s F)(t)
$$
In addition to the $L^2_x$ pairing $\langle,\rangle$ we define 
the space-time pairing $\langle,\rangle_{t,x}$ as usual: 
for any pair of functions $F\in L^q_t L^p_x$ and $G\in L^{q'}_t L^{p'}_x$
with $q,p\in [1,\infty]$ let 
$$
\langle F,G\rangle_{t,x} = \int_{\mathbb R} \int_{{\mathbb R}^3} 
F(t,x) \overline {G(t,x)}\,dx\,dt
$$ 
We now introduce the dual operator ${{\mathcal K}_R}^*_{s,t}: L^2_x\to L^2_t L^{\frac 65}_x$.
In addition, since ${{\mathcal K}_R}_s: L^2_t L^6_x\to L^2_t L^6_x$ we also define the 
dual of ${{\mathcal K}_R}_s$, ${{\mathcal K}_R}_s^*: L^2_t L^{\frac 65}_x\to 
L^2_t L^{\frac 65}_x$.
Therefore for $\ell\ge 1$ the left hand-side 
of~\eqref{eq:convKprod} can be written as
\be 
I_{R,s,t} &:=& \langle {{\mathcal K}_R}_{s,t} {{\mathcal K}_R}_s^{\ell-1} 
({\mathcal K}_s - {{\mathcal K}_R}_s) {{\mathcal K}}_s^{m-\ell-1}
 e^{i(\cdot -s)H_0}\psi_s,g\rangle \nn \\
&=& 
\langle 
({\mathcal K}_s - {{\mathcal K}_R}_s) {{\mathcal K}}_s^{m-\ell-1}
 e^{i(\cdot -s)H_0}\psi_s, {{{\mathcal K}_R}^*_s}^{\ell-1} {{\mathcal K}_R}^*_{s,t}g\rangle_{t,x}. \nn
\ee
We assume that $\eps$ and $m_0(\eps)$ are now fixed and invoke Egorov's theorem. 
According to~\eqref{eq:Vconv}
$$
\|V_R(s_1,\cdot)-V(s_1,\cdot)\|_{L^{\frac 32}_x} \to 0,\qquad
\mbox{as} \quad R\to \infty
$$
for a.e. $s_1\in [s,t]$.
Therefore, for any $\delta >0$ there exists a set 
${\mathcal B}\subset [s,t]$ such that $|\mathcal B|\le \delta$ and 
$$
\|V_R(s_1,\cdot)-V(s_1,\cdot)\|_{L^{\frac 32}_x} <\eps,\qquad 
\forall s_1\in [s,t]\setminus {\mathcal B}
$$
and all sufficiently large $R\ge R_0(\eps,\delta)$.
Let $\chi_{\mathcal B}$ be the characteristic function of the set ${\mathcal B}$.
We define operators 
\be
&&{\mathcal Y}_s = ({\mathcal K}_s - {{\mathcal K}_R}_s)\chi_{\mathcal B},\nn\\
&&{\mathcal Z}_s = ({\mathcal K}_s - {{\mathcal K}_R}_s)(1-\chi_{\mathcal B}),\nn 
\ee
It is easy to see that ${\mathcal Y}_s,  {\mathcal Z}_s: L^2_t L^6_x \to L^2_t L^6_x$.
Moreover, $\| {\mathcal Z}_s\|_{L^2_t L^6_x \to L^2_t L^6_x}\le C\eps$ for all 
 $R\ge R_0(\eps,\delta)$
and 
$\| {\mathcal Y}_s\|_{L^2_t L^6_x \to L^2_t L^6_x}\le C_1c_0<\half$, see~\eqref{eq:contract}.
Therefore,
$$
I_{R,s,t}=\langle {{\mathcal K}}_s^{m-\ell-1}
 e^{i(\cdot -s)H_0}\psi_s, {\mathcal Z}_s^* {{{\mathcal K}_R}^*_s}^{\ell-1} 
{{\mathcal K}_R}^*_{s,t}g\rangle_{t,x} + \langle \chi_{\mathcal B} {{\mathcal K}}_s^{m-\ell-1}
 e^{i(\cdot -s)H_0}\psi_s, {\mathcal Y}_s^* {{{\mathcal K}_R}^*_s}^{\ell-1} 
{{\mathcal K}_R}^*_{s,t}g\rangle_{t,x} 
$$
We can easily estimate the first term by 
$$
\|{{\mathcal K}}_s\|^{m-\ell-1}_{L^2_t L^6_x \to L^2_t L^6_x}
\|{{\mathcal Z}}^*_s\|_{L^2_t L^{\frac 65}_x \to L^2_t L^{\frac 65}_x}
\|{{\mathcal K}}^*_s\|^{\ell-1}_{L^2_t L^{\frac 65}_x \to L^2_t L^{\frac 65}_x}
\|{{\mathcal K}}^*_{s,t}\|_{L^2_x \to L^2_t L^{\frac 65}_x} \|\psi_s\|_{L^2_x}
\|g\|_{L^2_x}\le  \eps 2^{-(m-1)}. 
$$
For the second term we have the bound
$$
\aligned
\|\chi_{\mathcal B}{{\mathcal K}}_s^{m-\ell-1} e^{i(\cdot -s)H_0}\psi_s\|_
{L^2_t L^6_x}
\|{{\mathcal Y}}^*_s\|_{L^2_t L^{\frac 65}_x \to L^2_t L^{\frac 65}_x}
&\|{{\mathcal K}}^*_s\|^{\ell-1}_{L^2_t L^{\frac 65}_x \to L^2_t L^{\frac 65}_x}
\|{{\mathcal K}}^*_{s,t}\|_{L^2_x \to L^2_t L^{\frac 65}_x} \|\psi_s\|_{L^2_x}
\|g\|_{L^2_x}\\ &
\le \|\chi_{\mathcal B}{{\mathcal K}}_s^{m-\ell-1} e^{i(\cdot -s)H_0}\psi_s\|_
{L^2_t L^6_x} (C_1 c_0)^{\ell+1}  
\endaligned
$$
Observe that 
$$ 
\|{{\mathcal K}}_s^{m-\ell-1} e^{i(\cdot -s)H_0}\psi_s\|_
{L^2_t L^6_x}\le (C_1 c_0)^{m-\ell-1}\|\psi_s\|_{L^2_x}<\infty
$$
Therefore, we can chose $\delta=\delta(m_0)$ in Egorov's theorem in 
such a way that 
$$
\sum_{m=1}^{m_0(\eps)} \sum_{\ell=1}^{m-1} 
\|\chi_{\mathcal B} {{\mathcal K}}_s^{m-\ell-1} e^{i(\cdot -s)H_0}\psi_s\|_
{L^2_t L^6_x}\le \eps (C_1 c_0)^{m_0-\ell-1}
$$
Hence we have the desired bound 
$$
|I_{R,s,t}|\le \eps (C_1 c_0)^{m-1}
$$
for all $1\le \ell\le m-1$ and $m\le m_0$. 
To settle the remaining case of $\ell=0$ we observe that for $\ell=0$
$$
I_{R,s,t} = \langle ({\mathcal K}_s - {{\mathcal K}_R}_s) {{\mathcal K}}_s^{m-1}
 e^{i(\cdot -s)H_0}\psi_s,g\rangle (t) =
 \langle {\mathcal Y}_s {{\mathcal K}}_s^{m-1}
 e^{i(\cdot -s)H_0}\psi_s,g\rangle (t) +\langle {\mathcal Z}_s {{\mathcal K}}_s^{m-1}
 e^{i(\cdot -s)H_0}\psi_s,g\rangle (t) 
$$
Similarly to the operator ${{\mathcal K}_R}_{s,t}$ we can define 
the operators ${\mathcal Y}_{s,t},  {\mathcal Z}_{s,t}: L^2_t L^6_x\to L^2_x$.
Moreover, 
$$
\| {\mathcal Z}_{s,t}\|_{L^2_t L^6_x\to L^2_x}\le C\eps,\qquad
\| {\mathcal Y}_{s,t}\|_{L^2_t L^6_x\to L^2_x}\le C_1 c_0.
$$
Thus
$$
|I_{R,s,t}|\le C_1 c_0\,\|\chi_{\mathcal B}{{\mathcal K}}_s^{m-1}
 e^{i(\cdot -s)H_0}\|_{L^2_t L^6_x}  + 
\eps (C_1 c_0)^{m-1}\le  2\eps (C_1 c_0)^{m-1}
$$
by the choice of the constant $\delta$ in Egorov's theorem.
\end{proof}
Since the potentials $V_R(t,x)$ are smooth in both variables,
the solution operators $U_R(t,s)$ are unitary on $L^2_x$. Together 
with Lemma~\ref{Le:Conv} we have the following 
\begin{cor}
\label{co:conserv}
The $L^2$ norm of the solution $\psi(t,\cdot)$ of the Schr\"odinger
equation~\eqref{eq:schrt} is a non-increasing function of time,
i.e,
$$
\|U(t,s)\psi_s\|_{L^2_x}\le \|\psi_s\|_{L^2_x}
$$
for all $t\ge s$ and arbitrary functions $\psi_s\in L^2_x$.
\end{cor}
Lemma~\ref{Le:Conv} also implies that we can assume henceforth
that $V(t,x)$ is a smooth potential with compact support in the 
$x$-variable and the variable $\hat\tau$ of the Fourier transform relative to $t$.
We can also assume that $V$ satisfies the smallness assumption 
\eqref{eq:small}.
We shall show that the following estimates depend only on the 
norm of the potential in the space $Y$ defined in Defintion~\ref{def:Ydef} and the smallness 
constant~$c_0$.
 

\subsection{Functional calculus}
The goal of this section is to obtain the explicit representation of the integral kernels of the 
operators involved in the Neumann series expansion \eqref{eq:series} for $U(t,s)$, as some special
oscillatory integrals.

We introduce the notation 
$$
V(\hat\tau,\cdot):=\int e^{it\tau} V(t,\cdot)\,dt.
$$
The $m$-th term of the series~\eqref{eq:series}, which we denote by~${\mathcal I}_m$,
can then be written in the following form\footnote{Here we use the fact that $V(\hat\tau,\cdot)$ 
has compact support in $\hat\tau$ to interchange the integrals.}:
\begin{align}
\langle {\mathcal I}_m(t,s)\psi_s,g\rangle =\int_{{\mathbb R}^m} d\tau_1..d\tau_m 
\idotsint\limits_{s\le s_m\le..\le s_1\le t}ds_1...ds_m\;
&\langle e^{i(t-s_1)H_0} e^{is_1\tau_1} V(\hat\tau_1,\cdot) ...\label{funcIm}\\
&\cdot e^{is_m\tau_m} V(\hat\tau_m,\cdot) e^{i(s_m-s)H_0}\psi_s,\,g\rangle. \nn
\end{align}
The identity above is verified on arbitrary functions $\psi_s, g\in L^2_x$.

\noindent We shall also make use of the spectral representation of the 
operator $e^{itH_0}$,
$$
e^{it H_0} = \int_{\mathbb R} e^{it\la} dE(\la).
$$
Here, $dE(\la)$ is the spectral measure associated with the 
operator $H_0=-\Lap$. In dimension $n=3$, $dE(\la)$ 
has an explicit representation as an integral operator with
the kernel
$$
dE(\la)(x,y) = \begin{cases} 
\frac{\sin \sqrt{\la}|x-y|}{4\pi|x-y|}\,d\la & \la>0,\\
0& \la\le 0
\end{cases}
$$ 
Recall also that the resolvent $R(z)= (H_0-z)^{-1}$ is
an analytic function with values in the space of bounded operators 
in $z\in \mathbb C\setminus \mathbb R_+$.
In the above domain, 
\begin{equation}
\label{Resolv}
R(z) = \int_{\mathbb R} \frac {dE(\mu)}{\mu-z}
\end{equation}
We shall use the following simplified version of the limiting
absorption principle stating that
$R(z)= R(\la+ib)$ has well-defined operator limits 
$R_+(\la)$ and $R_{-}(\la)$, for $\la>0$, as $b\to 0^+$ and 
$b\to 0^-$ respectively. The operators $R_{\pm}(\la)$ map the 
space of Schwartz functions ${\mathcal S}$ into the space 
$C^\infty\cap L^4({\mathbb R}^3)$.

On the real axis, the resolvent
$R(\la)$ can be then described explicitly as the integral operators
with the kernels
\begin{align}
&R_+(\la)(x,y) = \lim_{\eps\to 0^+}R(\la+i\eps)(x,y)=
\frac {e^{i\sqrt \la|x-y|}}{4\pi|x-y|},\qquad \la\ge 0\label{Rrep}\\
&R_-(\la)(x,y) =  \lim_{\eps\to 0^+}R(\la-i\eps)(x,y)=
\overline {R_+(\la)(x,y)} = \frac {e^{-i\sqrt \la|x-y|}}{4\pi|x-y|},
\qquad \la \ge 0,\nn\\
&R(\la)(x,y) = \overline{R(\la)(x,y)}=
\frac {e^{-\sqrt {-\la}|x-y|}}{4\pi|x-y|},\qquad \la<0.\nn
\end{align}
In particular, we can write
\begin{equation}
\label{spmeas}
dE(\la) = \Im R(\la).
\end{equation}
We shall make repeated use of the following regularization:
$$
\int_a^b e^{i\a q}dq = \frac {e^{i(\a+i0) b}-e^{i(\a+i0) a}}{\a+i0}=
\lim_{\eps\to 0^+} \frac {e^{i(\a+i\eps) b}-e^{i(\a+i\eps)a}}{\a+i\eps}.
$$
which holds true for any finite $a,b\in \mathbb R$ and arbitrary $\a\in \mathbb R$.
\begin{proposition}
The function $\langle{\mathcal I}_m(t,s)\psi_s,\,g\rangle$ defined in~\eqref{funcIm}, 
the $m$-th term of the Born series 
\eqref{eq:series}, admits the following 
representation:
\begin{align}
{\mathcal I}_m(t,s)=i^m\int_{{\mathbb R}^m} d\tau_1..d\tau_m e^{it(\tau_1+..+\tau_m)}
\int_\la e^{i(t-s)\la}\sum_{k=0}^{m+1}\bigg\langle \bigg(\prod_{r=1}^{k-1}
R_+(\la+\tau_r+..+\tau_m)V(\hat\tau_r,\cdot)\bigg)&\nn\\
dE(\la +\tau_k+..+\tau_m)\bigg(\prod_{r=k+1}^{m+1} 
V(\hat\tau_{r-1},\cdot)R_-(\la+\tau_r+..+\tau_{m+1})\bigg) \psi_s,\,g\bigg\rangle&,\label{formIm}
\end{align}
where we formally set $\tau_{m+1}=0$. 
The representation holds true with arbitrary Schwartz functions $\psi_s, g\in {\mathcal S}$.
\label{Fcalc}
\end{proposition}
\begin{proof}
We start by verifying that the expression on the right hand-side of~\eqref{formIm}
defines an absolutely convergent integral.
Recall that the potential $V(\hat\tau,x)$ is smooth and has compact support in
both variables. Therefore, the variables $\tau_1,..,\tau_m$ are restricted
to a finite interval of $\mathbb R$. 
It also follows, with the help of our version of the limiting absorption principle,
that the operators $V(\hat\tau,\cdot) R_{\pm}(\la)$ map 
${\mathcal S}$ into  ${\mathcal S}$ for all $\hat\tau, \la \in \mathbb R$. 
In addition, we have that
$$
dE(\la) f = \la^{-N} dE(\la) (\Delta)^N f  
$$
for an arbitrary Schwartz function $f$.
Hence,
\be 
&& \bigg|\bigg\langle \bigg(\prod_{r=1}^{k-1}
R_+(\la+\tau_r+..+\tau_m)V(\hat\tau_r,\cdot)\bigg) \label{eq:firstblock} \\
&& \qquad\qquad dE(\la +\tau_k+..+\tau_m) \bigg(\prod_{r=k+1}^{m+1} 
V(\hat\tau_{r-1},\cdot)R_-(\la+\tau_r+..+\tau_m)\bigg) \psi_s, g\bigg\rangle \bigg| \nn \\ 
&& \le C (1+|\la|)^{-N} \label{eq:opmeas}
\ee
for arbitrary Schwartz functions $\psi_s$ and $g$  with a constant $C$ depending 
on $\psi_s$, $g$, and $V$ 
(in particular, on the size of the support of $V(\hat\tau,x)$ in $\hat\tau$). This can be seen most easily
by moving the operator in~\eqref{eq:firstblock} onto~$g$. 

\noindent In what follows we shall manipulate the operator valued expressions 
with the tacit understanding that all equalities are to be interpreted in the weak sense.
However, for ease of notation we suppress the pairing with the Schwartz functions $\psi_s$ and $g$.
The absolute convergence of all of integrals involved (after silent pairing with
$\psi_s, g$) will also allow us to freely interchange the order of integrations.
  
\noindent We replace each of the $e^{i(s_k-s_{k-1})H_0}$ in 
\eqref{funcIm} with its 
spectral representation:
\begin{align}
{\mathcal I}_m(t,s) = \int_{{\mathbb R}^m} d\tau_1..d\tau_m 
\idotsint\limits_{\la_1,..,\la_{m+1}} \,\,
\idotsint\limits_{s\le s_m\le..\le s_1\le t}ds_1...ds_m\;  
e^{i(t-s_1)\la_1} e^{is_1\tau_1}dE(\la_1)\,V(\hat\tau_1,\cdot) e^{i(s_1-s_2)\la_2}
&\nn \\e^{is_2\tau_2}dE(\la_2) V(\hat\tau_2,\cdot) ...
e^{i(s_{m-1}-s_m)\la_m} e^{is_m\tau_m}dE(\la_m) 
V(\hat\tau_m,\cdot)e^{i(s_m-s)\la_{m+1}}dE(\la_{m+1})&\nn\\=
\int_{{\mathbb R}^m} d\tau_1..d\tau_m 
\idotsint\limits_{\la_1,..,\la_{m+1}} \,\,
\idotsint\limits_{s\le s_m\le..\le s_1\le t}ds_1...ds_m\;
e^{it\la_1} dE(\la_1)\,V(\hat\tau_1,\cdot)e^{is_1(\tau_1-\la_1+\la_2)}
dE(\la_2) V(\hat\tau_2,\cdot)&\nn \\e^{is_2(\tau_2 -\la_2 +\la_3)} ...
\, dE(\la_m) V(\hat\tau_m,\cdot)e^{is_m(\tau_m -\la_m +\la_{m+1})}
dE(\la_{m+1}) \,e^{-is\la_{m+1}}.&\label{Im}
\end{align}
Consider the first term
$$
{\mathcal I}_1 = \int_{\mathbb R}d\tau_1 \int_{\la_1,\la_2} \int_s^t ds_1\;
e^{it\la_1} dE(\la_1)\,V(\hat\tau_1,\cdot)e^{is_1(\tau_1-\la_1+\la_2)} 
dE(\la_{2}) e^{-is\la_2}.\,
$$
Integrating explicitly relative to $s_1$ we infer that
$$
\aligned
{\mathcal I}_1(t,s)& = -i \int_{\mathbb R}d\tau_1 \int_{\la_1,\la_2} 
e^{it\la_1} dE(\la_1)\,V(\hat\tau_1,\cdot)\frac {e^{it(\tau_1-\la_1+\la_2+i0)}-
e^{is(\tau_1-\la_1+\la_2+i0)}}{\tau_1-\la_1+\la_2+i0}dE(\la_{2}) e^{-is\la_2}\\
&=-i \int_{\mathbb R}d\tau_1 e^{it(\tau_1+i0)} \int_{\la_1,\la_2} 
dE(\la_1)\,V(\hat\tau_1,\cdot)\frac {e^{i(t-s)\la_2}}{\tau_1-\la_1+\la_2+i0}dE(\la_{2})
\\ &+ i \int_{\mathbb R}d\tau_1 e^{is(\tau_1+i0)} \int_{\la_1,\la_2} 
dE(\la_1)\,V(\hat\tau_1,\cdot)\frac {e^{i(t-s)\la_1}}{\tau_1-\la_1+\la_2+i0}dE(\la_{2})
\\ &= i \int_{\mathbb R}d\tau_1 e^{it\tau_1} \int_{\la_2}{e^{i(t-s)\la_2}} 
R_{+}(\la_2+\tau_1)\,V(\hat\tau_1,\cdot) dE(\la_{2})\\
&+ i \int_{\mathbb R}d\tau_1 e^{is\tau_1} \int_{\la_1}{e^{i(t-s)\la_1}} 
dE(\la_1)\,V(\hat\tau_1,\cdot) R_{-}(\la_{1}-\tau_1)
\\ &=  i \int_{\mathbb R}d\tau_1 e^{it\tau_1} \int_{\la}{e^{i(t-s)\la}} 
\bigg(R_{+}(\la+\tau_1)\,V(\hat\tau_1,\cdot) dE(\la) + 
dE(\la+\tau_1)\,V(\hat\tau_1,\cdot) R_{-}(\la)\bigg).
\endaligned
$$
In the above calculation we have used the spectral representation
\eqref{Resolv} for the resolvent and~\eqref{Rrep}.
The proof now proceeds inductively. We shall assume that 
\begin{equation}\label{IdM}
{\mathcal I}_m(t,s) =i^m\int_{{\mathbb R}^m} d\tau_1\ldots d\tau_m\; e^{it(\tau_1+..+\tau_m)}
\int_\la e^{i(t-s)\la} dM_m(\la;\tau_1,..,\tau_m),
\end{equation}
where $dM_m(\la;\tau_1,..,\tau_m)$ is the operator valued measure
\footnote{Once again we make sense of $dM_m(\la;\tau_1,..,\tau_m)$ only after
pairing it with the Schwartz functions $\psi_s$ and $g$. 
Then $\langle dM_m(\la;\tau_1,..,\tau_m) \psi_s, g\rangle$ is a finite
measure relative to $\la$ --- in fact, rapidly decaying in $\la$, see~\eqref{eq:opmeas} ---  which depends
smoothly on $\tau_1,..,\tau_m$ and vanishes outside of a compact set in these variables.}
defined by 
$$
\aligned
dM_m(\la;\tau_1,..,\tau_m)=\sum_{k=0}^m \bigg[R_+(\la+\tau_1+..+\tau_m)V(\hat\tau_1,\cdot)
R_+(\la+\tau_2+..+\tau_m)V(\hat\tau_2,\cdot)...V(\hat\tau_{k-1},\cdot)&\\
dE(\la +\tau_k+..+\tau_m)
V(\hat\tau_k,\cdot)R_-(\la+\tau_{k+1}+..+\tau_m)V(\hat\tau_{k+1},\cdot)...
V(\hat\tau_m,\cdot)R_-(\la)\bigg].&
\endaligned
$$
Formally setting $\tau_{m+1}=0$, we can also write the above expression in
the following more concise form:
\begin{align}
dM_m(\la;\tau_1,..,\tau_m)=\sum_{k=0}^{m+1} &\bigg(\prod_{r=1}^{k-1}
R_+(\la+\tau_r+..+\tau_m)V(\hat\tau_r,\cdot)\bigg)
dE(\la +\tau_k+..+\tau_m) \nn\\
&\bigg(\prod_{r=k+1}^{m+1} V(\hat\tau_{r-1},\cdot)R_-(\la+\tau_r+..+\tau_{m+1})\bigg).
\label{dMm}
\end{align}
We have already verified~\eqref{IdM} for $m=1$. It remains to check that
$$
{\mathcal I}_{m+1}(t,s) =i^{m+1}\int_{{\mathbb R}^{m+1}} d\tau_1..d\tau_{m+1} \,
e^{it(\tau_1+..+\tau_m +\tau_{m+1})}\int_\la e^{i(t-s)\la} 
dM_{m+1}(\la;\tau_1,..,\tau_{m+1}).
$$
We can deduce from~\eqref{Im} the following recursive identity:
$$
{\mathcal I}_{m+1}(t,s) = \int_{\tau_1}d\tau_1\int_{\la_1} \int_{s}^t ds_1 e^{it\la_1}
dE(\la_1) V(\hat \tau_1,\cdot)\, e^{is_1(\tau_1 - \la_1)} {\mathcal I}_m(s_1,s).
$$
Substituting the expression for ${\mathcal I}_m$ from~\eqref{IdM} we obtain
$$
\aligned
{\mathcal I}_{m+1}(t,s) = i^m \int_{{\mathbb R}^{m+1}}d\tau_1... d\tau_{m+1}  
\int_{\la,\la_1} \int_{s}^t ds_1 e^{it\la_1}
dE(\la_1) V(\hat \tau_1,\cdot)\,&\\ e^{is_1(\tau_1 +.. + \tau_{m+1}- \la_1)}\,
e^{i(s_1-s)\la} dM_m(\la;\tau_2,..,\tau_{m+1}).&
\endaligned
$$
Integrating explicitly relative to $s_1$ we infer that
$$
\aligned
{\mathcal I}_{m+1}(t,s) &= -i^{m+1} \int_{{\mathbb R}^{m+1}}d\tau_1... d\tau_{m+1}  
\int_{\la,\la_1}e^{it\la_1}e^{-is\la}
dE(\la_1) V(\hat \tau_1,\cdot)\,\\ \qquad&\frac {e^{it(\tau_1 +..+\tau_{m+1}- \la_1+ \la +i0)}-
e^{is(\tau_1 +..+\tau_{m+1}- \la_1+ \la+i0)}}
{\tau_1 +..+\tau_{m+1} - \la_1+ \la+i0}  dM_m(\la;\tau_2,..,\tau_{m+1})
\\ &= -i^{m+1} \int_{{\mathbb R}^{m+1}}d\tau_1... d\tau_{m+1}\; 
e^{it(\tau_1+..+\tau_{m+1}+i0)}
\int_{\la,\la_1}e^{i(t-s)\la}\\
&\qquad\qquad\frac {dE(\la_1)}{\tau_1 +..+\tau_{m+1}- \la_1+ \la+i0}
V(\hat\tau_1,\cdot) dM_m(\la;\tau_2,..,\tau_{m+1})
\\&+ \,i^{m+1}  \int_{{\mathbb R}^{m+1}}d\tau_1... d\tau_{m+1}\;  
e^{is(\tau_1+..+\tau_{m+1}+i0)}\,
\int_{\la,\la_1}e^{i(t-s)\la_1}
dE(\la_1)\\ &\qquad\qquad V(\hat\tau_1,\cdot)
\frac{dM_m(\la;\tau_2,..,\tau_{m+1})}{\tau_1 +..+\tau_{m+1}- \la_1+ \la+i0}= J_1 + J_2. 
\endaligned
$$
According to~\eqref{Resolv} and~\eqref{Rrep}
$$
\int_{\la_1}\frac {dE(\la_1)}{\tau_1+..+\tau_{m+1} - \la_1+ \la+i0}= - 
R_+(\la+\tau_1+..+\tau_{m+1}).
$$
Therefore,
\begin{equation}
J_1= i^{m+1} \int_{{\mathbb R}^{m+1}}d\tau_1... d\tau_{m+1}  
e^{it(\tau_1+..+\tau_{m+1})}
\int_{\la}e^{i(t-s)\la} R_+(\la+\tau_1+..+\tau_{m+1}) V(\hat\tau_1,\cdot) 
dM_m(\la;\tau_2,..,\tau_{m+1}).\label{J1}
\end{equation}
Observe that, with the convention that $\tau_{m+2}=0$,
\begin{align}
R_+(\la+\tau_1+..+\tau_{m+1})  V(\hat\tau_1,\cdot) 
dM_m(\la;\tau_2,..,\tau_{m+1})  = 
\sum_{k=2}^{m+2} \bigg[\big(\prod_{r=1}^{k-1}
R_+(\la+\tau_r+..+\tau_{m+1})V(\hat\tau_r,\cdot)\big)&\nn\\
dE(\la +\tau_k+..+\tau_{m+1}) \big(\prod_{r=k+1}^{m+2} 
V(\hat\tau_{r-1},\cdot)R_-(\la+\tau_r+..+\tau_{m+2})\big)\bigg].&\label{RMm}
\end{align}
It remains to consider the integral $J_2$.
\begin{align}
J_2& = i^{m+1}  \int\limits_{{\mathbb R}^{m+1}}d\tau_1... d\tau_{m+1}  
e^{is(\tau_1+..+\tau_{m+1}+i0)}\, 
\int_{\la,\la_1}e^{i(t-s)\la_1}
dE(\la_1)V(\hat\tau_1,\cdot)\frac{dM_m(\la;\tau_2,..,\tau_{m+1})}
{\tau_1+..+\tau_{m+1} - \la_1+ \la+i0}\nn\\
&=  i^{m+1}  \int_{{\mathbb R}^{m+1}}d\tau_1... d\tau_{m+1}  
e^{it(\tau_1+..+\tau_{m+1})+is (i0) }\, 
\int_{\la_1}e^{i(t-s)\la_1}
dE(\la_1+\tau_1+..+\tau_{m+1})&\label{J2}\\
&\hskip 12pc V(\hat\tau_1,\cdot)\int_\la 
\frac{dM_m(\la;\tau_2,..,\tau_{m+1})}{\la - \la_1 +i0}\nn
\end{align}
Inspection of the desired expression for $dM_{m+1}(\la;\tau_1,..,\tau_{m+1})$ and
equations~\eqref{J1}-\eqref{J2} suffices to verify the following formula: 
$$
\int_\la \frac{dM_m(\la;\tau_2,..,\tau_{m+1})}{\la - \la_1+i0}=
\bigg(\prod_{r=2}^{m+1} R_-(\la_1+\tau_r+..+\tau_{m+1})V(\hat\tau_r,\cdot)\bigg)
R_-(\la_1)\,.
$$
This is accomplished in the following two lemmas, and we are done.
\end{proof}

We recall definition~\eqref{dMm} of the operator valued measure $dM_m$ and
prove the following more general result
\begin{lemma}
Let $a_1,..,a_m\in \mathbb R$ be a sequence of arbitrary real numbers and 
let $A_1,..,A_m$ be arbitrary operators\footnote{It suffices to assume that
that the operators $A_k$, $k=1,..,m$ map the space $C^\infty({\mathbb R}^3)\cap L^4({\mathbb R}^3)$
into the the space ${\mathcal S}$.}. Then
\begin{align}
\int_\la \frac 1{\la-\mu+i0}\sum_{k=1}^m \bigg(\prod_{r=1}^{k-1}R_+(\la+a_r)A_r\bigg)
dE(\la+a_k)\bigg(\prod_{r=k+1}^{m}A_{r-1}R_-(\la+a_r)\bigg)
&\nn\\ =\bigg(\prod_{r=1}^{m-1} R_-(\mu+a_r)A_r\bigg) 
R_-(\mu+a_m)&\label{fract}.
\end{align}
As before, the identity holds after pairing the above expressions with a pair
of Schwartz functions $\psi_s, g$.
\label{Funct}
\end{lemma}
\begin{proof}
We shall write each $R_\pm(\la+a_r)$, for all values of $r=1,..,m$ different
from $k$
using the spectral representation
$$
R_\pm(\la+a_r)= \int_{\la_r} \frac {dE(\la_r+a_r )}{\la_r-\la \mp i0}.
$$
We shall also rename the variable of integration $\la$ to 
$\la_k$ in each term of the sum in $k$.
The left hand-side of~\eqref{fract} then takes the following form:
$$
\idotsint\limits_{\la_1,..,\la_m} \sum_{k=1}^m
 \frac 1{\la_k-\mu+i0}\prod_{r=1}^{k-1} \frac 1{\la_r-\la_k-i0}
\prod_{r=k+1}^{m} \frac 1{\la_r-\la_k+i0}
\bigg(\prod_{j=1}^{m-1} dE(\la_j+a_j) A_j\bigg) dE(\la_m+a_m)
$$
The proof of Lemma~\ref{Funct} is finished provided that we can show 
that the 
following identity holds true:
$$
\sum_{k=1}^m
 \frac 1{\la_k-\mu+i0}\prod_{r=1}^{k-1} \frac 1{\la_r-\la_k-i0}
\prod_{r=k+1}^{m} \frac 1{\la_r-\la_k+i0}=
\prod_{r=1}^m  \frac 1{\la_r-\mu+i0}
$$
In the distributional sense 
$$
\aligned
\lim_{\eps_k\to 0^+}\lim_{\eps_1\to 0^{+}}...\lim_{\eps_{k-1}\to 0^+}
\lim_{\eps_{k+1}\to 0^-}...\lim_{\eps_m\to 0^-}
\frac 1{\la_k-\mu+i\eps_k}\prod_{r=1, r\ne k}^{m} 
\frac 1{\la_r-\la_k-i\eps_r}&\\=
\lim_{\eps\to 0^+} 
\frac 1{\la_k-\mu+i k\eps}\prod_{r=1, r\ne k}^{m} 
\frac 1{\la_r-\la_k+i(r-k)\eps}&
\endaligned
$$
Therefore, we can introduce the new variables $z_r=\la_r-\mu +i r\eps$, $r=1,..,m$
and prove instead the following statement.
\end{proof}
\begin{lemma}
For any pairwise distinct complex numbers $z_1,..,z_m\in \mathbb C$,
$$
\sum_{k=1}^m
 \frac 1{z_k}\prod_{r=1,r\ne k}^m \frac 1{z_r-z_k} =
\prod_{r=1}^m  \frac 1{z_r}.
$$
\end{lemma}
\begin{proof}
The key identity is the statement of the lemma for $m=2$
$$
\frac 1{z_1(z_2-z_1)} + \frac 1{z_2(z_1-z_2)}=\frac 1{z_1z_2}
$$
which follows immediately by inspection.
The general case then can be proved by induction. 
We shall assume that the identity holds true for $m-1$ 
and prove the result for $m$. 
We first note a simple equality
$$
\frac 1{(z_m-z_k)}=\frac 1{(z_m-z_1)} + \frac {z_k-z_1}{(z_m-z_k)(z_m-z_1)}.
$$
Therefore,
$$\aligned
\sum_{k=1}^m
 \frac 1{z_k}\prod_{r=1,r\ne k}^m \frac 1{z_r-z_k}& =
\frac {1}{z_m-z_1}\sum_{k=1}^{m-1}
 \frac 1{z_k}\prod_{r=1,r\ne k}^{m-1} \frac 1{z_r-z_k}
- \frac {1}{z_m-z_1} \sum_{k=2}^{m-1}\frac 1{z_k}
\prod_{r=2,r\ne k}^{m} \frac 1{z_r-z_k} \\&+ 
\frac 1{z_m} \prod_{r=1}^{m-1} \frac 1{z_r-z_m}. 
\endaligned
$$
According to the assumption $m-1$ with $z_1,..,z_{m-1}$ the first
term on the right hand-side gives 
$\frac {1}{(z_m-z_1)z_1\cdots z_{m-1}}$.
We also have
$$\aligned
\frac {1}{z_m-z_1} \sum_{k=2}^{m-1}\frac 1{z_k}
\prod_{r=2,r\ne k}^{m} \frac 1{z_r-z_k}& =
\frac {1}{z_m-z_1} \sum_{k=2}^{m}\frac 1{z_k}
\prod_{r=2,r\ne k}^{m} \frac 1{z_r-z_k} -
\frac {1}{z_m-z_1} \frac 1{z_m}
\prod_{r=2}^{m-1} \frac 1{z_r-z_m} \\ &=
\frac {1}{(z_m-z_1)z_2\cdots z_{m}} -
\frac 1{z_m}
\prod_{r=1}^{m-1} \frac 1{z_r-z_m}
\endaligned
$$
by the $m-1$ inductive assumption for $z_2,..,z_m$.
Finally,
$$
\frac {1}{(z_m-z_1)z_1\cdots z_{m-1}}-\frac {1}{(z_m-z_1)z_2\cdots z_{m}}
=\frac 1{z_1\cdots z_m},
$$
as desired.
\end{proof}
We shall now derive the explicit representation of the integral kernel
of the operator ${\mathcal I}_m(t,s)$ acting on the Schwartz functions
$\psi_s$. We start by noting the
following simple identity which holds for arbitrary real numbers $a_1,..,a_{m+1}$
with $m\ge 1$:
\begin{equation}\label{sinus}
\sum_{k=1}^{m+1} e^{i(a_1+..+a_{k-1}-a_{k+1}-..-a_{m+1})}\sin a_k =
\sin ({\sum_{k=1}^{m+1} a_k}).
\end{equation}
This identity can be easily proved by induction on $m$.
Recall that 
$$
R_{\pm}(\mu)(x,y)=\frac {e^{\pm i\sqrt{\mu}|x-y|}}{4\pi|x-y|},\quad \mu\in \mathbb R
$$
with $\sqrt\mu$ defined in such a way that 
$\mbox{Im}\sqrt{\mu}>0$ for $\mbox{Im} \mu >0$. We have
$R_+(\mu)=R_-(\mu)$ for $\mu<0$.
Also recall that the kernel of the spectral measure
$$
dE(\mu)(x,y) = \begin{cases} 
\frac{\sin \sqrt{\mu}|x-y|}{4\pi|x-y|}\,d\mu & \mu>0,\\
0& \mu\le 0
\end{cases}.
$$ 
We return to the representation~\eqref{formIm} for the ${\mathcal I}_m$.
Let (with $\tau_{m+1}=0$)
\[\tau_j+..+\tau_{m+1}=\min_{r\in [1,m]} (\tau_r+..+\tau_{m+1}).\]
To simplify the formulae we introduce a new operator 
${\mathcal J}_m(t,s)$, implicitly dependent on $\tau_1,..,\tau_m$,
\begin{align}
&{\mathcal I}_m(t,s)= i^m \int_{{\mathbb R}^m} d\tau_1..d\tau_m\; e^{it(\tau_1+..+\tau_m)}
e^{-i(t-s)(\tau_j+..+\tau_{m+1})}{\mathcal J}_m(t,s)(\tau_1,..,\tau_m),\nn\\
&{\mathcal J}_m(t,s):= \int_\la e^{i(t-s)(\la + \tau_j+..+\tau_{m+1})} 
\sum_{k=1}^{m+1}\bigg(\prod_{r=1}^{k-1}
R_+(\la+\tau_r+..+\tau_m)V(\hat\tau_r,\cdot)\bigg)
dE(\la +\tau_k+..+\tau_m)\label{Jm}\\&\hskip 10pc\bigg(\prod_{r=k+1}^{m+1} 
V(\hat\tau_{r-1},\cdot)R_-(\la+\tau_r+..+\tau_{m+1})\bigg)\nn.
\end{align}
Define non-negative numbers $\si_r$, $r=1,..,m+1$
$$
\si_r = (\tau_r+..+\tau_{m+1})-(\tau_j+..+\tau_{m+1}).
$$
After a change of variables we obtain the expression
$$
{\mathcal J}_m(t,s) = \int_\la e^{i(t-s)\la} \sum_{k=1}^{m+1}\bigg(\prod_{r=1}^{k-1}
R_+(\la+\si_r)V(\hat\tau_r,\cdot)\bigg)
dE(\la +\si_k)\bigg(\prod_{r=k+1}^{m+1} 
V(\hat\tau_{r-1},\cdot)R_-(\la+\si_r)\bigg).
$$
Observe that due to the presence of $dE(\la+\si_k)$ the $k^{th}$ 
term in the sum above vanishes for $\la\le - \si_k\le 0$. Therefore,
\begin{align}
\ &{\mathcal J}_m(t,s) = {\mathcal L}_m(t,s) + {\mathcal M}_m(t,s)\nn \\
\ &{\mathcal L}_m(t,s):=\int_{0}^\infty  e^{i(t-s)\la}\sum_{k=0}^{m+1}
\bigg(\prod_{r=1}^{k-1}
R_+(\la+\si_r)V(\hat\tau_r,\cdot)\bigg)
dE(\la +\si_k)\bigg(\prod_{r=k+1}^{m+1} 
V(\hat\tau_{r-1},\cdot)R_-(\la+\si_r)\bigg)\label{JLM} \\ 
\ &{\mathcal M}_m(t,s):= \int_{-\infty}^0  
e^{i(t-s)\la} \sum_{k=1}^{m+1} \bigg(\prod_{r=1}^{k-1}
R_+(\la+\si_r)V(\hat\tau_r,\cdot)\bigg)
dE(\la +\si_k)\bigg(\prod_{r=k+1}^{m+1} 
V(\hat\tau_{r-1},\cdot)R_-(\la+\si_r)\bigg).\nn
\end{align}
To obtain the explicit formula for the integral kernel of the operator 
${\mathcal L}_m(t,s)$ we make use of the following: the parameters $\si_k\ge 0$,
$\la\ge 0$ on the interval of integration, and the explicit representations 
for the kernels of $R_\pm(\mu)$ and $dE(\mu)$. We have
$$
\aligned
{\mathcal L}_m(t,s)&(x,y)=\int_{{\mathbb R}^{m}}dx_1..dx_{m}
\int_{0}^\infty  d\la \,e^{i(t-s)\la}\sum_{k=1}^{m+1}\bigg[
e^{i(\sqrt{\la+\si_1}|x-x_1|+..+\sqrt{\la+\si_{k-1}}|x_{k-2}-x_{k-1}|)}\\
&e^{-i(\sqrt{\la+\si_{k+1}}|x_{k}-x_{k+1}|-..-\sqrt{\la+\si_{m}}|x_{m-1}-y|)}
\sin {\big(\sqrt{\la+\si_{k}}|x_{k-1}-x_{k}|\big)} \prod_{r=1}^{m} 
\frac{V(\hat\tau_{r},x_r)}{4\pi|x_{r-1}-x_r|}\;\frac{1}{4\pi|x_m-y|}\bigg],   
\endaligned
$$
where we set $x_0=x$.
We now recall the identity~\eqref{sinus} to infer that 
$$
{\mathcal L}_m(t,s)(x,y)=\int_{{\mathbb R}^{m}}dx_1..dx_{m} \,
\prod_{r=1}^{m} \frac{V(\hat\tau_{r},x_r)}{4\pi|x_{r-1}-x_r|}\frac{1}{4\pi|x_m-y|}
\int_{0}^\infty d\la\,  e^{i(t-s)\la}
\sin{\bigg(\sum_{k=1}^{m+1}\sqrt{\la+\si_k}|x_{k-1}-x_k|\bigg)}.
$$
Changing variables in the $\la$-integral and integrating  by parts yield
$$
\aligned
\int_{0}^\infty  e^{i(t-s)\la}
\sin{\bigg(\sum_{k=1}^{m}\sqrt{\la+\si_k}|x_{k-1}-x_k|\bigg)}=
2 \int_{0}^\infty d\la\, \la \, e^{i(t-s)\la^2}
\sin{\bigg(\sum_{k=1}^{m}\sqrt{\la^2+\si_k}|x_{k-1}-x_k|\bigg)}&\\=
\frac {i}{t-s}\sum_{\ell=1}^m 
\int_{0}^\infty d\la\, e^{i(t-s)\la^2}
\cos{\bigg(\sum_{k=1}^{m}\sqrt{\la^2+\si_k}|x_{k-1}-x_k|\bigg)}
\frac {\la}{\sqrt{\la^2+\si_{\ell}}}|x_{\ell-1}-x_{\ell}|.&
\endaligned
$$
Therefore, finally
\begin{align}
\ &{\mathcal L}_m(t,s)(x,y)=\frac {i}{t-s}\sum_{\ell=1}^{m} {\mathcal L}^{\ell}_m(t,s)(x,y),\nn\\
\ & {\mathcal L}^{\ell}_m(t,s)(x,y):= \int_{{\mathbb R}^{m}}dx_1\ldots dx_{m}
\prod_{r=1}^{m} \frac{V(\hat\tau_{r},x_r)}{4\pi|x_{r-1}-x_r|}\,\,
\frac{|x_{\ell-1}-x_{\ell}|}{4\pi|x_m-y|}\label{finLm}\\
&\hskip 8pc \int_{0}^\infty d\la\, e^{i(t-s)\la^2}
\cos{\bigg(\sum_{k=1}^{m}\sqrt{\la^2+\si_k}|x_{k-1}-x_k|\bigg)}
\frac {\la}{\sqrt{\la^2+\si_{\ell}}}.\nn
\end{align}
To describe the integral kernels of the operators ${\mathcal M}_m^k(t,s)$ 
we shall first order and rename the parameters $\si_k$, $k=1,..,m+1$.
In fact, define inductively
$$
\omega_d=\max \{\si_k\}_{k\in [1,m+1]}\setminus \{\omega_{\ell}\}_{\ell\in [1,d-1]}, 
$$
and set $k=k(c)$ and $c=c(k)$ iff $\si_k=\omega_c$.
We shall split the interval of integration in $\la$ in $(-\infty,0]$ 
into the subintervals 
\[ (-\infty,-\sqrt{\omega_1}\;],
[-\sqrt{\omega_{d-1}},-\sqrt{\omega_{d}}\;]\text{\ \  for\ \ }d\in [2,m+1],\text{\ \  and \ \ }
[-\sqrt{\omega}_{m+1},0]. \]
For $\la\in [-\sqrt{\omega}_{d-1},-\sqrt{\omega}_{d}]$, the 
spectral measures $dE(\la+\si_{k(c)})=dE(\la+\omega_c)$ vanish
for all $c\ge d$. 
Therefore, with the convention that $\omega_0=\infty$ and $\omega_{m+2}=0$,
we have
\begin{align}
&{\mathcal M}_m(t,s)= \sum_{d=1}^{m+2} {\mathcal M}_m^d,\nn\\
&{\mathcal M}_m^d(t,s):= \int\limits_{-\sqrt{\omega_{d-1}}}^{-\sqrt{\omega_{d}}}
e^{i(t-s)\la} \sum_{c=1}^{d} \bigg(\prod_{r=1}^{k(c)-1}
R_+(\la+\omega_{a(r)}) V(\hat\tau_r,\cdot)\bigg)
dE(\la +\omega_c )\label{Md}\\& \hskip 12pc \bigg(\prod_{r=k(c)+1}^{m+1} 
V(\hat\tau_{r-1},\cdot)R_-(\la+\omega_{a(r)})\bigg). \nn
\end{align}
The integral kernels of $R_\pm(\la+\omega_{a(r)})$ for 
$a(r)\le d-1$ contribute oscillating exponential phases
while for $a(r)\ge d$ they produce exponentially decaying
factors.
Hence,
$$
\aligned
{\mathcal M}_m^d(t,s)(x,y)=
\int_{{\mathbb R}^m}& dx_1..dx_m 
\prod_{r=1}^{m} 
\frac{V(\hat\tau_{r},x_r)}{4\pi|x_{r-1}-x_r|}\frac{1}{4\pi|x_m-y|}
\int\limits_{-\sqrt{\omega_{d-1}}}^{-\sqrt{\omega_{d}}}\;d\lambda
e^{i(t-s)\la}\\ &\sum_{c=1}^{d-1} e^{i(\sqrt{\la + \omega_1}|x_{k(1)-1}-x_{k(1)}|+..+
\sqrt{\la + \omega_{c-1}}|x_{k(c)-2}-x_{k(c)-1}|)}\\
\ &e^{-i(\sqrt{\la + \omega_{c+1}}|x_{k(c+1)-1}-x_{k(c+1)}|+..+
\sqrt{\la + \omega_{d-1}}|x_{k(d-1)-1}-x_{k(d-1)}|)}\\
\ &\sin \big(\sqrt{\la + \omega_c}
|x_{k(c)-1}-x_{k(c)}|\big) 
e^{-\sum_{a=d+1}^{m+1}\sqrt{-\omega_a-\la}|x_{k(a)-1}-x_{k(a)}|}.
\endaligned
$$
Once again we recall the identity~\eqref{sinus} to infer that 
$$
\aligned
\sum_{c=1}^{d-1} e^{i\big(\sum_{a=1}^{c-1}\sqrt{\la + \omega_a}|x_{k(a)-1}-x_{k(a)}|-
\sum_{b=c+1}^{d-1}\sqrt{\la + \omega_{b+1}}|x_{k(b+1)-1}-x_{k(b+1)}|\big)}
\,\, \sin \big(\sqrt{\la + \omega_c}|x_{k(c)-1}-x_{k(c)}|\big)
& \\ =\sin\bigg(\sum_{c=1}^{d-1} \sqrt{\la + \omega_c}
|x_{k(c)-1}-x_{k(c)}|\bigg).&
\endaligned
$$
Therefore,
$$
\aligned
{\mathcal M}_m^d(t,s)(x,y)= 
&\int_{{\mathbb R}^m} dx_1..dx_m 
\prod_{r=1}^{m} 
\frac{V(\hat\tau_{r},x_r)}{4\pi|x_{r-1}-x_r|}\frac{1}{4\pi|x_m-y|}
\\ &\int\limits_{-\sqrt{\omega_{d-1}}}^{-\sqrt{\omega_{d}}}d\la\,\,
e^{i(t-s)\la}\,\sin\bigg(\sum_{c=1}^{d-1} \sqrt{\la + \omega_c}|x_{k(c)-1}-x_{k(c)}|\bigg)
e^{-\sum_{a=d}^{m+1}\sqrt{-\omega_a-\la}|x_{k(a)-1}-x_{k(a)}|}.
\endaligned
$$
We would like to change variables $\la\to \la^2$ and integrate by parts 
relative to $\la$, as we did for ${\mathcal L}_m$. 
Denote the $\la$-integrand $F_d(\la)$ in each of the kernels ${\mathcal M}_m^d(t,s)(x,y)$.
It is not difficult to see that $F_d(-\omega_d) = F_{d+1}(-\omega_d)$ for 
$d=1,..,m+1$. Therefore, the boundary terms will cancel each other telescopically. The boundary terms at 
the two endpoints $\omega_{m+2}=0$ and $\omega_0=\infty$ will also disappear
as in ${\mathcal L}_m$.
This allows us, in what follows, to ignore the boundary terms altogether.
We now make a change of variables $\la\to \la^2 - \omega_{d-1}$. 
We also re-introduce the notation $\si_a$ in the new capacity:
$$
\aligned
&0\le \si_a = \omega_a -\omega_{d-1},\qquad a=0,..,d-1,\\
&0\le \rho_a = \omega_{d-1} -\omega_a,\qquad a=d,..,m+2.
\endaligned
$$
Thus
$$
\aligned
{\mathcal M}_m^d(t,s)(x,y)= 
&\int_{{\mathbb R}^m} dx_1..dx_m 
\prod_{r=1}^{m} 
\frac{V(\hat\tau_{r},x_r)}{4\pi|x_{r-1}-x_r|}\frac{1}{4\pi|x_m-y|}
\\ &\int\limits_{0}^{\sqrt{\rho_d}}d\la\,\,\la
e^{i(t-s)\la^2}\sin\bigg(\sum_{c=1}^{d-1} \sqrt{\la^2 + \si_c}|x_{k(c)-1}-x_{k(c)}|\bigg)
e^{-\sum_{a=d}^{m+1}\sqrt{\rho_a-\la^2}|x_{k(a)-1}-x_{k(a)}|}.
\endaligned
$$
Integrating by parts relative to $\la$ and canceling the contribution from 
the boundary terms as explained above, we finally obtain
$$
{\mathcal M}_m^d(t,s)(x,y)= \frac{i}{t-s}\sum_{\ell=1}^{d-1} {\mathcal M}_m^{d,\ell}(t,s)(x,y) +
\frac{i}{t-s} \sum_{\ell=d}^{m+1} \widetilde{\mathcal M}_m^{d,\ell}(t,s)(x,y),
$$
\begin{align}
{\mathcal M}_m^{d,\ell}(t,s)&(x,y):=\int_{{\mathbb R}^m} dx_1..dx_m\; 
\prod_{r=1}^{m} 
\frac{V(\hat\tau_{r},x_r)}{4\pi|x_{r-1}-x_r|}\frac{|x_{k(l)-1}-x_{k(l)}|}{4\pi|x_m-y|}
\label{Mdl}\\ &\int\limits_{0}^{\sqrt{\rho_d}}d\la\,\,
e^{i(t-s)\la^2}\cos\bigg(\sum_{c=1}^{d-1} \sqrt{\la^2 + \si_c}|x_{k(c)-1}-x_{k(c)}|\bigg)
e^{-\sum_{a=d}^{m+1}\sqrt{\rho_a-\la^2}|x_{k(a)-1}-x_{k(a)}|}\,
\frac {\la}{\sqrt{\la^2+\si_{\ell}}}\nn
\end{align}
\begin{align}
\widetilde{\mathcal  M}_m^{d,\ell}(t,s)&(x,y):=
- \int_{{\mathbb R}^m}  dx_1..dx_m 
\prod_{r=1}^{m} 
\frac{V(\hat\tau_{r},x_r)}{4\pi|x_{r-1}-x_r|}\frac{|x_{k(l)-1}-x_{k(l)}|}{4\pi|x_m-y|}
\label{Mbardl}\\ &\int\limits_{0}^{\sqrt{\rho_d}}d\la\,\,
e^{i(t-s)\la^2}\cos\bigg(\sum_{c=1}^{d-1} \sqrt{\la^2 + \si_c}|x_{k(c)-1}-x_{k(c)}|\bigg)
e^{-\sum_{a=d}^{m+1}\sqrt{\rho_a-\la^2}|x_{k(a)-1}-x_{k(a)}|}\,
\frac {\la}{\sqrt{\rho_{\ell}-\la^2}}.\nn
\end{align}
Combining~\eqref{Jm}-\eqref{Mbardl} we can state the following

\begin{proposition}
The integral kernel of ${\mathcal I}_m(t,s)$, 
the $m$-th term of the Born series~\eqref{eq:series},
can be written in the following form:
\be 
{\mathcal I}_m(t,s)(x,y) = \frac {i^{m+1}}{t-s}\int_{{\mathbb R}^m} d\tau_1..d\tau_m 
e^{i(\tau_1+..+\tau_m)}&& \nn \\
\bigg(\sum_{\ell=1}^m {\mathcal L}_m^{\ell}(t,s)(x,y)(\si_1,..,\si_m) +
&& \sum_{d=0}^{m+2} \sum_{\ell=1}^{d-1}{\mathcal M}_m^{d,\ell}(t,s)(x,y)
(\si_1,..,\si_{d-1},\rho_d,..,\rho_{m+1}) + \nn \\
 && \sum_{d=0}^{m+2} \sum_{\ell=1}^{d-1}\widetilde{\mathcal M}_m^{d,\ell}(t,s)(x,y)
(\si_1,..,\si_{d-1},\rho_d,..,\rho_{m+1})\bigg). \label{eq:kernrep}
\ee
We interpret ${\mathcal I}_m(t,s)(x,y)$ as follows: for any pair of
Schwartz functions $\psi_s$ and $g$
$$
\langle {\mathcal I}_m(t,s)\psi_s,\,g\rangle = 
\int_{{\mathbb R}^6} {\mathcal I}_m(t,s)(x,y)\,\psi_s (y)\,g(x)\,dx\,dy.
$$
The functions  
\[{\mathcal L}_m^{\ell}(t,s)(x,y),\quad {\mathcal M}_m^{d,\ell}(t,s)(x,y),\quad 
 \widetilde{\mathcal  M}_m^{d,\ell}(t,s)(x,y)\]  
are defined in~\eqref{finLm}, \eqref{Mdl}, and~\eqref{Mbardl} correspondingly with implicit 
dependence on the parameters $\si_k, \rho_{\ell}$. The latter are 
positive and depend exclusively and in a linear fashion on 
$\tau_1,..,\tau_m$.
\label{KernRep}
\end{proposition}

\section{Estimates for oscillatory integrals I}\label{sec:time2}

The purpose of this section is to prove the following lemma.
Up to a change of variables in $\lambda$ it provides the estimate on
the oscillatory integral in~\eqref{finLm} that we need.
Strictly speaking, the second bound in~\eqref{eq:bound1} suffices for the dispersive estimate, 
but we include the first for the sake of completeness and possible future applications.

\begin{lemma} \label{lem:osc1}
Let $\omega(\lambda)\ge0$ be a twice differentiable, monotone function such that 
\[ |\omega^{(j)}(\lambda)|\le a_0\,\lambda^{-j}\]
for $j=0,1,2$ and all $\lambda>0$. Then 
there exists a constant $C_0$ which only depends on the constant $a_0$
so that for any positive integer $m$ and any $1\le k\le m$, 
\be
\label{eq:bound1}
\left| 
\int_0^\infty e^{\half i\lambda^2} \,e^{\pm i\sum_{j=1}^m b_j\,\sqrt{\lambda^2+\sigma_j}}
\,\frac{\lambda}{\sqrt{\lambda^2+\sigma_k}}\,\omega(\lambda)\,d\lambda
\right| &\le& C_0\,\min\big[(1+\sigma_1)^{\frac14},m^{\frac32}\,b_k^{-1}\,\max_{\ell} b_\ell\big] \\
\left| 
\int_0^\infty e^{\half i\lambda^2} \,e^{\pm i\sum_{j=1}^m b_j\,\sqrt{\lambda^2+\sigma_j}}
\,\omega(\lambda)\,d\lambda
\right| &\le& C_0\,(1+\sigma_1)^{\frac14} \label{eq:bound1'}
\ee
for any choice of $\sigma_1\ge \sigma_2\ge \ldots\ge \sigma_m \ge0$ and $b_j>0$. 
\end{lemma}

\noindent The proof of this lemma will be broken up into several sublemmas. 
We start with an elementary lemma that establishes a basic estimate on the functions
$\frac{\lambda}{\sqrt{\lambda^2+\tau}}$. We shall use this bound repeatedly throughout this section.

\begin{lemma}
\label{lem:trivlem}
For all nonnegative integers $j$ there exist constants $C_j$ so that
\begin{equation}
\left|\frac{d^j}{d\lambda^j} \frac{\lambda}{\sqrt{\lambda^2+\tau}} \right| \le C_j\,\min(\tau^{-\frac{j}{2}},
\lambda^{-j}\frac{\lambda}{\sqrt{\lambda^2+\tau}}) \label{eq:trivest}
\end{equation}
for all $\tau>0$ and $\lambda>0$.
\end{lemma}
\begin{proof} The case $j=0$ in \eqref{eq:trivest} is obvious. Hence it suffices to consider $j\ge1$.
Let $f(\lambda)=\frac{\lambda}{\sqrt{\lambda^2+1}}$. 
For any $\tau>0$ define $\rho_\tau(\lambda):=\frac{\lambda}{\sqrt{\lambda^2+\tau}}$. Then
$\rho_\tau(\lambda)=f(\lambda\tau^{-\half})$. Since $|f^{(j)}(\lambda)|\le C_j\min(1,\lambda^{-j})$, it follows from the chain rule that
\[ |\rho_\tau^{(j)}(\lambda)| \le C_j\,\tau^{-\frac{j}{2}}\min\Bigl(1,\Bigl(\frac{\sqrt{\tau}}{\lambda}\Bigr)^j\Bigr)
= C_j\min(\tau^{-\frac{j}{2}}, \lambda^{-j})\le C_j\,\lambda^{-j}\frac{\lambda}{\sqrt{\lambda^2+\tau}},\]
as desired. The second inequality here follows 
by considering the cases $\lambda\ge\sqrt{\tau}$ and $\lambda\le\sqrt{\tau}$.
\end{proof}

\noindent We now dispense with the easy case of the phase with the $+$ sign in~\eqref{eq:bound1}. 
By our positivity assumptions this phase has no critical point, but some care is needed in terms of 
upper bounds on higher derivatives.

\begin{lemma}
\label{lem:easycase1}
Let $\omega(\lambda)$ be a differentiable function such that $|\omega^{(j)}(\lambda)|\le a_0\,\lambda^{-j}$
for $j=0,1$ and all $\lambda>0$. 
Then there exists a constant $C_0$ which only depends on the constant $a_0$, 
so that for any positive integer $m$ and any $1\le k\le m$, 
\begin{equation} \nonumber
\left| 
\int_0^\infty e^{\half i\lambda^2 + i\sum_{j=1}^m b_j\,\sqrt{\lambda^2+\sigma_j}}
\,\frac{\lambda}{\sqrt{\lambda^2+\sigma_k}}\,\omega(\lambda)\,d\lambda
\right| \le C_0
\end{equation}
for any choice of $\sigma_1\ge \sigma_2\ge \ldots\ge \sigma_m \ge0$ and $b_j>0$. 
\end{lemma}
\begin{proof}
Let $\phi(\lambda)=\half\lambda^2 + \sum_{j=1}^m b_j\,\sqrt{\lambda^2+\sigma_j}$. Then 
\be
\phi'(\lambda) &=& \lambda\Bigl[1+\sum_{j=1}^m \frac{b_j}{\sqrt{\lambda^2+\sigma_j}}\Bigr] \nonumber\\
\phi''(\lambda) &=& \frac{\phi'(\lambda)}{\lambda} - \eta(\lambda)
\ee
where
\[ \eta(\lambda) = \sum_{j=1}^m \frac{b_j\lambda^2}{(\lambda^2+\sigma_j)^{\frac32}} \le \sum_{j=1}^m 
\frac{b_j}{\sqrt{\lambda^2+\sigma_j}}.\]
Hence $\phi'(\lambda)\ge\lambda$, $\frac{\eta(\lambda)}{\phi'(\lambda)^2}\le \lambda^{-2}$ and therefore also
\[ \frac{|\phi''(\lambda)|}{\phi'(\lambda)^2} \le \frac{1}{\lambda\phi'(\lambda)}+\frac{\eta}{\phi'(\lambda)^2}\le  2\lambda^{-2}.\]
Let $\chi$ be a smooth non-decreasing function with $\chi(\lambda)=0$ for $\lambda\le1$ and $\chi(\lambda)=1$ for $\lambda\ge2$. Then
\begin{equation} 
\left|\int_0^\infty e^{i\phi(\lambda)}\; \omega(\lambda)\frac{\lambda}{\sqrt{\lambda^2+\sigma_k}}\,d\lambda \right| 
\le C + \limsup_{L\to\infty} \left|\int_0^\infty e^{i\phi(\lambda)}\;g_L(\lambda)\,d\lambda \right| 
\label{eq:easysplit}
\end{equation}
where we have set
\[
g_L(\lambda):=\chi(\lambda)(1-\chi(\lambda/L))\,\omega(\lambda)\frac{\lambda}{\sqrt{\lambda^2+\sigma_k}}.
\]
By our assumptions and Lemma~\ref{lem:trivlem}, $|g_L^{(j)}(\lambda)|\le C_j\,\lambda^{-j}$ for $j=0,1$ uniformly in $L$. 
Integrating by parts once inside the integral on the right-hand side of \eqref{eq:easysplit} yields an upper bound of the form
\[ 
\int_0^\infty \Bigl|\frac{d}{d\lambda}\Bigl[\frac{1}{\phi'(\lambda)}g_L(\lambda)\Bigr]\Bigr|\,d\lambda \le 
\int_1^\infty \Bigl[\frac{|\phi''(\lambda)|}{\phi'(\lambda)^2}+\frac{1}{\lambda\phi'(\lambda)}\Bigr]\,d\lambda 
\lesssim \int_1^\infty  \lambda^{-2}\,d\lambda \le C, 
\]
as claimed.
\end{proof}

\noindent We now turn to the phase 
$\phi(\lambda)=\half\lambda^2 - \sum_{j=1}^m b_j\,\sqrt{\lambda^2+\sigma_j}$. 
The following example shows that this phase can vanish to the third order.

\begin{example} Let $b_j>0$ be arbitrary positive numbers for $1\le j\le m$. Set
\[ \phi(\lambda) = \half\lambda^2 - \sum_{j=1}^m b_j\sqrt{\lambda^2 + m^2 b_j^2}.\]
Then
\be
\frac{\phi'(\lambda)}{\lambda} &=& 1 -\sum_{j=1}^m \frac{b_j}{\sqrt{\lambda^2+m^2 b_j^2}} \nonumber\\
\phi''(\lambda) &=& \frac{\phi'(\lambda)}{\lambda} + \sum_{j=1}^m \frac{b_j\lambda^2}{(\lambda^2+m^2 b_j^2)^{\frac32}} \nonumber\\
\Bigl(\frac{\phi'(\lambda)}{\lambda}\Bigr)' &=& \sum_{j=1}^m \frac{b_j\lambda}{(\lambda^2+m^2 b_j^2)^{\frac32}}.
\nonumber
\ee
It follows that $\phi'(0)=\phi''(0)=\phi'''(0)=0$, but
\[ \phi^{(4)}(0) = \frac{3}{m^2}\sum_{j=1}^m b_j^{-2}.\]
If $\phi^{(4)}(0)\ge1$, then heuristically speaking
\[ \Bigl|\int_0^\infty e^{i\phi(\lambda)} d\lambda\Bigr| \asymp \Bigl|\int_0^\infty e^{i\phi^{(4)}(0)\lambda^4} d\lambda\Bigr| \lesssim (\phi^{(4)}(0))^{-\frac14}.\]
In case all $b_j's$ are comparable, this bound agrees with the $\sigma_1^{\frac14}$ estimate in~\eqref{eq:bound1}.
On the other hand, if $\tau^2\phi^{(4)}(0)\ge1$, then again heuristically speaking
\[ \Bigl|\int_0^\infty e^{i\phi(\lambda)} \frac{\lambda}{\sqrt{\lambda^2+\tau}} \, d\lambda\Bigr| \asymp \sqrt{\tau} \Bigl|\int_0^\infty e^{i\tau^2\phi^{(4)}(0)\lambda^4}\frac{\lambda}{\sqrt{\lambda^2+1}} \, d\lambda\Bigr| \lesssim \sqrt{\tau}(\tau^2\phi^{(4)}(0))^{-\frac12}= (\tau\phi^{(4)}(0))^{-\frac12}.\]
In case all $b_j's$ are comparable and one chooses $\tau=\sigma_k$ for an arbitrary $k$, 
then the right-hand side is $O(1)$. This  agrees with the second term in~\eqref{eq:bound1} (ignoring powers of $m$).
\end{example}

\begin{lemma}
\label{lem:lem1} The phase 
\[ \phi(\lambda) = \half\lambda^2 - \sum_{j=1}^m b_j\sqrt{\lambda^2+\sigma_j}\]
has the following properties for any choice of parameters $b_j>0$ and $\sigma_j\ge0$:
\begin{enumerate}
\item There exists at most one critical point $\lambda_0>0$ of $\phi$ on $(0,\infty)$,
with $\phi'(\lambda)>0$ for $\lambda>\lambda_0$, and $\phi'(\lambda)<0$ for $\lambda<\lambda_0$.
Moreover, if $\lambda_0$ exists, then the second derivative $\phi''$ satisfies $\phi''(\lambda)>0$ for $\lambda>\lambda_0$ and there is at most one zero of $\phi''$ on the interval $(0,\lambda_0)$, which we denote by $\lambda_1$. If $\lambda_0$ does not exist, then $\phi''(\lambda)\ge0$ for all $\lambda\ge0$.
\item One has
\begin{equation}
\label{eq:claim} 
\Bigl(\frac{\phi''(\lambda)}{{\phi'}^3(\lambda)}\Bigr)' \le 0
\end{equation}
for all $\lambda\not=\lambda_0$.
\item 
If $\lambda_0>0$ exists, then
\begin{equation}
\label{eq:fund}
\phi'(\lambda) = \lambda\int_{\lambda_0}^\lambda \sum_{j=1}^m \frac{b_js}{(s^2+\sigma_j)^{\frac32}}\,ds,
\end{equation}
whereas if $\lambda_0>0$ does not exist, then
\begin{equation}
\label{eq:fundne}
\phi'(\lambda) \ge \lambda \int_0^\lambda \sum_{j=1}^m \frac{b_js}{(s^2+\sigma_j)^{\frac32}}\,ds.
\end{equation}
\end{enumerate}
\end{lemma}
\begin{proof} One has
\be
\frac{\phi'(\lambda)}{\lambda} &=& 1 -\sum_{j=1}^m \frac{b_j}{\sqrt{\lambda^2+\sigma_j}} \label{eq:ein'}\\
\phi''(\lambda) &=& \frac{\phi'(\lambda)}{\lambda} + \sum_{j=1}^m \frac{b_j\lambda^2}{(\lambda^2+\sigma_j)^{\frac32}} = 1 - \sum_{j=1}^m \frac{b_j\sigma_j}{(\lambda^2+\sigma_j)^{\frac32}}. \label{eq:zwei'}
\ee
It is clear that $\lambda_0>0$ is uniquely given by
\begin{equation}
\label{eq:lam0def}
 1=\sum_{j=1}^m \frac{b_j}{\sqrt{\lambda_0^2+\sigma_j}}.
\end{equation}
Moreover, $\phi''$ is non decreasing by the second equality in \eqref{eq:zwei'}, and also
\[ \phi''(\lambda_0) > 0.\]
If all $\sigma_j=0$, then $\phi''(\lambda)=1$, whereas if one $\sigma_j\not=0$, then $\phi''$
is strictly increasing and can have at most one zero~$\lambda_1$ which then necessarily falls 
into the interval $[0,\lambda_0]$. 
To motivate~\eqref{eq:claim}, observe that for large~$\lambda$ one has $\phi''(\lambda)\asymp1$ and $\phi'(\lambda)\asymp\lambda$ so that 
$\phi''(\lambda){\phi'}^{-3}(\lambda)$ is decreasing. 
However, the situation is not so clear for smaller~$\lambda$ and
we will need to use the specific structure of~$\phi'$ and~$\phi''$. 
Firstly, consider the case $\lambda<\lambda_0$, which is precisely the range where $\phi'<0$. Since 
\[ \phi'''(\lambda) = 3\sum_{j=1}^m \frac{\lambda\sigma_j b_j}{(\lambda^2+\sigma_j)^{\frac52}} \ge 0,\]
it follows that
\[  \frac{d}{d\lambda}\Bigl[\frac{\phi''(\lambda)}{\phi'(\lambda)^{3}} \Bigr] = \frac{\phi'''(\lambda)}{\phi'(\lambda)^3}-3\frac{\phi''(\lambda)^2}{\phi'(\lambda)^4} \le 0
\]
in that range. To deal with the range $\lambda\ge\lambda_0$, we use the identity
\begin{equation} \label{eq:phi2} 
\phi''(\lambda) = \frac{\phi'(\lambda)}{\lambda}+\psi(\lambda)
\end{equation}
where 
\[ \psi(\lambda) = \sum_{j=1}^m \frac{\lambda^2 b_j}{(\lambda^2+\sigma_j)^{\frac32}}.\]
One therefore has
\be
\frac{d}{d\lambda}\Bigl[\frac{\phi''(\lambda)}{\phi'(\lambda)^{3}} \Bigr] &=& 
\frac{d}{d\lambda}\Bigl[\frac{1}{\phi'(\lambda)^{2}}\,\frac{\phi''(\lambda)}{\phi'(\lambda)} \Bigr]
=\frac{d}{d\lambda}\Bigl[\frac{1}{\phi'(\lambda)^{2}}\Bigl(\frac{1}{\lambda}+\frac{\psi(\lambda)}{\phi'(\lambda)^2}\Bigr) \Bigr]
\nonumber \\
&=& -2\frac{\phi''(\lambda)^2}{\phi'(\lambda)^4}
-\frac{1}{\lambda^2\phi'(\lambda)^2} + \frac{\psi'(\lambda)}{\phi'(\lambda)^3} - \frac{\phi''(\lambda)}{\phi'(\lambda)^4} \psi(\lambda).
\label{eq:mess1}
\ee
Next we compute $\psi'$:
\be
\psi'(\lambda) &=& 2\frac{\psi(\lambda)}{\lambda} - \sum_{j=1}^m \frac{3\lambda^3 b_j}{(\lambda^2 + \sigma_j)^{\frac52}} \nonumber\\
&=& \frac{2\phi''(\lambda)}{\lambda} - \frac{2\phi'(\lambda)}{\lambda^2} - \eta(\lambda) \label{eq:psi'}
\ee
where $\eta\ge0$ in \eqref{eq:psi'} is defined to be the sum in the preceding line, and we have used \eqref{eq:phi2} to obtain the first
expression in~\eqref{eq:psi'}. Combining~\eqref{eq:mess1} with~\eqref{eq:psi'} leads to
\begin{equation}  
\frac{d}{d\lambda}\Bigl[\frac{\phi''(\lambda)}{\phi'(\lambda)^{3}} \Bigr] 
= -2\frac{\phi''(\lambda)^2}{\phi'(\lambda)^4}
-\frac{3}{\lambda^2\phi'(\lambda)^2} + \frac{2\phi''(\lambda)}{\lambda\phi'(\lambda)^3} - \frac{\eta(\lambda)}{\phi'(\lambda)^3} - \frac{\phi''(\lambda)}{\phi'(\lambda)^4} \psi(\lambda).
\label{eq:mess2}
\end{equation}
Since both $\phi'(\lambda)\ge0$ and $\phi''(\lambda)\ge0$ in the range $\lambda\ge\lambda_0$, 
all but the third term are negative. We now show that the first and third terms in~\eqref{eq:mess2} result in a negative expression:
\[ -2\frac{\phi''(\lambda)^2}{\phi'(\lambda)^4} + \frac{2\phi''(\lambda)}{\lambda\phi'(\lambda)^3} = \frac{2\phi''(\lambda)}{\phi'(\lambda)^4}\Bigl(\frac{\phi'(\lambda)}{\lambda}-\phi''(\lambda)\Bigr) = - \frac{2\phi''(\lambda)\psi(\lambda)}{\phi'(\lambda)^4},
\]
and \eqref{eq:claim} has been established. 

\noindent As for the final assertion, if $\lambda_0$ exists, then \eqref{eq:fund} follows by integrating
\[ \Bigl(\frac{\phi'(s)}{s}\Bigr)' = \sum_{j=1}^m \frac{b_js}{(s^2+\sigma_j)^{\frac32}}.\]
If $\lambda_0>0$ does not exist, then necessarily 
\begin{equation}
\label{eq:lam0ne}
 1\ge \sum_{j=1}^m \frac{b_j}{\sqrt{\sigma_j}},
\end{equation}
so that $\phi'(s)s^{-1}\ge0$ for all $s\ge0$. Integrating $s$ from $0$ to  $\lambda$ therefore proves
\eqref{eq:fundne}.
\end{proof}

\noindent In the proof of Lemma~\ref{lem:osc1} we consider the intervals $[\lambda_0,\infty)$ and
$[0,\lambda_0)$ separately. We start with the former case, and  present a lemma that reduces 
matters to establishing suitable lower bounds on~$\phi'$. The exact form of the assumed lower bounds
might be unmotivated at this point, but it will become clearer later, see Lemma~\ref{lem:genlower}.

\begin{lemma}
\label{lem:right}
Let 
\[ \phi(\lambda) = \half\lambda^2 - \sum_{j=1}^m b_j\sqrt{\lambda^2+\sigma_j}\]
with arbitrary $b_j>0$ and $\sigma_j\ge0$, 
and suppose $\omega(\lambda)\ge0$ is a twice differentiable, monotone function with
\[ |\omega^{(j)}(\lambda)| \le a_j\,\lambda^{-j}\]
for $j=0,1,2$. 
\begin{enumerate}
\item Suppose a critical point $\lambda_0>0$ of $\phi$ exists, and assume that
\begin{equation}
\label{eq:A} 
\phi'(\lambda)\ge A(\lambda-\lambda_0) \text{\ \ for all\ \ }\lambda>\lambda_0+A^{-\half}
\end{equation}
for some $A>0$.  Then
\begin{equation}
\label{eq:rechts1}
\left|\int_{\lambda_0}^\infty e^{i\phi(\lambda)} \omega(\lambda)\,d\lambda \right| \le C_1\,A^{-\half}.
\end{equation}
If $\lambda_0>0$ does not exist, then \eqref{eq:A} with $\lambda_0=0$ implies \eqref{eq:rechts1}
with $\lambda_0=0$.
\item Let $\tau>0$ and suppose the critical point $\lambda_0>0$ exists and satisfies $\lambda_0<\sqrt{\tau}$.
Assume the lower bounds
\begin{equation}
\label{eq:R}
\phi'(\lambda) \ge \left\{ 
\begin{array}{lcl} c\,R^{-1}\lambda^2(\lambda-\lambda_0) & \text{for} & \lambda_0\le\lambda\le\sqrt{R} \\
                   c (\lambda-\lambda_0) & \text{for} & \lambda \ge \sqrt{R}
\end{array}
\right.
\end{equation}
with constants $1>c>0$ and $R\ge\lambda_0^2$.
Then
\begin{equation}
\label{eq:rechts2}
\left|\int_{\lambda_0}^\infty e^{i\phi(\lambda)} \frac{\lambda}{\sqrt{\lambda^2+\tau}}\,\omega(\lambda)\,d\lambda \right| \le C_1\,\Bigl[c^{-\frac32}+\sqrt{\frac{R}{c\tau}}\Bigr].
\end{equation}
If $\lambda_0>0$ does not exist, then the assumption \eqref{eq:R} with $\lambda_0=0$ implies~\eqref{eq:rechts2}
with $\lambda_0=0$.
\end{enumerate}
In all cases $C_1$ only depends on the constants $a_0,a_1,a_2$. 
\end{lemma}
\begin{proof} Setting $\lambda_0=0$ if a positive critical point does not exist allows us to ignore the issue whether or not such a point exists. Let $\chi$ be a smooth non-decreasing function with $\chi(\lambda)=0$ for $\lambda\le1$ and $\chi(\lambda)=1$ for $\lambda\ge2$.  With an arbitrary parameter $L>0$ define
\[ g_L(\lambda) := \chi(\sqrt{A}(\lambda-\lambda_0))(1-\chi(\lambda/L)) \omega(\lambda).\]
Then for any $j=0,1,2$, 
\be
|g_L^{(j)}(\lambda)| &\le& C\, A^{\frac{j}{2}} \text{\ \ for all\ \ }\lambda
\label{eq:gderiv}
\ee
Clearly,  $\chi(\sqrt{A}(\lambda-\lambda_0))$ provides a cut-off to an interval $[\lambda_0,\lambda_0+2A^{-\half}]$. Thus
\begin{equation}
\label{eq:rechtsteil}
\left|\int_{\lambda_0}^\infty e^{i\phi(\lambda)} (1-\chi(\lambda/L))\,\omega(\lambda)\,d\lambda\right|
\le C\,A^{-\half} +  \left|\int_{\lambda_0}^\infty e^{i\phi(\lambda)} g_L(\lambda)\,d\lambda\right|.
\end{equation}
Since
\[ \Bigl(\frac{d}{d\lambda}\frac{1}{\phi'(\lambda)} \Bigr)^2= 
\phi'(\lambda)^{-2}\frac{d^2}{d\lambda^2} - 3\frac{\phi''(\lambda)}{\phi'(\lambda)^3}\frac{d}{d\lambda} 
- \frac{d}{d\lambda}\Bigl(\frac{\phi''(\lambda)}{\phi'(\lambda)^3}\Bigr),
\]
integrating by parts twice shows that the integral on the right-hand side of \eqref{eq:rechtsteil} is no larger than
\begin{equation} \label{eq:Ione}
 \int_0^\infty |\phi'(\lambda)|^{-2}\, |g_L''(\lambda)|\,d\lambda 
+  3\int_0^\infty \Bigl|\frac{\phi''(\lambda)}{\phi'(\lambda)^{3}}\Bigr|\, |g_L'(\lambda)|\,d\lambda + 
\int_0^\infty \Bigl|\frac{d}{d\lambda}\Bigl[\frac{\phi''(\lambda)}{\phi'(\lambda)^{3}}\Bigr]\Bigr|\, |g_L(\lambda)|\,d\lambda .
\end{equation}
Setting $\gamma=A^{-\half}$ for convenience, \eqref{eq:gderiv} and \eqref{eq:A} imply that the first term on the right-hand side of \eqref{eq:Ione} does not exceed
\[ \int_{\lambda_0+\gamma}^\infty |A(\lambda-\lambda_0)|^{-2}\, \gamma^{-2}\,d\lambda \le C\,A^{-2}\gamma^{-3}
= C\, A^{-\half}\]
uniformly in~$L$. 
In view of \eqref{eq:claim} in Lemma~\ref{lem:lem1} the absolute values in the third term of~\eqref{eq:Ione}
can be taken outside the integral. Integrating by parts therefore reduces the third term to the second. As for the latter, Lemma~\ref{lem:lem1} shows that $\phi''(\lambda){\phi'}^{-3}(\lambda)>0$ for all $\lambda>\lambda_0$.
Whereas $g_L'$ does not have a definite sign, it is the sum of three terms each of which does have a definite sign, say $g_L'=f_1+f_2+f_3$ where $f_k(\lambda)\ge0$ for all $\lambda>\lambda_0$ or $f_k(\lambda)\le0$ for all $\lambda>\lambda_0$. We can therefore move the absolute values outside the second integral provided we
consider each of the $f_k$ separately. Since 
\be 
\int_0^\infty \frac{\phi''(\lambda)}{\phi'(\lambda)^{3}}\, f_j(\lambda)\,d\lambda &=& -\int_0^\infty \phi'(\lambda)\, \Bigl(\frac{f_j(\lambda)}{\phi'(\lambda)^3}\Bigr)'\,d\lambda \nonumber\\
&=& -\int_0^\infty \frac{1}{\phi'(\lambda)^2}\, f_j'(\lambda)\,d\lambda + 3\int_0^\infty \frac{\phi''(\lambda)}{\phi'(\lambda)^{3}}\, f_j(\lambda)\,d\lambda, \nonumber
\ee
and $|f_k'(\lambda)|\le C\gamma^{-2}$ for each $k$ and all $\lambda$ (see~\eqref{eq:gderiv} with $j=2$), the second term satisfies the same estimates as the first.
Passing to the limit $L\to\infty$ in~\eqref{eq:rechtsteil} proves~\eqref{eq:rechts1}.

As for the second part of the lemma, we first consider the case $R\ge1$. With $L>0$ arbitrary and some $\gamma>0$ that will be specified below, set
\[ g_L(\lambda) := \chi\Bigl(\frac{\lambda-\lambda_0}{\gamma}\Bigr)(1-\chi(\lambda/L))\omega(\lambda)\frac{\lambda}{\sqrt{\lambda^2+\tau}}.\]
Since $\lambda\ge\gamma+\lambda_0\ge\gamma$ on the support of~$g_L$, our assumptions on $\omega$ and Lemma~\ref{lem:trivlem} imply the derivative bounds
\begin{equation}
\label{eq:gederiv2} 
|g_L^{(j)}(\lambda)| \le C\,\gamma^{-j}\,\frac{\lambda}{\sqrt{\lambda^2+\tau}}
\end{equation}
for $j=0,1,2$ and uniformly in $L>0$. By definition of $g_L$, 
\begin{equation}
\label{eq:teil2}
\left|\int_{\lambda_0}^\infty e^{i\phi(\lambda)}\,\frac{\lambda}{\sqrt{\lambda^2+\tau}}(1-\chi(\lambda/L))\omega(\lambda)\,d\lambda
\right| \le a_0\int_{\lambda_0}^{\lambda_0+2\gamma} \frac{\lambda}{\sqrt{\lambda^2+\tau}}\,d\lambda
+  \left|\int_{\lambda_0}^\infty e^{i\phi(\lambda)}\,\,g_L(\lambda)\,d\lambda
\right|.
\end{equation}
The first integral on the right-hand side of \eqref{eq:teil2} is 
\begin{equation}
\label{eq:cutoff1}
\sqrt{(\lambda_0+2\gamma)^2+\tau}-\sqrt{\lambda_0^2+\tau}\le C\,\frac{\gamma(\gamma+\lambda_0)}{\sqrt{\tau}}.
\end{equation}
This bound is wasteful if $\gamma(\gamma+\lambda_0)>\tau$ but that does not concern us. Integrating by parts twice in the integral on the right-hand side of~\eqref{eq:teil2} and repeating the arguments that reduced~\eqref{eq:Ione} to the first term, allows us to bound the final term in~\eqref{eq:teil2} by
\be
&& \int_{\lambda_0+\gamma}^{\sqrt{R}\vee (\lambda_0+\gamma)} \frac{c^{-2}\,R^2}{\lambda^4(\lambda-\lambda_0)^2}\gamma^{-2}
\frac{\lambda}{\sqrt{\tau}}\,d\lambda + c^{-2}\int_{\sqrt{R}\vee (\lambda_0+\gamma)}^\infty \frac{1}{(\lambda-\lambda_0)^2}\,\gamma^{-2}\frac{\lambda}{\sqrt{\lambda^2+\tau}}\,d\lambda \nonumber\\
&\lesssim & c^{-2}\frac{R^2\gamma^{-3}}{(\lambda_0+\gamma)^3\sqrt{\tau}}+ c^{-2} R^{-\half}\gamma^{-2}. 
\label{eq:Rbound}
\ee
One now chooses $\gamma$ in such a way that the first term in \eqref{eq:Rbound} equals the right-hand side
of~\eqref{eq:cutoff1}. This leads to
\begin{equation}
\label{eq:gammadef}
   \gamma = \left\{ \begin{array}{lcr} (R/c)^{\frac14} & \text{if}& (R/c)^{\frac14}>\lambda_0 \\
                                       \frac{\sqrt{R/c}}{\lambda_0} & \text{if}& (R/c)^{\frac14}\le \lambda_0.
                    \end{array}
            \right.
\end{equation}
Note that $R\ge1$ implies that $\gamma\ge c^{-\frac14}$ in the first case, whereas in the second case 
$\gamma\ge c^{-\half}$ because of $\sqrt{R}\ge\lambda_0$. Therefore, the final term in~\eqref{eq:Rbound} is at most~$c^{-\frac32}$. 
Inserting the preceding bounds into~\eqref{eq:teil2} shows that
\begin{equation}
\label{eq:fertig1} 
\left|\int_{\lambda_0}^\infty e^{i\phi(\lambda)}\,\frac{\lambda}{\sqrt{\lambda^2+\tau}}(1-\chi(\lambda/L))\omega(\lambda)\,d\lambda
\right| \le C\,\Bigl[\sqrt{\frac{R}{c\tau}}+c^{-\frac32}\Bigr]
\end{equation}
for both choices of $\gamma$ in \eqref{eq:gammadef}, provided $R\ge1$. Finally, if $R\le1$, then also $\lambda_0\le1$. Hence
\be
&& \left|\int_{\lambda_0}^\infty e^{i\phi(\lambda)}\,\frac{\lambda}{\sqrt{\lambda^2+\tau}}(1-\chi(\lambda/L))\omega(\lambda)\,d\lambda\right| 
\le \left|\int_{\lambda_0}^\infty e^{i\phi(\lambda)}\,\frac{\lambda}{\sqrt{\lambda^2+\tau}}(1-\chi(c\lambda))\omega(\lambda)\,d\lambda\right|  \label{eq:leftover}\\ 
&& \qquad +\left|\int_{\lambda_0}^\infty e^{i\phi(\lambda)}\,\frac{\lambda}{\sqrt{\lambda^2+\tau}}
\chi(c\lambda)(1-\chi(\lambda/L))\omega(\lambda)\,d\lambda\right|  \nonumber\\
&& \lesssim c^{-1}+ c^{-2}\int_{c^{-1}}^\infty (\lambda-\lambda_0)^{-2}\,|g_L''(\lambda)|\,d\lambda \label{eq:newg} 
\ee
where $g_L(\lambda):=  \frac{\lambda}{\sqrt{\lambda^2+\tau}}\chi(c\lambda)(1-\chi(\lambda/L))\omega(\lambda)$. To pass to~\eqref{eq:newg} we used the lower bound on $\phi'$ from~\eqref{eq:R} as well as the usual reduction arguments involving~\eqref{eq:Ione}. Since $|g_L''(\lambda)|\lesssim 1$, say, and $\lambda_0\le1$, \eqref{eq:newg} is at most $c^{-1}$ uniformly in~$L$. Passing to the limit $L\to\infty$ in \eqref{eq:fertig1} and~\eqref{eq:leftover} finishes the proof.
\end{proof}

\noindent We now turn to the contribution of the interval $[0,\lambda_0]$ to the oscillatory integral.

\begin{lemma}
\label{lem:left}
Let $\phi(\lambda),\omega(\lambda)$ be as in the previous lemma and suppose that $\phi$ has
a positive critical point~$\lambda_0$. 
If
\begin{equation}
\label{eq:leftcond}
|\phi'(\lambda)|\ge B\lambda(\lambda_0-\lambda) \text{\ \ for all\ \ }(B\lambda_0)^{-\half}<\lambda<\lambda_0,
\end{equation}
then 
\begin{equation}
\label{eq:links}
\left|\int_0^{\lambda_0} e^{i\phi(\lambda)} \frac{\lambda}{\sqrt{\lambda^2+\tau}}\,\omega(\lambda)\,d\lambda \right| \le C_1\,(B\lambda_0)^{-\half}\,\frac{\lambda_0}{\sqrt{\lambda_0^2+\tau}}.
\end{equation}
The constant $C_1$ only depends on $\omega$.
\end{lemma}
\begin{proof} As in the previous lemma, we will cut off intervals of size $\gamma$ from the endpoints $0$ and~$\lambda_0$ and integrate by parts between $\gamma$ and $\lambda_0-\gamma$. More precisely, we assume that  $(B\la_0)^{-\frac 12}\le \gamma\le\half\lambda_0$ and define
\[ h(\lambda) := \chi(\lambda/\gamma)\chi((\lambda_0-\lambda)/\gamma)\frac{\lambda}{\sqrt{\lambda^2+\tau}}\,\omega(\lambda).
\]
Since $\frac{\lambda}{\sqrt{\lambda^2+\tau}}$ is an increasing function and $\lambda\ge\gamma$ on the support of~$h$, Lemma~\ref{lem:trivlem} and our assumptions on~$\omega$ imply
\begin{equation}
\label{eq:hderiv} 
|h^{(j)}(\lambda)| \le C\,\gamma^{-j}\frac{\lambda_0}{\sqrt{\lambda_0^2+\tau}} 
\end{equation}
for $j=0,1,2$. 
By the definition of $h$, 
\begin{equation}
\label{eq:teil3}
\left| \int_0^{\lambda_0} e^{i\phi(\lambda)}\;\frac{\lambda}{\sqrt{\lambda^2+\tau}}\,\omega(\lambda)\,d\lambda
\right| \le 2a_0\gamma\frac{\lambda_0}{\sqrt{\lambda_0^2+\tau}} +  
\left| \int_0^{\lambda_0} e^{i\phi(\lambda)}\;h(\lambda)\,d\lambda\right|.
\end{equation}
Integrating by parts twice in the last terms yields the upper bound
\begin{equation} 
\label{eq:Itwo}
\int_0^\infty |\phi'(\lambda)|^{-2}\, |h''(\lambda)|\,d\lambda 
+  3\int_0^\infty \Bigl|\frac{\phi''(\lambda)}{\phi'(\lambda)^{3}}\Bigr|\, |h'(\lambda)|\,d\lambda + 
\int_0^\infty \Bigl|\frac{d}{d\lambda}\Bigl[\frac{\phi''(\lambda)}{\phi'(\lambda)^{3}}\Bigr]\Bigr|\, |h(\lambda)|\,d\lambda,
\end{equation}
see \eqref{eq:Ione}.
By \eqref{eq:claim} in Lemma~\ref{lem:lem1}, the absolute values in the third term can be taken outside, which reduces it to the second term in~\eqref{eq:Itwo}. To deal with the second term, recall from Lemma~\ref{lem:lem1} that there can be at most one zero of $\phi''$ on the support of $h$, which we denote by $\lambda_1$ (if there is no such zero, then apply the following argument with~$\lambda_1=0$). Since 
 $\phi''(\lambda)\ge0$ for $\lambda\ge\lambda_1$, $\phi''(\lambda)\le0$ for $\lambda\le\lambda_1$, and $\phi'(\lambda)\le0$ on the support of~$h$,
\begin{equation} \label{eq:secterm}
\int_0^\infty \Bigl|\frac{\phi''(\lambda)}{\phi'(\lambda)^{3}}\Bigr|\, |h'(\lambda)|\,d\lambda = 
- \int_{\lambda_1}^\infty \frac{\phi''(\lambda)}{\phi'(\lambda)^{3}}\, |h'(\lambda)|\,d\lambda + \int_0^{\lambda_1} \frac{\phi''(\lambda)}{\phi'(\lambda)^{3}}\, |h'(\lambda)|\,d\lambda
\end{equation}
In view of the definition of $h$ and the assumed monotonicity of $\omega$ one has $h'=\sum_{i=1}^4 f_i$ where each function $f_i$ has a definite sign. Integrating by parts once one obtains 
\[ 
J := \int_{\lambda_1}^\infty \frac{\phi''(\lambda)}{\phi'(\lambda)^{3}}\, f_i(\lambda)\,d\lambda = -\int_{\lambda_1}^\infty 
\frac{1}{\phi'(\lambda)^{2}}\, f_i'(\lambda)\,d\lambda + 3J - \frac{1}{\phi'(\lambda_1)^{2}}\, f_i(\lambda_1)
\]
and similarly for the second integral in~\eqref{eq:secterm}. Hence the second term in~\eqref{eq:secterm} is no larger than
\begin{equation}
\label{eq:lam1} 
C\sum_{i=1}^4 \left\{\int_0^\infty |\phi'(\lambda)|^{-2}\, |f_i'(\lambda)|\,d\lambda + 
\frac{1}{\phi'(\lambda_1)^{2}}\, |f_i(\lambda_1)|\right\}.
\end{equation}
As $\gamma\le \lambda_1\le \lambda_0-\gamma$, one concludes from \eqref{eq:leftcond} and \eqref{eq:hderiv} with $j=1$ that
\begin{equation}
\label{eq:boundterm}
\sum_{i=1}^4 \frac{1}{\phi'(\lambda_1)^{2}}\, |f_i(\lambda_1)| \lesssim B^{-2}\frac{1}{\lambda_1^2(\lambda_0-\lambda_1)^2}\,\gamma^{-1}\frac{\lambda_0}{\sqrt{\lambda_0^2+\tau}} \lesssim B^{-2}\lambda_0^{-2}\gamma^{-3}\,\frac{\lambda_0}{\sqrt{\lambda_0^2+\tau}}.
\end{equation}
On the other hand, since $f_i'$ and $h''$ both satisfy \eqref{eq:hderiv}, the first term in~\eqref{eq:Itwo} as well as the contribution by the integrals in~\eqref{eq:lam1} does not exceed
\[ B^{-2}\int_{\gamma}^{\lambda_0-\gamma} \frac{1}{\lambda^2(\lambda_0-\lambda)^2}\,\gamma^{-2}\,\frac{\lambda_0}{\sqrt{\lambda_0^2+\tau}}\,d\lambda \lesssim B^{-2}\lambda_0^{-2}\gamma^{-3}\,\frac{\lambda_0}{\sqrt{\lambda_0^2+\tau}},
\] which is the same as \eqref{eq:boundterm}. 
Inserting this bound with the choice of $\gamma=(B\lambda_0)^{-\half}$ into~\eqref{eq:teil3} yields
the desired estimate~\eqref{eq:links}. Recall that this argument assumed that $\gamma\le \half\lambda_0$, i.e.,
$B\gtrsim \lambda_0^{-3}$ by our choice of $\gamma$. However, if $B\lesssim \lambda_0^{-3}$, then
one has
\[ 
\left| \int_0^{\lambda_0} e^{i\phi(\lambda)}\;\frac{\lambda}{\sqrt{\lambda^2+\tau}}\,\omega(\lambda)\,d\lambda
\right| \le \frac{\lambda_0^2}{\sqrt{\lambda^2_0+\tau}} \lesssim (B\lambda_0)^{-\half}\,\frac{\lambda_0}{\sqrt{\lambda_0^2+\tau}},
\]
and we are done.
\end{proof}

\noindent In the next lemma we obtain general lower bounds on $\phi'$ and $\phi''$.

\begin{lemma}
\label{lem:genlower}
Let $\phi$ and $\lambda_0$ be as in Lemma~\ref{lem:lem1}.
\begin{enumerate}
\item 
If $\lambda_0>0$ exists, then
\be
\label{eq:rightder}
\phi'(\lambda) &\ge& \half(1+\sqrt{\sigma_1})^{-1}(\lambda-\lambda_0)\text{\ \ for all\ \ }\lambda\ge\lambda_0+2\sigma_1^{\frac14} \\
|\phi'(\lambda)| &\ge& \bigl[\lambda_0(1+\sqrt{\sigma_1})\bigr]^{-1}\lambda(\lambda-\lambda_0)\text{\ \ for all\ \ }\sigma_1^{\frac14}\le\lambda<\lambda_0. \label{eq:leftder}
\ee
If $\lambda_0>0$ does not exist, then 
\eqref{eq:rightder} holds for all $\lambda\ge0$.
\item
Now suppose $\lambda_0$ exists. Then 
\be
\phi'(\lambda) &\ge& \phi''(\lambda_0) (\lambda-\lambda_0) \text{\ \ for all\ \ }\lambda\ge\lambda_0\label{eq:right2'} \\
|\phi'(\lambda)| &\ge& \half\phi''(\lambda_0)\frac{\lambda}{\lambda_0} (\lambda_0 -\lambda)\text{\ \ for all\ \ }0\le \lambda\le\lambda_0. \label{eq:left2'}
\ee
If $\lambda_0\ge\sqrt{\sigma_k}$, then $\phi''(\lambda_0)\ge\frac{b_k}{4m \max_\ell b_\ell}$. 
\end{enumerate}
\end{lemma}
\begin{proof} By \eqref{eq:zwei'}, for all $\lambda\ge\lambda_0$ 
(or $\lambda\ge0$ if $\lambda_0$ does not exist),
\be
\phi''(\lambda) &\ge& 1-\sum_{j=1}^m \frac{b_j}{\sqrt{\lambda^2+\sigma_j}} + \frac{\lambda^2}{\lambda^2+\sigma_1} \sum_{j=1}^m \frac{b_j}{\sqrt{\lambda^2+\sigma_j}} \nonumber\\
&\ge& 1 - \frac{\sigma_1}{\lambda^2+\sigma_1} \sum_{j=1}^m \frac{b_j}{\sqrt{\lambda^2+\sigma_j}}
\ge 1 - \frac{\sigma_1}{\lambda^2+\sigma_1} = \frac{\lambda^2}{\lambda^2+\sigma_1}. 
\label{eq:secder}
\ee
For the second inequality sign in \eqref{eq:secder} we used the fact that if $\lambda_0$ exists, then
\[ \sum_{j=1}^m \frac{b_j}{\sqrt{\lambda^2+\sigma_j}} \le 1 \text{\ \ for all\ \ }\lambda\ge\lambda_0,\]
and for all $\lambda\ge0$ if $\lambda_0$ does not exist. See~\eqref{eq:lam0def} and~\eqref{eq:lam0ne}. 
Hence $\phi''(\lambda)\ge (1+\sqrt{\sigma_1})^{-1}$ provided $\lambda\ge\sigma_1^{\frac14}$.
Recall from Lemma~\ref{lem:lem1} that $\phi''(\lambda)\ge0$ for all $\lambda\ge\lambda_0$ (or for all $\lambda\ge0$ in case a positive critical point does not exist).
Setting $\lambda_0=0$ if a positive critical point does not exist, it follows that
\be
\phi'(\lambda) &\ge& \int_{\lambda_0}^\lambda \phi''(s)\,ds \ge \int_{\lambda_0+\sigma_1^{\frac14}}^\lambda \phi''(s)\,ds \nonumber\\
&\ge& (1+\sqrt{\sigma_1})^{-1} (\lambda-\lambda_0-\sigma_1^{\frac14}) \ge
\half (1+\sqrt{\sigma_1})^{-1} (\lambda-\lambda_0) \nonumber
\ee
for all $\lambda\ge \lambda_0+2\sigma_1^{\frac14}$, as claimed. 

\noindent Now assume that $\lambda_0>0$ exists. If $\lambda<\lambda_0$, then~\eqref{eq:fund} implies that
\be
|\phi'(\lambda)| &=& \lambda\int_{\lambda}^{\lambda_0} \sum_{j=1}^m \frac{b_js}{(s^2+\sigma_j)^{\frac32}}\,ds
\ge \lambda \int_{\lambda}^{\lambda_0} \sum_{j=1}^m \frac{b_j}{\sqrt{s^2+\sigma_j}}\,\frac{s^2}{s^2+\sigma_1}\,\frac{ds}{s} \nonumber\\
&\ge& \frac{\lambda^3}{\lambda^2+\sigma_1} \int_{\lambda}^{\lambda_0} \sum_{j=1}^m \frac{b_j}{\sqrt{\lambda^2_0+\sigma_j}}\,\frac{ds}{\lambda_0} = \frac{\lambda^2}{\lambda^2+\sigma_1} \frac{\lambda}{\lambda_0} (\lambda_0-\lambda),
\label{eq:linksabl}
\ee 
which is \eqref{eq:leftder} provided $\lambda\ge \sigma_1^{\frac14}$. 
To pass to \eqref{eq:linksabl} we used that $\frac{s^2}{s^2+\sigma_1}$ is increasing, whereas the final expression uses \eqref{eq:lam0def}.

\noindent The proof of \eqref{eq:right2'} is an immediate consequence of the fact that $\phi''(\lambda)$ is increasing. For~\eqref{eq:left2'} one again uses~\eqref{eq:fund}. Indeed, for $0\le\lambda\le\lambda_0$, 
\be
|\phi'(\lambda)| &\ge& \lambda\int_\lambda^{\lambda_0} \sum_{j=1}^m \frac{b_js}{(\lambda_0^2+\sigma_j)^{\frac32}}\,ds \nonumber\\
&\ge& \lambda \int_\lambda^{\lambda_0} \frac{\phi''(\lambda_0)}{\lambda_0^2}\,s\,ds \ge \half 
\frac{\phi''(\lambda_0)}{\lambda_0^2}\la (\lambda_0^2-\lambda^2), \nonumber
\ee
as claimed. Assume $\lambda_0\ge\sqrt{\sigma_k}$. By \eqref{eq:lam0def} one has 
\[ \max_{1\le\ell\le m} \frac{b_\ell}{\sqrt{\lambda_0^2+\sigma_\ell}} \ge \frac{1}{m}
\text{\ \ and thus also\ \ } m\max_{\ell} b_\ell \ge \lambda_0\]
Furthermore, \eqref{eq:zwei'} implies that
\[ \phi''(\lambda_0) \ge \lambda_0^2 \frac{b_k}{(\lambda_0^2+\sigma_k)^{\frac32}} \ge 
\frac{b_k}{4m\max_{\ell} b_\ell},\]
and the lemma follows.
\end{proof}

\begin{proof}[Proof of Lemma~\ref{lem:osc1}]
The bound $(1+\sigma_1)^{\frac14}$ on the right-hand side of \eqref{eq:bound1} follows from \eqref{eq:leftder}, \eqref{eq:rightder}, \eqref{eq:rechts1} with $A\asymp (1+\sigma_1)^{-\frac12}$, 
and~\eqref{eq:links} with $B\asymp [\lambda_0^2(1+\sigma_1)]^{-\frac12}$. As for the other bound in~\eqref{eq:bound1}, fix some $1\le k\le m$. If $\lambda_0>0$ exists and $\lambda_0 > \sqrt{\sigma_k}$, then one applies \eqref{eq:rechts1} with $A\asymp \frac{b_k}{m\max_\ell b_\ell}$ and \eqref{eq:links} with $B\asymp \lambda_0^{-1} A$, see the second part of Lemma~\ref{lem:genlower}. Observe that this leads to the bound 
\[ m^{\half}\sqrt{\frac{\max_\ell b_\ell}{b_k}},\]
which is smaller than the one in~\eqref{eq:bound1}.
It therefore remains to consider the case  $\lambda_0 \le \sqrt{\sigma_k}$, where we set $\lambda_0=0$
if a positive critical point does not exist.
Choose a maximal $n$ so that $k\le n\le m$ and for which
$\sigma_{n+1}\le \lambda_0^2\le \sigma_n$. Here we have set $\sigma_{m+1}:=0$. 
If $\lambda_0>0$ exists, then \eqref{eq:lam0def} implies that
\begin{equation}
\label{eq:teilung}
1=\sum_{j=1}^m \frac{b_j}{\sqrt{\lambda_0^2+\sigma_j}} \asymp \frac{1}{\lambda_0}\sum_{j>n} b_j
+ \sum_{j\le n}\frac{b_j}{\sqrt{\sigma_j}}.
\end{equation}  
If  $\lambda_0=0$ (which includes the case that the positive critical point does not exist), then
\begin{equation} 
\nonumber
1\ge \sum_{j=1}^m \frac{b_j}{\sqrt{\sigma_j}},
\end{equation}
see \eqref{eq:lam0ne}. If in fact
\[ \sum_{j=1}^m \frac{b_j}{\sqrt{\sigma_j}}\le\half,\]
then \eqref{eq:ein'} implies that $\phi'(\lambda)\ge\half\lambda$ for all $\lambda\ge0$. 
Hence~\eqref{eq:rechts1} with $A=\half$ yields a~$O(1)$ bound for the entire oscillatory integral.
Hence we can assume without loss of generality that 
\begin{equation} 
\label{eq:ne2}
1\ge \sum_{j=1}^m \frac{b_j}{\sqrt{\sigma_j}}\ge\half.
\end{equation}
This again implies that the right-hand side of \eqref{eq:teilung} is comparable to~$1$. 
Suppose that the first sum on the right-hand side of~\eqref{eq:teilung}
dominates. Then
\[ \phi''(\lambda_0) = \sum_{j=1}^m \frac{b_j\lambda_0^2}{(\lambda_0^2+\sigma_j)^{\frac32}}
\gtrsim \lambda_0^{-1}\,\sum_{j>n} b_j \gtrsim 1.\]
In view of the second part of Lemma~\ref{lem:genlower} one can apply~\eqref{eq:rechts1} and~\eqref{eq:links} with
$A\asymp B\asymp 1$, which leads to a~$O(1)$ bound in this case. Now suppose that the second sum 
in~\eqref{eq:teilung} is the larger one, i.e., 
\[ 1 \asymp  \sum_{j\le n}\frac{b_j}{\sqrt{\sigma_j}}.\]
Then there is some $p\le n$ for which
\begin{equation}
\label{eq:bpsigp}
 b_p \gtrsim \frac{\sqrt{\sigma_p}}{m}.
\end{equation}
By \eqref{eq:zwei'} therefore
\begin{equation}
\label{eq:2'inter}
\phi''(\lambda_0)\ge \lambda_0^2\frac{b_p}{(\lambda_0^2+\sigma_p)^{\frac32}} \gtrsim \lambda_0^2 \frac{b_p}{\sigma_p^{\frac32}} \gtrsim \frac{\lambda_0^2}{m\sigma_p}.
\end{equation}
Combining \eqref{eq:2'inter}, \eqref{eq:left2'}, and~\eqref{eq:links} with $B\asymp \lambda_0^{-1}\phi''(\lambda)$ leads to an upper bound 
\begin{equation}
\label{eq:endlich}
\left|\int_0^{\lambda_0} e^{i\phi(\lambda)}\;\frac{\lambda}{\sqrt{\lambda^2+\sigma_k}}\,\omega(\lambda)\,d\lambda \right| \lesssim \frac{\sqrt{m\sigma_p}}{\lambda_0}\frac{\lambda_0}{\sqrt{\sigma_k}} \le m^{\frac32}\frac{b_p}{b_k},
\end{equation}
which is precisely the second bound in \eqref{eq:bound1}. To pass to the final expression in \eqref{eq:endlich}
one uses~\eqref{eq:bpsigp} as well as $b_k\le\sqrt{\sigma_k+\lambda_0^2}\lesssim \sqrt{\sigma_k}$, which follows
immediately from~\eqref{eq:lam0def} or~\eqref{eq:ne2}. We will now estimate the oscillatory integral on the interval~$[\lambda_0,\infty)$ by means of~\eqref{eq:rechts2}. Recall that $\lambda_0\leq\sqrt{\sigma_p}$. If $\lambda_0\le\lambda\le 2\sqrt{\sigma_p}$, then by~\eqref{eq:fund} 
\begin{equation}
\label{eq:stuck1}
\phi'(\lambda) \ge \lambda \int_{\lambda_0}^\lambda \frac{b_ps}{(s^2+\sigma_p)^{\frac32}}\,ds
\gtrsim \lambda \int_{\lambda_0}^\lambda \frac{b_ps}{\sigma_p^{\frac32}}\,ds \ge \frac{\lambda^2}{m\sigma_p}(\lambda-\lambda_0),
\end{equation}
whereas for $2\sqrt{\sigma_p}\le \lambda$, 
\begin{equation}
\label{eq:stuck2}
\phi'(\lambda) \ge \lambda \int_{\lambda_0}^{2\sqrt{\sigma_p}} \frac{b_ps}{(s^2+\sigma_p)^{\frac32}}\,ds
\gtrsim \frac{\lambda}{m}.
\end{equation}
It follows from \eqref{eq:stuck1} and \eqref{eq:stuck2} that the condition \eqref{eq:R} in 
Lemma~\ref{lem:right} holds with~$R=\sigma_p$ and $c\asymp m^{-1}$. 
Hence~\eqref{eq:rechts2} implies that 
\begin{equation}
\label{eq:endlich2}
\left|\int_{\lambda_0}^\infty e^{i\phi(\lambda)}\;\frac{\lambda}{\sqrt{\lambda^2+\sigma_k}}\,\omega(\lambda)\,d\lambda \right| \lesssim m^{\frac32}+\sqrt{\frac{m\sigma_p}{\sigma_k}} \lesssim m^{\frac32}\frac{b_p}{b_k},
\end{equation}
as desired. To pass to the final inequality we used \eqref{eq:bpsigp} and that $b_k\le \sqrt{\sigma_k+\lambda_0^2}\lesssim \sqrt{\sigma_k}$. The lemma is proven in all cases.
\end{proof}

\begin{remark} It is possible to obtain the same bounds as in Lemma~\ref{lem:osc1} by means of a single integration
by parts. Whereas this would lead to some minor simplifications, as avoiding~\eqref{eq:claim}, we believe that the
previous argument involving two integrations by parts might have some interest in its own right. For example, 
it applies to potentials whose first derivative does not decay (but is bounded, say). In the next section we
show how to implement an argument based on a single integration by parts.
\end{remark}

\section{Estimates for oscillatory integrals II}\label{sec:time3}

The purpose of this section is to prove Lemmas \ref{lem:singlem1} and  \ref{lem:singlem2} below. 
These lemmas are needed to control the oscillatory integrals arising in~\eqref{Mdl} and~\eqref{Mbardl}. 
The arguments are similar to the ones from
the previous section, the main difference being the singularity at the point~$\sqrt{\rho_\ell}$.
To overcome it, we change variables $u=\sqrt{\rho_\ell-\lambda^2}$. This leads to a new class of phase functions~$\psi(u)$.
We present some useful properties of this class in the following lemma.

\begin{lemma}
\label{lem:uphase}
Let $\psi(u)=\half u^2 +\sum_{j=1}^m b_j\,\sqrt{\tau_j-u^2}$ where $\tau_1\ge\tau_2\ge\ldots\ge \tau_m>0$
and $b_j>0$ are arbitrary. For any choice of these parameters the following properties hold:
\begin{enumerate}
\item
$\psi(u)$ has at most one critical point $u\in(0,\sqrt{\tau_m})$. This point will always be denoted by~$u_0$.
If $u_0<u<\sqrt{\tau_m}$, then~$\psi'(u)<0$ and $\psi''(u)<0$, whereas $\psi'(u)>0$ on~$(0,u_0)$. 
There is at most one point $u_1\in(0,u_0)$ such that $\psi''(u_1)=0$. 
\item 
If $u_0>0$ exists, then
\begin{equation}
\label{eq:ufund}
\psi'(u) = -u\int_{u_0}^u \sum_{j=1}^m \frac{b_j\,s}{(\tau_j-s^2)^{\frac32}} \,ds
\end{equation}
for any $0<u<\sqrt{\tau_m}$, whereas if $u_0>0$ does not exist, then
\begin{equation}
\label{eq:ufundne}
|\psi'(u)| \ge u\int_{0}^u \sum_{j=1}^m \frac{b_j\,s}{(\tau_j-s^2)^{\frac32}} \,ds
\end{equation}
for all $0<u<\sqrt{\tau_m}$. 
\item
If $\sum_{j=1}^m \frac{b_j}{\sqrt{\tau_j}}\le\half$, then $\psi'(u)\ge\frac14 u$ for all~$0<u<\half\sqrt{\tau_m}$.
If $\sum_{j=1}^m \frac{b_j}{\sqrt{\tau_j}}\ge 2$, then $\psi'(u)\le -u$ for all~$0<u<\sqrt{\tau_m}$.
\end{enumerate}
\end{lemma}
\begin{proof}
The first assertion follows from the explicit expressions
\be
\frac{\psi'(u)}{u} &=& 1-\sum_{j=1}^m \frac{b_j}{\sqrt{\tau_j-u^2}} \nonumber\\
\psi''(u) &=&  \frac{\psi'(u)}{u} - \sum_{j=1}^m \frac{b_j\,u^2}{(\tau_j-u^2)^{\frac32}}= 1-\sum_{j=1}^m\frac{b_j\tau_j}{(\tau_j-u^2)^{\frac32}}. \nonumber
\ee
For the second assertion, viz.~\eqref{eq:ufund} and \eqref{eq:ufundne}, integrate out
\[ \Bigl(\frac{\psi'(u)}{u}\Bigr)' = -\sum_{j=1}^m \frac{b_j\,u}{\sqrt{\tau_j-u^2}}.\]
If $u_0>0$ does not exist, then use that $\frac{\psi'(u)}{u}\le0$ for all $u\ge0$. 
If $\sum_{j=1}^m \frac{b_j}{\sqrt{\tau_j}}\le\half$, then for all~$0<u<\half\sqrt{\tau_m}$
\[ \sum_{j=1}^m \frac{b_j}{\sqrt{\tau_j-u^2}} \le \sqrt{\frac43}\sum_{j=1}^m \frac{b_j}{\sqrt{\tau_j}}\le \frac{1}{\sqrt{3}}\]
which implies that $\psi'(u)\ge u(1-1/\sqrt{3})\ge u/4$ for those $u$. 
If on the other hand $\sum_{j=1}^m \frac{b_j}{\sqrt{\tau_j}}\ge 2$, then 
\[ -\psi'(u)\ge u\,\sum_{j=1}^m \frac{b_j}{\sqrt{\tau_j}} - u\ge u\]
for all $0<u<\sqrt{\tau_m}$. 
\end{proof}

\noindent In what follows let $\chi$ be a smooth non-decreasing function with $\chi(\lambda)=0$ for $\lambda\le\frac14$ 
and $\chi(\lambda)=1$ for $\lambda\ge\frac12$. Furthermore, we shall use the notation $\chi_E$ for the indicator 
of a set~$E$. 

\begin{lemma} 
\label{lem:urechts}
Let $\psi(u)=\half u^2 +\sum_{j=1}^m b_j\,\sqrt{\tau_j-u^2}$ where $\tau_1\ge\tau_2\ge\ldots\ge \tau_m>0$
and $b_j>0$. Assume that $\omega$ is a differentiable function such that
\[ |\omega^{(j)}(u)| \le a_0\,u^{-j}\]
for all $u>0$ and $j=0,1$. 
\begin{enumerate}
\item 
Suppose that $A>0,u_0\ge 0$ (where $u_0$ is not necessarily a critical point of $\psi$) have the property that  
$|\psi'(u)|\ge Au^2(u-u_0)$ and $\psi''(u)\le0$ for all $\sqrt{\tau_m} > u\ge u_0$. Let $\tau\ge\tau_m$. Then
\be
\label{eq:uright1}
\left| \int_{u_0}^\infty e^{i\psi(u)}\;\Bigl(1-\chi\Bigl(\frac{u}{\sqrt{\tau_m}}\Bigr)\Bigr)\omega(u)\frac{u}{\sqrt{\tau}}\,du
\right| 
&\le&\, C_0\,\min(A^{-\frac14},(A\tau)^{-\half}) \\
\left| \int_{u_0}^\infty e^{i\psi(u)}\;\Bigl(1-\chi\Bigl(\frac{u}{\sqrt{\tau_m}}\Bigr)\Bigr)\omega(u)\,du
\right| 
&\le&\, C_0\,A^{-\frac14}. \label{eq:uright2}
\ee
\item
Suppose that $|\psi'(u)|\ge\half u$ and $\psi''(u)\le0$ for all $\half\sqrt{\tau_m}>u\ge0$. Then
\begin{equation}
\label{eq:O1uright}
\left| \int_{0}^\infty e^{i\psi(u)}\;\Bigl(1-\chi\Bigl(\frac{u}{\sqrt{\tau_m}}\Bigr)\Bigr)\omega(u)\,du
\right| 
\le\, C_0.
\end{equation}
\end{enumerate}
In both cases the constant $C_0$ only depends on $a_0$.
\end{lemma}
\begin{proof} 
By our assumptions the function
\[ h(u):= \Bigl(1-\chi\Bigl(\frac{u}{\sqrt{\tau_m}}\Bigr)\Bigr)\omega(u)\frac{u}{\sqrt{\tau}}\]
satisfies the derivative estimates
\[ |h^{(j)}(u)|\lesssim u^{-j}\frac{u}{\sqrt{\tau}} \text{\ \ for all\ \ }u>0 \text{\ and\ }j=0,1.\]
Thus
\be
&& \left| \int_{u_0}^\infty e^{i\psi(u)}\;\Bigl(1-\chi\Bigl(\frac{u}{\sqrt{\tau_m}}\Bigr)\Bigr)\omega(u)\frac{u}{\sqrt{\tau}}\,du
\right| \nonumber\\
&\lesssim& \frac{\gamma(u_0+\gamma)}{\sqrt{\tau}} + \left| \int_{u_0}^\infty e^{i\psi(u)}\; h(u)\chi((u-u_0)/\gamma)\,du\right| \nonumber\\
&\lesssim& \frac{\gamma(u_0+\gamma)}{\sqrt{\tau}} + \int_{u_0}^\infty \frac{|\psi''(u)|}{\psi'(u)^2}\,h(u)\chi((u-u_0)/\gamma)\,du 
+ \int_{u_0}^\infty \frac{1}{|\psi'(u)|}\,|[h(u)\chi((u-u_0)/\gamma)]'|\,du.  \label{eq:ubreak}
\ee
The second term here can be reduced to the third by means of a further integration by parts:
\be 
&& \int_{u_0}^\infty \frac{\psi''(u)}{\psi'(u)^2}\,h(u)\chi((u-u_0)/\gamma)\,du = 
-\int_{u_0}^\infty \psi'(u)\Bigl[\frac{1}{\psi'(u)^2}\,h(u)\chi((u-u_0)/\gamma)\Bigr]'\,du \nonumber\\
&& = 2\int_{u_0}^\infty \frac{\psi''(u)}{\psi'(u)^2}\,h(u)\chi((u-u_0)/\gamma)\,du 
- \int_{u_0}^\infty \frac{1}{\psi'(u)}\,[h(u)\chi((u-u_0)/\gamma)]'\,du. \nonumber
\ee
The final integral in \eqref{eq:ubreak} is no larger than
\begin{equation}
\label{eq:ur2}
\int_{u_0+\gamma}^\infty \frac{1}{A\,u^2(u-u_0)}\Bigl[\gamma^{-1}\frac{u_0+\gamma}{\sqrt{\tau}}\chi_{[\frac14\gamma\le u-u_0\le \frac12\gamma]}+\frac{1}{\sqrt{\tau}}\Bigr]\,du 
\lesssim \frac{1}{A(u_0+\gamma)\gamma}\frac{1}{\sqrt{\tau}},
\end{equation}
where $\chi$  refers to a usual indicator function.
As always, we now choose $\gamma$ so that the first term in \eqref{eq:ubreak} equals the third. This leads to the 
estimate
\[ \eqref{eq:ubreak} \lesssim (A\tau)^{-\half}.\]
To obtain the $A^{-\frac14}$ bound in \eqref{eq:uright1}, one includes the factor $\frac{u}{\sqrt{\tau}}$ into the function~$\omega$. Since $u\le\sqrt{\tau}$ on the support of~$h$ this can be done without violating the conditions on~$\omega$. The same arguments as in the previous case lead to
\be
\left| \int_{u_0}^\infty e^{i\psi(u)}\;\Bigl(1-\chi\Bigl(\frac{u}{\sqrt{\tau_m}}\Bigr)\Bigr)\omega(u)\frac{u}{\sqrt{\tau}}\,du
\right| &\lesssim& \gamma + \int_{u_0+\gamma}^\infty \frac{1}{Au^2(u-u_0)}\,\gamma^{-1}\,du \nonumber\\
&\lesssim& \gamma + A^{-1}\gamma^{-3}. \nonumber
\ee
Setting $\gamma=A^{-\frac14}$ finishes the proof of both \eqref{eq:uright1} and \eqref{eq:uright2}.
The proof of \eqref{eq:O1uright} is similar. More precisely, by the same arguments involving \eqref{eq:ubreak} one obtains
\be
&& \left| \int_{0}^\infty e^{i\psi(u)}\;\Bigl(1-\chi\Bigl(\frac{u}{\sqrt{\tau_m}}\Bigr)\Bigr)\omega(u)\,du
\right| 
\lesssim 1  + \left| \int_{0}^\infty e^{i\psi(u)}\; h(u)\chi(u)\,du\right| \nonumber \\
&\lesssim& 1 + \int_{0}^\infty \frac{1}{\psi'(u)} |[h(u)\chi(u)]'|\,du 
\lesssim 1 + \int_{\frac14}^\infty \frac{1}{u}\Bigl[\chi_{[\frac14,\frac12]}(u)+\frac{1}{u}\Bigr]\,du \lesssim 1 
\nonumber
\ee
and we are done.
\end{proof}

\begin{lemma}
\label{lem:ulinks}
Let $\psi$ and $\omega$ be as above. Suppose that there are constants $A$ and $u_0$ so that 
\[ \psi'(u)\gtrsim Au_0\,u(u-u_0) \text{\ \ for all \ \ } 0\le u\le u_0.\]
Then for any $\tau$ so that $\tau\ge\tau_m$, 
\be
\label{eq:uleft1}
\left| \int_0^{u_0} e^{i\psi(u)}\;\Bigl(1-\chi\Bigl(\frac{u}{\sqrt{\tau_m}}\Bigr)\Bigr)\omega(u)\frac{u}{\sqrt{\tau}}\,du
\right| 
&\le& C_0\,\min(A^{-\frac14},(A\tau)^{-\half}) \\ 
\label{eq:uleft2}
\left| \int_0^{u_0} e^{i\psi(u)}\;\Bigl(1-\chi\Bigl(\frac{u}{\sqrt{\tau_m}}\Bigr)\Bigr)\omega(u)\,du
\right| 
&\le& C_0\,A^{-\frac14}. 
\ee
\end{lemma}
\begin{proof} 
As in the previous proof we define
\[ h(u):= \Bigl(1-\chi\Bigl(\frac{u}{\sqrt{\tau_m}}\Bigr)\Bigr)\omega(u)\frac{u}{\sqrt{\tau}}.\]
Then the derivative estimates 
\[ |h^{(j)}(u)|\lesssim u^{-j}\frac{u}{\sqrt{\tau}}\]
hold for all $u>0$ and $j=0,1$.
Let  $0<\gamma<u_0$ be arbitrary but fixed. Introducing the usual cut-off functions at $0$ and $u_0$ gives
\be
&& \left| \int_0^{u_0} e^{i\psi(u)}\;\Bigl(1-\chi\Bigl(\frac{u}{\sqrt{\tau_m}}\Bigr)\Bigr)\omega(u)\frac{u}{\sqrt{\tau}}\,du
\right| \nonumber\\
&\lesssim& \gamma \frac{u_0}{\sqrt{\tau}} + \left| \int_0^{u_0} e^{i\psi(u)}\;\chi(u/\gamma)\chi((u_0-u)/\gamma)\,h(u)\,du 
\right| \nonumber \\
&\lesssim& \gamma \frac{u_0}{\sqrt{\tau}} + \int_0^{u_0} \frac{|\psi''(u)|}{\psi'(u)^2} \chi(u/\gamma)\chi((u_0-u)/\gamma)h(u)\,du + 
\nonumber\\
&& \qquad + \int_0^{u_0} \frac{1}{\psi'(u)} |[\chi(u/\gamma)\chi((u_0-u)/\gamma)h(u)]'|\,du.
\label{eq:ubreak2}
\ee 
Suppose that $\psi''(u_1)=0$ for some $u_1\in\supp(h)$. If this $u_1$ does not exist, then apply the following
arguments with $u_1=0$. Integrating by parts in the first integral in~\eqref{eq:ubreak2} yields
\be
&& \int_{u_1}^{u_0} \frac{\psi''(u)}{\psi'(u)^2}\, \chi(u/\gamma)\chi((u_0-u)/\gamma)h(u)\,du =
2\int_{u_1}^{u_0} \frac{\psi''(u)}{\psi'(u)^2}\, \chi(u/\gamma)\chi((u_0-u)/\gamma)h(u)\,du \nonumber \\
&& - \int_{u_1}^{u_0} \frac{1}{\psi'(u)}\, [\chi(u/\gamma)\chi((u_0-u)/\gamma)h(u)]'\,du - \frac{1}{\psi'(u_1)}  \chi(u_1/\gamma)\chi((u_0-u_1)/\gamma)h(u_1) 
\label{eq:ubdaryterm}
\ee
and similarly for the integral over the interval $[0,u_1]$. It follows that the estimate in~\eqref{eq:ubreak2} reduces to
\begin{equation}
\label{eq:reduceu}
\gamma \frac{u_0}{\sqrt{\tau}} + \int_0^{u_0} \frac{1}{\psi'(u)} |[\chi(u/\gamma)\chi((u_0-u)/\gamma)h(u)]'|\,du
+ \frac{1}{\psi'(u_1)}  \chi(u_1/\gamma)\chi((u_0-u_1)/\gamma)h(u_1).
\end{equation}
The final term in \eqref{eq:reduceu} is no larger than 
\begin{equation}
\label{eq:bdryest}
\frac{1}{A\,u_0u_1(u_0-u_1)} \frac{u_1}{\sqrt{\tau}} \lesssim \frac{1}{A\,\gamma u_0\sqrt{\tau}},
\end{equation}
since we can assume that $u_1\in\supp(h)$. On the other hand, the integral in~\eqref{eq:reduceu} is bounded by
\be
&& \int_0^{u_0} \frac{1}{\psi'(u)} |[\chi(u/\gamma)\chi((u_0-u)/\gamma)h(u)]'|\,du \nonumber \\
&& \lesssim \int_{\gamma/4}^{u_0-\gamma/4} \frac{1}{Au_0\, u(u_0-u)}\, \Bigl[ \gamma^{-1}\chi_{[\frac{\gamma}{4},\frac{\gamma}{2}]}(u)\,\frac{u}{\sqrt{\tau}} + \gamma^{-1}\chi_{[\frac{\gamma}{4},\frac{\gamma}{2}]}(u_0-u)\,\frac{u}{\sqrt{\tau}} + \frac{1}{\sqrt{\tau}}\Bigr]\,du \nonumber \\
&& \lesssim \frac{1}{A\gamma u_0\sqrt{\tau}}.\label{eq:2nduint}
\ee 
We now choose $\gamma$ so that the first term in \eqref{eq:reduceu} equals the estimates from~\eqref{eq:bdryest} and~\eqref{eq:2nduint}.
This leads to $\gamma=A^{-\half}u_0^{-1}$ which in turn yields the estimate~$(A\tau)^{-\half}$ on the oscillatory integral over~$[0,u_0]$. Recall that this required the condition $\gamma\le u_0$, i.e., $u_0\ge A^{-\frac14}$. If $u_0\le A^{-\frac14}$, then the oscillatory integral over $[0,u_0]$ is clearly no larger than
\[ \int_0^{u_0} \frac{u}{\sqrt{\tau}}\,du \lesssim \frac{u_0^2}{\sqrt{\tau}} \lesssim (A\tau)^{-\half},\]
as claimed.

To obtain the $A^{-\frac14}$ estimate in \eqref{eq:uleft1} and \eqref{eq:uleft2} we incorporate the $\frac{u}{\sqrt{\tau}}$--factor into the function~$\omega$. Since $\frac{u}{\sqrt{\tau}}\leq1$ on the support of the integrand in~\eqref{eq:uleft1} by our assumption $\tau_m\le \tau$, this does not violate the conditions on~$\omega$. Let again $0\le \gamma\le u_0$ be arbitary but fixed. With
\[ h(u):= \Bigl(1-\chi\Bigl(\frac{u}{\sqrt{\tau_m}}\Bigr)\Bigr)\omega(u) \]
the arguments leading to \eqref{eq:ubreak2} and~\eqref{eq:reduceu} yield in this context
\be
&& \left| \int_0^{u_0} e^{i\psi(u)}\;\Bigl(1-\chi\Bigl(\frac{u}{\sqrt{\tau_m}}\Bigr)\Bigr)\omega(u)\,du
\right| \nonumber\\
&\lesssim& \gamma  + \int_0^{u_0} \frac{1}{\psi'(u)} |[\chi(u/\gamma)\chi((u_0-u)/\gamma)h(u)]'|\,du
+ \frac{1}{\psi'(u_1)}  \chi(u_1/\gamma)\chi((u_0-u_1)/\gamma)h(u_1). \label{eq:sameu}
\ee
The final term in \eqref{eq:sameu} is no larger than
\begin{equation}
\label{eq:bdryest2}
\frac{1}{A\,u_0u_1(u_0-u_1)}  \lesssim \frac{1}{A\,\gamma u_0^2} \lesssim \gamma^{-3}\, A^{-1},
\end{equation}
whereas the integral in \eqref{eq:sameu} is at most
\[ 
\int_{\gamma/4}^{u_0-\gamma/4} \frac{1}{Au_0\, u(u_0-u)}\, \Bigl[ \gamma^{-1}\chi_{[\frac{\gamma}{4},\frac{\gamma}{2}]}(u) + 
\gamma^{-1}\chi_{[\frac{\gamma}{4},\frac{\gamma}{2}]}(u_0-u) + \frac{1}{u}\Bigr]\,du 
 \lesssim \frac{1}{A\gamma u_0^2},
\]
which is the same as \eqref{eq:bdryest2}. In view of~\eqref{eq:sameu} and~\eqref{eq:bdryest2}, choosing $\gamma= A^{-\frac14}$ leads to the estimate of $\gamma=A^{-\frac14}$ for the oscillatory integral over~$[0,u_0]$, as desired. Recall that we made the assumption that $\gamma\le u_0$, i.e., $A^{-\frac14} \le  u_0$. If $A^{-\frac14}\gtrsim u_0$, then the oscillatory integral over~$[0,u_0]$ is trivially bounded by~$u_0\lesssim A^{-\frac14}$, as desired.
\end{proof}

\noindent We can now state a bound on oscillatory integrals involving the phase function $\psi(u)$. The reader should note the similarity with
Lemma~\ref{lem:osc1} above. 

\begin{cor}
\label{cor:uosc1}
Let $\psi(u)=\half u^2 +\sum_{j=1}^m b_j\,\sqrt{\tau_j-u^2}$ where $\tau_1\ge\tau_2\ge\ldots\ge \tau_p>0$
and $b_j>0$. Assume that $\omega$ is a differentiable function such that
\[ |\omega^{(j)}(u)| \le a_0\,u^{-j}\]
for all $u>0$ and $j=0,1$. 
Then for any $1\le k\le m$, 
\begin{equation}
\label{eq:umain1}
\left| \int_{0}^\infty e^{i\psi(u)}\;\Bigl(1-\chi\Bigl(\frac{u}{\sqrt{\tau_m}}\Bigr)\Bigr)\omega(u)\frac{u}{\sqrt{\tau_k}}\,du
\right| 
\le\, C_0\,\min\bigl(\tau_1^{\frac14},m^{\frac32}\,b_k^{-1}\max_{\ell} b_\ell\bigr)
\end{equation}
\end{cor}
\begin{proof}
If either $\sum_{j=1}^m \frac{b_j}{\sqrt{\tau_j}}\le\half$ or $\sum_{j=1}^m \frac{b_j}{\sqrt{\tau_j}}\ge 2$, then the final assertion
of Lemma~\ref{lem:uphase} and~\eqref{eq:O1uright} of Lemma~\ref{lem:urechts} imply that the left-hand side of~\eqref{eq:umain1}
is~$O(1)$. We can therefore assume that 
\[ \sum_{j=1}^m \frac{b_j}{\sqrt{\tau_j}}\asymp 1.\]
This implies that there is some $p$ so that $m b_p\gtrsim \sqrt{\tau_m}$, whereas $\sqrt{\tau_k}\lesssim b_k$ for
all~$k$. Let 
\[ A:= \sum_{j=1}^m \frac{b_j}{{\tau_j}^{3/2}}.\]
Then 
\begin{equation}
\label{eq:Abound}
A\gtrsim \max\big(\tau_1^{-1}, (m\tau_p)^{-1}\big).
\end{equation}
If $u_0>0$ exists, then by \eqref{eq:ufund}
\be
 |\psi'(u)| &\ge& u\int_{u_0}^u A\,s\,ds \gtrsim A\,u^2(u-u_0) \text{\ \ for all\ \ }\sqrt{\tau_m}>u>u_0\nonumber\\
 |\psi'(u)| &\ge& u\int_{u}^{u_0} A\,s\,ds \gtrsim A\,u_0\,u(u_0-u) \text{\ \ for all\ \ }0<u<u_0.\nonumber
\ee
By Lemma~\ref{lem:urechts} and~\ref{lem:ulinks} the oscillatory integral in~\eqref{eq:umain1} is therefore no larger than
\[ \min\Bigl(\tau_1^{\frac14}, \Bigl(\frac{m\tau_p}{\tau_k}\Bigr)^{\half}\Bigr) \lesssim 
    \min\Bigl(\tau_1^{\frac14}, m^{\frac32}\,\frac{b_p}{b_k} \Bigr),
\]
as claimed. If $u_0>0$ does not exist, then~\eqref{eq:ufundne} implies that
\[ |\psi'(u)| \gtrsim A\,u^3 \text{\ \ for all\ \ }0<u<\sqrt{\tau_m}.\]
Now apply \eqref{eq:uright1} of Lemma~\ref{lem:urechts} with $u_0=0$.
\end{proof}

\noindent The previous estimates allow us to deal with the singularity at $\lambda=\sqrt{\rho_\ell}$ in the 
following situation.

\begin{lemma} 
\label{lem:singlem1}
There exists a constant $C_0$ so that for any choice of $\sigma_1\ge \sigma_2\ge \ldots\ge \sigma_m \ge0$, 
$\rho_1\ge \rho_{2}\ge \ldots \ge \rho_{\ell}>0$,  $b_j>0$, and $c_i>0$, one has
\be
&& \left| 
\int_0^{\sqrt{\rho_\ell}} e^{\half i\lambda^2} \,e^{\pm i\sum_{j=1}^m b_j\,\sqrt{\lambda^2+\sigma_j}}
\,\exp\Bigl(-\sum_{i=1}^\ell c_i\sqrt{\rho_i-\lambda^2}\Bigr)\,\frac{\lambda}{\sqrt{\lambda^2+\sigma_k}}\,d\lambda \right|  \nonumber \\
&& \qquad\qquad \le C_0\min\Bigl[(1+\sigma_1+\rho_{\ell})^{\frac14}, m^{\frac32} b_k^{-1}\max_{1\le j\le m} 
b_j\Bigr]
\label{eq:bound2}
\ee
for any $1\le k\le m$.
\end{lemma}
\begin{proof}
For the purposes of this proof, set $w_k(\lambda)=\frac{\lambda}{\sqrt{\lambda^2+\sigma_k}}$. 
We claim that
\begin{equation}
\label{eq:claimder}
\Bigl| \frac{d^j}{d\lambda^j} \Bigl[ \exp\Bigl(-\sum_{i=1}^\ell c_i\sqrt{\rho_i-\lambda^2}\Bigr)\,w_k(\lambda)\Bigl(1-\chi\Bigl(\frac{\lambda}{\sqrt{\rho_\ell}}\Bigr)\Bigr) \Bigr] \Bigr| \le C_j\, \lambda^{-j}
\end{equation}
for every $j\ge0$ and $\lambda>0$. In view of Lemma~\ref{lem:lem1} it suffices to show that 
\[ 
\Bigl| \frac{d^j}{d\lambda^j} \exp\Bigl(- c_i\sqrt{\rho_i-\lambda^2}\Bigr) \Bigr| \le C_j\,\lambda^{-j}
\]
for every $j\ge0$, $1\le i\le \ell$,  and $0<\lambda<\half\sqrt{\rho_\ell}$. Fix some such~$i$ and note that
$\lambda\mapsto \sqrt{\rho_i-\lambda^2}$ is a diffeomorphism on the interval~$0<\lambda<\half\sqrt{\rho_\ell}$.
Since
\[ \Bigl| \frac{d^j}{d u^j} e^{-c u} \Bigr| = \Bigl| c^j e^{-c u} \Bigr| \le u^{-j} \sup_{v>0} (cv)^j e^{-cv} = C_j\,u^{-j}\]
irrespective of the choice of $c$, \eqref{eq:claimder} follows. Lemma~\ref{lem:osc1}, with $\omega$ equal to
the function in~\eqref{eq:claimder}, therefore implies that
\be
&& \left| 
\int_0^{\sqrt{\rho_\ell}} e^{\half i\lambda^2} \,e^{\pm i\sum_{j=1}^m b_j\,\sqrt{\lambda^2+\sigma_j}}
\,\exp\Bigl(-\sum_{i=1}^\ell c_i\sqrt{\rho_i-\lambda^2}\Bigr)\,w_k(\lambda)\Bigl(1-\chi\Bigl(\frac{\lambda}{\sqrt{\rho_\ell}}\Bigr)\Bigr)\,d\lambda
\right| \nonumber\\ 
&& \qquad\qquad \le C_0\,\min\Bigl[(1+\sigma_1)^{\frac14}, m^{\frac32}\,b_k^{-1}\,\max_{1\le j\le m} b_j\Bigr]. \label{eq:osc1reduc}
\ee
It remains to deal with the integral over the interval close to $\sqrt{\rho_\ell}$. Fix any $1\le k\le m$. The
change of variables $u=\sqrt{\rho_\ell-\lambda^2}$ leads to the identity
\be 
&& \int_0^{\sqrt{\rho_\ell}} e^{\half i\lambda^2} \,e^{\pm i\sum_{j=1}^m b_j\,\sqrt{\lambda^2+\sigma_j}}
\,\exp\Bigl(-\sum_{i=1}^\ell c_i\sqrt{\rho_i-\lambda^2}\Bigr)\,w_k(\lambda)\chi\Bigl(\frac{\lambda}{\sqrt{\rho_\ell}}\Bigr)\,d\lambda 
\nonumber\\
&=&  e^{i\rho_\ell/2}\int_0^\infty e^{-i\psi_{\mp}(u)} \; \omega(u)\,\frac{u}{\sqrt{\sigma_k+\rho_\ell}}\,du,
\label{eq:secondhalf}
\ee
where 
\begin{equation}
\label{eq:psipm}
\psi_{\mp}(u)=\half u^2\mp \sum_{j=1}^m b_j\sqrt{\sigma_j+\rho_\ell-u^2}
\end{equation}
and 
\[ 
\omega(u) = \exp\Bigl(-\sum_{i=1}^\ell c_i\sqrt{\rho_i-\rho_\ell+u^2}\Bigr) \chi\Bigl(\sqrt{1-u^2/\rho_\ell}\Bigr) 
\sqrt{\frac{\sigma_k+\rho_\ell}{\sigma_k+\rho_\ell-u^2}}.
\]
By the same arguments that lead to \eqref{eq:claimder} one sees that $|\omega^{(j)}(u)|\le C_j\,u^{-j}$ for every~$j\ge0$. 
Moreover, $\omega(u)=0$ for all $u\ge \sqrt{\frac{15}{16}}\sqrt{\rho_\ell}$. 
In case of the more difficult phase $\psi_+(u)$, Corollary~\ref{cor:uosc1} with $\tau_j:=\sigma_j+\rho_\ell$ therefore yields the
desired bound~\eqref{eq:bound2}. If, on the other hand, $\psi(u)=\psi_{-}(u)$, then 
\[
\psi'(u) = u + \sum_{j=1}^m \frac{b_ju}{\sqrt{\tau_j-u^2}} \qquad\qquad
\psi''(u) = 1 + \sum_{j=1}^m \frac{b_j\tau_j}{(\tau_j-u^2)^{\frac32}} 
\]
shows that $\psi'(u)\ge u$ and $\psi''(u)\ge0$ on the support of $\omega$. Thus integrating by parts once leads to
\be
 && \left| \int_0^\infty e^{-i\psi(u)} \; \omega(u)\,\frac{u}{\sqrt{\sigma_k+\rho_\ell}}\,du \right| \nonumber\\
 && \lesssim 1+ \int_0^\infty \frac{\psi''(u)}{\psi'(u)^2}\,\omega(u)\chi(u)\,\frac{u}{\sqrt{\sigma_k+\rho_\ell}}\,du +
 \int_0^\infty \frac{1}{\psi'(u)}\,\Bigl|\Bigl[\omega(u)\chi(u)\,\frac{u}{\sqrt{\sigma_k+\rho_\ell}}\Bigr]'\Bigr|\,du. 
\label{eq:uwieder}
\ee
Integrating by parts once more in the first integral, see~\eqref{eq:ubdaryterm}, reduces it to the second. Hence 
\[ \eqref{eq:uwieder} \lesssim \int_1^\infty \frac{1}{u\psi'(u)}\,du \lesssim 1,\]
and the lemma follows.
\end{proof}

\noindent To conclude this section, we turn to oscillatory integrals with singular weights.

\begin{lemma} 
\label{lem:singlem2}
There exists a constant $C_0$ so that for any choice of $\sigma_1\ge \sigma_2\ge \ldots\ge \sigma_m \ge0$, 
$\rho_1\ge \rho_{2}\ge \ldots \ge \rho_{\ell}>0$,  $b_j>0$, and $c_i>0$, one has
\be
&& \left| 
\int_0^{\sqrt{\rho_\ell}} e^{\half i\lambda^2} \,e^{\pm i\sum_{j=1}^m b_j\,\sqrt{\lambda^2+\sigma_j}}
\,\exp\Bigl(-\sum_{i=1}^\ell c_i\sqrt{\rho_i-\lambda^2}\Bigr)\,\frac{\lambda}{\sqrt{\rho_k-\lambda^2}}\,d\lambda \right|  \nonumber \\
&& \qquad\qquad \le C_0\min\Bigl[(1+\sigma_1+\rho_{\ell})^{\frac14}, m^{\frac32} c_k^{-1}\max_{j,i} (b_j+c_i)\Bigr]
\label{eq:bound3}
\ee
for any $1\le k\le m$.
\end{lemma}
\begin{proof} 
We start with the elementary comment that we can assume that
\begin{equation}
\label{eq:rhoell1}
\rho_\ell\ge 1.
\end{equation}
Indeed, if \eqref{eq:rhoell1} fails, then the oscillatory integral in \eqref{eq:bound3} is
\[ \le \int_0^{\sqrt{\rho_\ell}} \frac{\lambda}{\sqrt{\rho_k-\lambda^2}}\,d\lambda 
   \le \int_0^{\sqrt{\rho_\ell}} \frac{\lambda}{\sqrt{\rho_\ell-\lambda^2}}\,d\lambda = \sqrt{\rho_\ell}\le1.
\]
As in the previous proof, we first consider the case  $0\le \lambda\le \half\sqrt{\rho_\ell}$. 
Set $w_k(\lambda):=\frac{\lambda}{\sqrt{\rho_k-\lambda^2}}$ and fix some $1\le k\le \ell$. 
Then 
\[ \omega(\lambda):= \exp\Bigl(-\sum_{i=1}^\ell c_i\sqrt{\rho_i-\lambda^2}\Bigr)\,w_k(\lambda)\Bigl(1-\chi\Bigl(\frac{\lambda}{\sqrt{\rho_\ell}}\Bigr)\Bigr) \]
satsifies the derivative bounds
\begin{equation}
\label{eq:omder2} 
|\omega^{(j)}(\lambda)| \le C_j\,\lambda^{-j}
\end{equation}
for all $\lambda>0$ and $j\ge0$. Hence \eqref{eq:bound1'} of Lemma~\ref{lem:osc1} implies that
\begin{equation}
\label{eq:lampart} 
\left| 
\int_0^{\sqrt{\rho_\ell}} e^{\half i\lambda^2} \,e^{\pm i\sum_{j=1}^m b_j\,\sqrt{\lambda^2+\sigma_j}}
\,\omega(\lambda)\,d\lambda
\right| 
\lesssim (1+\sigma_1)^{\frac14}.
\end{equation}
For the part close to $\sqrt{\rho_\ell}$ we again use the change of variables $u=\sqrt{\rho_\ell-\lambda^2}$, which yields
\be 
&& \int_0^{\sqrt{\rho_\ell}} e^{\half i\lambda^2} \,e^{\pm i\sum_{j=1}^m b_j\,\sqrt{\lambda^2+\sigma_j}}
\,\exp\Bigl(-\sum_{i=1}^\ell c_i\sqrt{\rho_i-\lambda^2}\Bigr)\,w_k(\lambda)\chi\Bigl(\frac{\lambda}{\sqrt{\rho_\ell}}\Bigr)\,d\lambda 
\nonumber\\
&=&  e^{i\rho_\ell/2}\int_0^\infty e^{-i\psi_{\mp}(u)} \; \omega(u)\,du.
\label{eq:apjat}
\ee
Here $\psi_{\mp}(u)$  are as in \eqref{eq:psipm}, whereas
\begin{equation}
\label{eq:uomegadef} 
\omega(u)=\exp\Bigl(-\sum_{i=1}^\ell c_i\sqrt{\rho_i-\rho_\ell+u^2}\Bigr) \chi\Bigl(\sqrt{1-u^2/\rho_\ell}\Bigr)
\frac{u}{\sqrt{\rho_k-\rho_\ell+u^2}}.
\end{equation}
Observe that $\omega(u)=0$ if $u\ge\sqrt{\frac{15}{16}}\sqrt{\rho_\ell}$. Moreover, the reasoning that lead to~\eqref{eq:claimder}
yields in this case that
\[ \Bigl|\frac{d^j}{du^j}\, \omega(u)\Bigr| \le C_j\,u^{-j}\]
for all $j\ge0$ and $u>0$. 
Hence \eqref{eq:uright2} of Lemma~\ref{lem:urechts} and \eqref{eq:uleft2} of Lemma~\ref{lem:ulinks}
imply that 
\[
\left| \int_0^\infty e^{-i\psi_{+}(u)} \; \omega(u)\,du \right|
\lesssim (1+\sigma_1+\rho_\ell)^{\frac14},
\]
whereas the case of $\psi_{-}$ is the dealt with in the same way as in the previous proof, see~\eqref{eq:uwieder}. 
In conjunction with \eqref{eq:lampart} this yields the first bound in~\eqref{eq:bound3}. To obtain the second estimate of~\eqref{eq:bound3}, we first consider the case of $\lambda\le\half\sqrt{\rho_\ell}$, cf.~\eqref{eq:lampart}. In contrast to the latter formula, though, we seek bounds on
\begin{equation}
\label{eq:lampart2}
\left| 
\int_0^{\sqrt{\rho_\ell}} e^{i\phi_{\pm}(\lambda)}\;\omega(\lambda)\,\frac{\lambda}{\sqrt{\rho_k+\lambda^2}}\,d\lambda
\right| 
\end{equation}
with $\phi_{\pm}=\half\lambda^2\pm \sum_{j=1}^m b_j\,\sqrt{\lambda^2+\sigma_j}$, and 
\[ \omega(\lambda):= \exp\Bigl(-\sum_{i=1}^\ell c_i\sqrt{\rho_i-\lambda^2}\Bigr)\,\Bigl(1-\chi\Bigl(\frac{\lambda}{\sqrt{\rho_\ell}}\Bigr)\Bigr)\,\sqrt{\frac{\rho_k+\lambda^2}{\rho_k-\lambda^2}}. \]
In view of these choices \eqref{eq:lampart2} is equal to \eqref{eq:lampart}.
Since $\lambda\le\half\sqrt{\rho_\ell}\le\half\sqrt{\rho_k}$ on the support of $\omega$, this function satisfies the derivative bounds~\eqref{eq:omder2}. Suppose that a critical point $\lambda_0>0$ of $\phi_{-}$ exists. 
By definition of $\lambda_0$, see~\eqref{eq:lam0def}, there exists $1\le p\le m$ such that
\begin{equation}
\label{eq:pb2}
m\,b_p\ge \sqrt{\sigma_p+\lambda_0^2}.
\end{equation}
Hence, see \eqref{eq:zwei'},
\be
 \phi''(\lambda_0) &=& \sum_{j=1}^m \frac{b_j\lambda_0^2}{(\sigma_j+\lambda_0^2)^{\frac32}} \ge
\lambda_0^2 \frac{b_p}{\sqrt{\lambda_0^2+\sigma_p}} \, \frac{1}{\lambda_0^2+\sigma_p} \nonumber\\
&\ge& \frac{\lambda_0^2}{m(\lambda_0^2+\sigma_p)}. \label{eq:phi''5}
\ee
In view of \eqref{eq:left2'} one can apply Lemma~\ref{lem:left} with 
\[ B=\frac{\phi''(\lambda_0)}{\lambda_0}.\]
By \eqref{eq:phi''5} this leads to the bound
\begin{equation}
\label{eq:left5}
\left| 
\int_0^{\lambda_0} e^{i\phi_{-}(\lambda)}\;\omega(\lambda)\,\frac{\lambda}{\sqrt{\rho_k+\lambda^2}}\,d\lambda
\right| \lesssim \phi''(\lambda_0)^{-\half} \frac{\lambda_0}{\sqrt{\lambda_0^2+\rho_k}} \lesssim 
\sqrt{m}\sqrt{\frac{\sigma_p+\lambda_0^2}{\rho_k+\lambda_0^2}} \lesssim m^{\frac32} \frac{b_p}{\sqrt{\rho_k}}.
\end{equation}
We will now estimate the contribution by the interval $[\lambda_0,\infty)$ by means of Lemma~\ref{lem:right}
with $\tau=\rho_k$. The condition $\lambda_0\le\tau$ in that lemma can be assumed to hold, as otherwise~$\omega$
vanishes on the interval~$[\lambda_0,\infty)$. 
If the positive critical point $\lambda_0$ does not exist, then we formally set $\lambda_0=0$. Recall that
\[ 1\ge\sum_{j=1}^m \frac{b_j}{\sqrt{\sigma_j}}\]
if $\lambda_0>0$ does not exist, see \eqref{eq:lam0ne}. If 
\[ \half \ge\sum_{j=1}^m \frac{b_j}{\sqrt{\sigma_j}},\]
then $\phi'(\lambda)\ge\half\lambda$. Thus \eqref{eq:A} holds with $A=\half$, and therefore~\eqref{eq:rechts1}
leads to a~$O(1)$ bound. In view of \eqref{eq:lam0def} and the preceding we can thus assume that 
\[ 1 \asymp \sum_{j=1}^m \frac{b_j}{\sqrt{\sigma_j+\lambda_0^2}}\]
regardless of the value of $\lambda_0$. 
In particular, there exists some $1\le p\le m$ with the property that $m\,b_p\gtrsim\sqrt{\sigma_p+\lambda_0^2}$. 
We now proceed to verify the conditions~\eqref{eq:R} with~$R\asymp \sigma_p+\lambda_0^2$. More precisely, let $\lambda_0\le\lambda\le 2\sqrt{\sigma_p+\lambda_0^2}$. Then~\eqref{eq:fund} or~\eqref{eq:fundne} imply that
\be 
 \phi'(\lambda) &\gtrsim& \lambda \int_{\lambda_0}^\lambda \frac{b_p\,s}{(b_p+s^2)^{\frac32}}\,ds 
\gtrsim \lambda^2 \frac{b_p}{(b_p+\lambda_0^2)^{\frac32}}\,(\lambda-\lambda_0) \nonumber\\
&\gtrsim& \frac{\lambda^2}{m(\sigma_p+\lambda_0^2)}\,(\lambda-\lambda_0), \nonumber
\ee
which is the first condition in~\eqref{eq:R} with $R\asymp \sigma_p+\lambda_0^2$ and $c=m^{-1}$. If on the
other hand $\lambda \ge 2\sqrt{\sigma_p+\lambda_0^2}$. Then 
\be 
\phi'(\lambda) &\gtrsim& \lambda \int_{\sqrt{\sigma_p+\lambda_0^2}}^\lambda \frac{b_p\,s}{(b_p+s^2)^{\frac32}}\,ds 
\gtrsim \lambda \int_{\sqrt{\sigma_p+\lambda_0^2}}^\lambda \frac{b_p}{s^2}\,ds \nonumber \\
&\gtrsim& \frac{b_p}{\sqrt{\sigma_p+\lambda_0^2}}\,(\lambda-\sqrt{\sigma_p+\lambda_0^2}) 
\gtrsim m^{-1}\,(\lambda-\lambda_0). \label{eq:weitrechts}
\ee
To pass to the last line we used that $b_p\gtrsim m^{-1}\sqrt{\sigma_p+\lambda_0^2}$ as well as 
\[ 2(\lambda-\sqrt{\sigma_p+\lambda_0^2}) \ge \lambda-\lambda_0.\] 
Since \eqref{eq:weitrechts} verifies the second condition in~\eqref{eq:R} we conclude from~\eqref{eq:rechts2} 
 that 
\begin{equation}
\nonumber
\left| 
\int_{\lambda_0}^\infty e^{i\phi_{-}(\lambda)}\;\omega(\lambda)\,\frac{\lambda}{\sqrt{\rho_k+\lambda^2}}\,d\lambda
\right| \lesssim m^{\frac32}+ \sqrt{m}\sqrt{\frac{\sigma_p+\lambda_0^2}{\rho_k}} \lesssim m^{\frac32}+
m^{\frac32} \frac{b_p}{\sqrt{\rho_k}}.
\end{equation}
In conjunction with \eqref{eq:left5} this shows that 
\begin{equation}
\label{eq:great}
\left| 
\int_0^{\sqrt{\rho_\ell}} e^{i\phi_{-}(\lambda)}\;\omega(\lambda)\,\frac{\lambda}{\sqrt{\rho_k+\lambda^2}}\,d\lambda
\right|  \lesssim m^{\frac32} + m^{\frac32} \frac{\max_j b_j}{\sqrt{\rho_k}}.
\end{equation}
Since the phase $\phi_{+}$ satisfies $\phi_{+}'(\lambda)\ge \lambda$, the integral in~\eqref{eq:lampart2}
with phase $\phi_{+}$ is~$O(1)$ by the usual argument involving a single integration by parts. Thus~\eqref{eq:lampart2} satisfies the estimate~\eqref{eq:great} in all cases.
Clearly, both~\eqref{eq:lampart2} and the integral in~\eqref{eq:bound3} are no larger than
\begin{equation}
\label{eq:trivci}
\min_{1\le i\le \ell}\int_0^{\sqrt{\rho_\ell}} e^{-c_i\sqrt{\rho_i-\lambda^2}}\; \frac{\lambda}{\sqrt{\rho_k-\lambda^2}}\,d\lambda \le \min_{1\le i\le \ell}\int_0^{\sqrt{\rho_\ell}} e^{-c_i\sqrt{\rho_\ell-\lambda^2}}\; \frac{\lambda}{\sqrt{\rho_\ell-\lambda^2}}\,d\lambda = \min_{1\le i\le \ell} c_i^{-1}.
\end{equation} 
The desired bound \eqref{eq:bound3} now follows easily from \eqref{eq:great} and \eqref{eq:trivci}. Indeed,
if~$c_k\le \sqrt{\rho_k}$, then~\eqref{eq:great} implies~\eqref{eq:bound3}. On the other hand, if~$c_k\ge \sqrt{\rho_k}$, then in view of~\eqref{eq:rhoell1} one has $c_k\ge \sqrt{\rho_k}\ge\sqrt{\rho_\ell}\ge 1$. In that 
case~\eqref{eq:trivci} yields a~$O(1)$ bound.

\noindent It remains to consider the case of $\lambda$ close to $\sqrt{\rho_\ell}$. In that case the change of
variables $u=\sqrt{\rho_k-\lambda^2}$ reduces matters to~\eqref{eq:apjat} with~$\omega(u)$ as in~\eqref{eq:uomegadef}. 
Suppose first that 
\[\sum_{j=1}^m \frac{b_j}{\sqrt{\sigma_j+\rho_\ell}} \le \half.\]
Then Lemma~\ref{lem:uphase} implies that $\psi_{+}'(u)\ge\half u$ on the support of~$\omega$, 
whereas in all cases $\psi_{-}'(u)\ge u$. In case of~$\psi_{+}(u)$ a~$O(1)$ bound is obtained via~\eqref{eq:O1uright}, whereas
for~$\psi_{-}(u)$ this bound follows by a single integration by parts in the usual way. Hence it suffices to assume that
\[\sum_{j=1}^m \frac{b_j}{\sqrt{\sigma_j+\rho_\ell}} \ge \half.\]
There exists some $1\le p\le m$ so that $m\,b_p \ge\sqrt{\sigma_p+\rho_\ell}$. By \eqref{eq:rhoell1} this implies that~$b_p\ge m^{-1}$. 
Since the estimate~\eqref{eq:trivci} also applies to the oscillatory integrals after changing variables to~$u$, one concludes from
the preceding that there is an upper bound of the form 
\[ m\;\frac{\max_{1\le j\le m} b_j}{c_k},\]
and we are done.
\end{proof}

\section{Putting it all together} \label{sec:time4}

By combining the results of the previous three sections we are now able to prove our main result.

\begin{theorem}
Let $V(t,x)$ be a real-valued measurable function on~$\R^4$ such that 
\begin{equation}
\sup_{t}\|V(t,\cdot)\|_{L^{\frac32}(\R^3)} < c_0 \text{\ \ and\ \ } \sup_{y\in\R^3} \int_{\R^3}\frac{\|V(\hat{\tau},x)\|_{\mathcal M}}{|x-y|}\,dx < 4\pi
\end{equation}
for some small constant $c_0>0$, see Definition~\ref{def:Ydef}. Then
\begin{equation}
\nn 
\|U(t,s)\psi_s\|_{\infty} \le C|t-s|^{-\frac32}\,\|\psi_s\|_1 \text{\ \ for all times $t,s$ and any $\psi_s\in L^1$} ,
\end{equation}
where $U(t,s)$ is the weak propagator constructed in Lemma~\ref{lem:weaksol}. 
\end{theorem}
\begin{proof}
Recall from Proposition~\ref{Fcalc} that 
\[ \langle U(t,s)\psi_s, g\rangle = \sum_{m=0}^\infty \langle {\mathcal I}_m \psi_s,g \rangle \]
for any pair $\psi_s,g\in \calS(\R^3)$. Furthermore, Proposition~\ref{KernRep} provides a
representation of the kernel of $ {\mathcal I}_m(t,s)$ in terms of three kinds of oscillatory integrals,
which are defined in~\eqref{finLm}, \eqref{Mdl}, and~\eqref{Mbardl}. Suppose $t>s$. Changing variables
$\lambda\mapsto \frac{\lambda}{\sqrt{t-s}}$ in each of these integrals brings out one factor of~$(t-s)^{-\frac12}$, 
whereas~\eqref{eq:kernrep} already contains the factor~$(t-s)^{-1}$. This leads to the desired power~$(t-s)^{-\frac32}$. More precisely, for the oscillatory integrals from~\eqref{finLm} this process leads to
\be
&& \int_{0}^\infty  e^{i(t-s)\la^2} \cos{\bigg(\sum_{k=1}^{m}\sqrt{\la^2+\si_k}|x_{k-1}-x_k|\bigg)} \frac {\la}{\sqrt{\la^2+\si_{\ell}}} \, d\la \nn \\
&=& (t-s)^{-\frac12}\,\int_{0}^\infty  e^{i\la^2} \cos{\bigg(\sum_{k=1}^{m}\sqrt{\la^2+\si_k(t-s)}\frac{|x_{k-1}-x_k|}
{\sqrt{t-s}}\bigg)} \frac {\la}{\sqrt{\la^2+\si_{\ell}(t-s)}} \, d\la, \label{eq:rescLm} 
\ee
and similarly for \eqref{Mdl} and~\eqref{Mbardl}. Thus the parameters $\sigma_j$ and $\rho_k$ and $|x_{i+1}-x_i|$ in these expressions are rescaled to $\sigma_j\,(t-s)$, $\rho_k\,(t-s)$, and $\frac{|x_{i+1}-x_i|}{\sqrt{t-s}}$, respectively. We now estimate~\eqref{eq:rescLm} and the analogous integrals from~\eqref{Mdl}, and~\eqref{Mbardl} by means of Lemma~\ref{lem:osc1}, \ref{lem:singlem1}, \ref{lem:singlem2}, respectively. Using the second bound in each of these lemmas, which is invariant under the aforementioned rescaling of the parameters, one arrives at the upper bound (setting $x=x_0$ and $y=x_{m+1}$) 
\[ m^{\frac32}\; \frac{\max\limits_{0\le j\le m+1} |x_{j+1}-x_j|}{|x_\ell-x_{\ell-1}|}\]
in case of ${\mathcal L}_m^\ell(t,s)(x,y)$, and 
\[ m^{\frac32}\; \frac{\max\limits_{0\le j\le m+1} |x_{j+1}-x_j|}{|x_{k(\ell)}-x_{k(\ell)-1}|}\]
in case of ${\mathcal M}_m^{d,\ell}(t,s)(x,y),\; \widetilde{\mathcal  M}_m^{d,\ell}(t,s)(x,y)$. 
Inserting these bounds into the (rescaled) definitions~\eqref{finLm}, \eqref{Mdl}, and~\eqref{Mbardl} 
finally leads to the estimate
\be
&& |\langle {\mathcal I}_m \psi_s,g \rangle| \nn\\
&\le& C\, m^{\frac72}\;|t-s|^{-\frac32}\, \int_{\R^m}\int_{{\mathbb R}^{m}}
\prod_{r=1}^{m} \frac{|V(\hat\tau_{r},x_r)|}{4\pi|x_{r-1}-x_r|}\,\,
\frac{\max\limits_{0\le j\le m+1} |x_{j+1}-x_j|}{|x_m-y|}\;dx_1\ldots dx_{m}\,d\tau_1\ldots d\tau_m.
\nn
\ee
In view of Lemma~\ref{lem:iter} this is no larger than  
\[ C\,m^{\frac92}\;|t-s|^{-\frac32}\,  \, \left(\sup_{y\in\R^3}\int_{\R^3} \int\frac{|V(\hat\tau,x)|}{4\pi|x-y|}\,d\tau\,dx\right)^{m},\]
and we are done.
\end{proof}

\bibliographystyle{amsplain}

\noindent
\textsc{Rodnianski: Department of Mathematics, Princeton University, Fine Hall, Princeton N.J. 08544, U.S.A.}\\
{\em email: }\textsf{\bf irod@math.princeton.edu}

\medskip\noindent
\textsc{Schlag: Division of Astronomy, Mathematics, and Physics, 253-37 Caltech, Pasadena, CA 91125, U.S.A.}\\
{\em current address: } \textsc{Institute for Advanced Study, Olden Lane, Princeton, NJ 08540}\\
{\em email: }\textsf{\bf schlag@ias.edu}

\end{document}